\documentclass[letter,11pt,reqno]{article}
\usepackage[margin = 1in]{geometry}
\usepackage{amsmath,amssymb,amsthm,bm,enumitem,mathrsfs,mathtools,mleftright,setspace,soul,tensor,tikz-cd}

\usepackage{titlesec}

\allowdisplaybreaks
\setlist[itemize]{leftmargin = *}
\setlist[enumerate]{leftmargin = *}

\usepackage[colorlinks]{hyperref}
\usepackage{amsrefs}

\makeatletter
\DeclareFontFamily{OMX}{MnSymbolE}{}
\DeclareSymbolFont{MnLargeSymbols}{OMX}{MnSymbolE}{m}{n}
\SetSymbolFont{MnLargeSymbols}{bold}{OMX}{MnSymbolE}{b}{n}
\DeclareFontShape{OMX}{MnSymbolE}{m}{n}{
    <-6>  MnSymbolE5
   <6-7>  MnSymbolE6
   <7-8>  MnSymbolE7
   <8-9>  MnSymbolE8
   <9-10> MnSymbolE9
  <10-12> MnSymbolE10
  <12->   MnSymbolE12
}{}
\DeclareFontShape{OMX}{MnSymbolE}{b}{n}{
    <-6>  MnSymbolE-Bold5
   <6-7>  MnSymbolE-Bold6
   <7-8>  MnSymbolE-Bold7
   <8-9>  MnSymbolE-Bold8
   <9-10> MnSymbolE-Bold9
  <10-12> MnSymbolE-Bold10
  <12->   MnSymbolE-Bold12
}{}

\let\llangle\@undefined
\let\rrangle\@undefined
\DeclareMathDelimiter{\llangle}{\mathopen}{MnLargeSymbols}{'164}{MnLargeSymbols}{'164}
\DeclareMathDelimiter{\rrangle}{\mathclose}{MnLargeSymbols}{'171}{MnLargeSymbols}{'171}

\DeclareFontFamily{U}{matha}{\hyphenchar\font45}
\DeclareFontShape{U}{matha}{m}{n}{<5> <6> <7> <8> <9> <10> gen * matha <10.95> matha10 <12> <14.4> <17.28> <20.74> <24.88> matha12}{}
\DeclareSymbolFont{matha}{U}{matha}{m}{n}
\DeclareFontSubstitution{U}{matha}{m}{n}

\DeclareFontFamily{U}{mathx}{\hyphenchar\font45}
\DeclareFontShape{U}{mathx}{m}{n}{<5> <6> <7> <8> <9> <10> <10.95> <12> <14.4> <17.28> <20.74> <24.88> mathx10}{}
\DeclareSymbolFont{mathx}{U}{mathx}{m}{n}
\DeclareFontSubstitution{U}{mathx}{m}{n}

\DeclareMathDelimiter{\vvvert}{0}{matha}{"7E}{mathx}{"17}

\newtheorem{Lem}{Lemma}
\newtheorem{Thm}{Theorem}
\newtheorem{ThmDef}{Theorem-Definition}
\newtheorem{Prop}{Proposition}
\newtheorem{Cor}{Corollary}

\theoremstyle{definition}
\newtheorem{Def}{Definition}
\newtheorem{Eg}{Example}
\newtheorem{Rmk}{Remark}

\makeatletter
\let\orgautoref\autoref
\providecommand{\Identity}[1]{\def\equationautorefname{Identity}\orgautoref{#1}}
\providecommand{\Inequality}[1]{\def\equationautorefname{Inequality}\orgautoref{#1}}
\makeatother

\BibSpec{collection.article}{%
    +{}  {\PrintAuthors}                {author}
    +{,} { \textit}                     {title}
    +{.} { }                            {part}
    +{:} { \textit}                     {subtitle}
    +{,} { \PrintContributions}         {contribution}
    +{,} { \PrintConference}            {conference}
    +{}  {\PrintBook}                   {book}
    +{,} { }                            {booktitle}
    +{,} { }                            {publisher}
    +{,} { }                            {address}
    +{,} { \PrintDateB}                 {date}
    +{,} { pp.~}                        {pages}
    +{,} { }                            {status}
    +{,} { \PrintDOI}                   {doi}
    +{,} { available at \eprint}        {eprint}
    +{}  { \parenthesize}               {language}
    +{}  { \PrintTranslation}           {translation}
    +{;} { \PrintReprint}               {reprint}
    +{.} { }                            {note}
    +{.} {}                             {transition}
    +{}  {\SentenceSpace \PrintReviews} {review}
}

\newcommand{\A}{\mathscr{A}}
\newcommand{\B}[1]{\func{\mathscr{B}}{#1}}
\newcommand{\C}[1]{\mathbb{C}^{#1}}
\newcommand{\D}{\displaystyle}
\newcommand{\E}{\mathcal{E}}
\newcommand{\F}{\mathcal{F}}
\newcommand{\M}{\mathcal{M}}
\newcommand{\N}[1]{\mathbb{N}^{#1}}
\newcommand{\R}{\mathcal{R}}
\newcommand{\T}[1]{\mathbb{T}^{#1}}
\newcommand{\U}[1]{\func{\mathcal{U}}{#1}}
\newcommand{\X}{\mathsf{X}}
\newcommand{\Y}{\mathsf{Y}}
\newcommand{\Z}[1]{\mathbb{Z}^{#1}}
\newcommand{\Ad}[1]{\func{\operatorname{Ad}}{#1}}
\newcommand{\Br}[1]{\( #1 \)}
\newcommand{\Cb}[1]{\func{C_{b}}{#1}}
\newcommand{\Cc}[1]{\func{C_{c}}{#1}}
\newcommand{\CC}[2]{\SqBr{#1,#2}}
\newcommand{\Cl}[2]{\overline{#1}^{#2}}

\newcommand{\CP}[3]{\func{C^{*}}{#1,#2,#3}}
\newcommand{\df}{\vcentcolon =}
\newcommand{\FF}{\mathscr{F}}
\newcommand{\GG}{\mathscr{G}}
\newcommand{\Id}{\operatorname{Id}}
\newcommand{\IM}[2]{\Br{#1} \SqBr{#2}}
\newcommand{\LL}[2]{\func{L^{#1}}{#2}}

\newcommand{\om}[2]{\func{\omega}{#1,#2}}
\newcommand{\RM}[3]{\bm{\R}_{\Pair{#1}{#2}}^{#3}}
\newcommand{\RR}[2]{\mathbb{R}_{#1}^{#2}}
\newcommand{\si}{\operatorname{si}}
\newcommand{\tr}{\operatorname{tr}}
\newcommand{\Abs}[1]{\mleft| #1 \mright|}
\newcommand{\Adj}[1]{\func{\mathbb{L}}{#1}}
\newcommand{\Aut}[1]{\func{\operatorname{Aut}}{#1}}
\newcommand{\bra}[1]{\Langle #1 \mright|}
\newcommand{\Bra}[1]{\Llangle #1 \mright|}
\newcommand{\Del}[2]{\func{\Delta}{#1}^{#2}}

\newcommand{\Esi}{\E_{\si}}
\newcommand{\FCP}[3]{\func{C^{*}}{#1,#2,#3}}
\newcommand{\Fix}[2]{\func{\mathsf{Fix}}{#1;#2}}
\newcommand{\Fsi}{\F_{\si}}
\newcommand{\gam}[1]{\gamma^{#1}}
\newcommand{\Imp}[2]{\func{\bm{\F}}{#1;#2}}
\newcommand{\Int}[3]{\int_{#1} #2 ~ \d{#3}}
\newcommand{\Inv}[1]{#1^{-1}}
\newcommand{\INV}[1]{\Br{#1}^{-1}}
\newcommand{\ket}[1]{\mleft| #1 \Rangle}
\newcommand{\Ket}[1]{\mleft| #1 \Rrangle}
\newcommand{\Map}[4]{\SSet{\begin{matrix} #1 & \to & #2 \\ #3 & \mapsto & \D #4 \end{matrix}}}
\newcommand{\Orb}[3]{#1 \backslash_{#3} #2}
\newcommand{\RCP}[3]{\func{C^{*}_{\operatorname{r}}}{#1,#2,#3}}
\newcommand{\Seq}[2]{\Br{#1}_{#2}}
\newcommand{\Set}[2]{\SSet{#1 ~ \middle| ~ #2}}
\newcommand{\Ball}[1]{\mathbb{B}_{#1}}
\newcommand{\BBra}[1]{\tensor[_{2}]{\Bra{#1}}{}}
\newcommand{\Comp}[1]{\func{\mathbb{K}}{#1}}
\newcommand{\Conj}[1]{\overline{#1}}
\newcommand{\Cont}[2]{\func{C_{#1}}{#2}}
\newcommand{\FixR}[4]{\func{\texttt{Fix}_{\Trip{#1}{#2}{#3}}}{#4}}
\newcommand{\FixT}[2]{\func{\mathfrak{Fix}}{#1;#2}}
\newcommand{\Flat}[2]{\func{#1^{\flat}}{#2}}

\newcommand{\FTCP}[4]{\func{C^{*}}{#1,#2,#3,#4}}
\newcommand{\func}[2]{#1 \Br{#2}}
\newcommand{\Func}[2]{\func{\Br{#1}}{#2}}
\newcommand{\FUNC}[2]{\func{\SqBr{#1}}{#2}}
\newcommand{\gamm}[2]{\gamma^{#1}_{#2}}
\newcommand{\Hilb}[4]{\func{\mathbf{Hilb}}{#1,#2,#3,#4}}
\newcommand{\ImpR}[4]{\func{\X_{\Trip{#1}{#2}{#3}}}{#4}}
\newcommand{\Isom}[1]{\func{\operatorname{Isom}}{#1}}
\newcommand{\KKet}[1]{\tensor[]{\Ket{#1}}{_{2}}}
\newcommand{\Mult}[1]{\func{M}{#1}}
\newcommand{\Norm}[2]{\mleft\| #1 \mright\|_{#2}}
\newcommand{\Odot}[1]{\odot_{#1}}
\newcommand{\Pair}[2]{\Br{#1,#2}}
\newcommand{\Quad}[4]{\Br{#1,#2,#3,#4}}
\newcommand{\RTCP}[4]{\func{C^{*}_{\operatorname{r}}}{#1,#2,#3,#4}}
\newcommand{\Span}[1]{\func{\operatorname{Span}}{#1}}
\newcommand{\SqBr}[1]{\[ #1 \]}
\newcommand{\SSet}[1]{\mleft\{ #1 \mright\}}
\newcommand{\Supp}[1]{\func{\operatorname{Supp}}{#1}}
\newcommand{\Trip}[3]{\Br{#1,#2,#3}}
\newcommand{\AdjEq}[1]{\func{\mathbb{L}_{\operatorname{eq}}}{#1}}
\newcommand{\Graph}[1]{\func{\operatorname{Graph}}{#1}}
\newcommand{\Inner}[3]{\Langle #1 \middle| #2 \Rangle_{#3}}
\newcommand{\NNorm}[1]{\mleft\vvvert #1 \mright\vvvert}
\newcommand{\Range}[1]{\func{\operatorname{Range}}{#1}}
\newcommand{\Sharp}[2]{\func{#1^{\sharp}}{#2}}
\newcommand{\UMult}[1]{\func{\mathcal{U}}{\Mult{#1}}}
\newcommand{\BraArg}[2]{\func{\Bra{#1}}{#2}}
\newcommand{\KetArg}[2]{\func{\Ket{#1}}{#2}}
\newcommand{\Langle}{\mleft\langle}
\newcommand{\Otimes}[1]{\otimes_{#1}}
\newcommand{\Rangle}{\mright\rangle}
\newcommand{\AdjPair}[2]{\mathbb{L} \Br{#1,#2}}
\newcommand{\alphArg}[2]{\func{\alpha_{#1}}{#2}}
\newcommand{\alphExt}[1]{\overline{\alph{#1}}}
\newcommand{\BigNorm}[2]{\Big\| #1 \Big\|_{#2}}
\newcommand{\BBraArg}[2]{\func{\BBra{#1}}{#2}}
\newcommand{\csiHilb}[4]{\func{\mathbf{c.s.i.Hilb}}{#1,#2,#3,#4}}
\newcommand{\FixFull}[5]{\func{\mathsf{Fix}_{\Trip{#1}{#2}{#3}}}{#4;#5}}
\newcommand{\gammArg}[3]{\func{\gamma^{#1}_{#2}}{#3}}
\newcommand{\GammArg}[2]{\func{\Gamma_{#1}}{#2}}
\newcommand{\HilbMod}[1]{\func{\mathbf{HilbMod}}{#1}}
\newcommand{\HModule}[3]{{_{#1}} #2_{#3}}
\newcommand{\ImpFull}[5]{\func{\bm{\F}_{\Trip{#1}{#2}{#3}}}{#4;#5}}
\newcommand{\KKetArg}[2]{\func{\KKet{#1}}{#2}}
\newcommand{\Llangle}{\mleft\llangle}
\newcommand{\Rrangle}{\mright\rrangle}
\newcommand{\BBraKKet}[2]{\tensor[_{2}]{\Llangle #1 \middle| #2 \Rrangle}{_{2}}}
\newcommand{\BigInner}[3]{\Big\langle #1 \Big| #2 \Big\rangle_{#3}}

\newcommand{\AdjEqPair}[2]{\func{\mathbb{L}_{\operatorname{eq}}}{#1,#2}}
\newcommand{\LeftInner}[3]{\tensor[_{#3}]{\Langle #1 \middle| #2 \Rangle}{}}
\newcommand{\alphArgExt}[2]{\func{\overline{\alpha_{#1}}}{#2}}

\newcommand{\numberthis}{\addtocounter{equation}{1}\tag{\theequation}}

\renewcommand{\(}{\mleft(}
\renewcommand{\)}{\mright)}
\renewcommand{\[}{\mleft[}
\renewcommand{\]}{\mright]}
\renewcommand{\d}[1]{\mathrm{d}{#1}}
\renewcommand{\L}{\mathcal{L}}
\renewcommand{\S}{\mathcal{S}}
\renewcommand{\Im}[2]{#1 \SqBr{#2}}
\renewcommand{\MR}[1]{~\href{http://www.ams.org/mathscinet-getitem?mr=#1}{MR #1}}
\renewcommand{\alph}[1]{\alpha_{#1}}

\begin{document}


\begin{titlepage}

{\setstretch{1.0}
\begin{center}
{\Large \textbf{GENERALIZED FIXED-POINT ALGEBRAS FOR \\ TWISTED $ C^{*} $-DYNAMICAL SYSTEMS}}

\mbox{}

BY

\mbox{}

LEONARD TRISTAN HUANG ZHILIANG

\vspace{1in}

Submitted to the graduate degree program in the Department of Mathematics and \\
the Graduate Faculty of the University of Kansas in partial fulfillment of \\
the requirements for the degree of Doctor of Philosophy.
\end{center}

\vspace{1in}

\hfill
\parbox{2in}{
\underline{\hspace{2in}} \\
Chairperson, Albert Sheu \\ [0.5in]
\underline{\hspace{2in}} \\
David Lerner             \\ [0.5in]
\underline{\hspace{2in}} \\
Rodolfo Torres           \\ [0.5in]
\underline{\hspace{2in}} \\
Jeremy Martin            \\ [0.5in]
\underline{\hspace{2in}} \\
John Peter Ralston
}

\vspace{1in}

\begin{center}
Date Defended: May 3, 2016
\end{center}
}

\end{titlepage}


\pagenumbering{roman}
\setcounter{page}{2}

\newpage

{\setstretch{1.0}
\begin{center}
The Dissertation Committee for Leonard Tristan Huang Zhiliang \\
certifies that this is the approved version of the following dissertation:

\vspace{1in}

GENERALIZED FIXED-POINT ALGEBRAS FOR TWISTED $ C^{*} $-DYNAMICAL SYSTEMS
\end{center}

\vspace{1in}

\hfill
\parbox{2in}{
\underline{\hspace{2in}} \\
Chairperson, Albert Sheu
}

\vspace{1in}

\begin{center}
Date approved: May 3, 2016
\end{center}
}


\newpage

\begin{center}
{\LARGE \textbf{Abstract}}
\end{center}

In his seminal paper \emph{Generalized Fixed Point Algebras and Square-Integrable Group Actions} \cite{Meyer2}, Ralf Meyer showed how to construct generalized fixed-point algebras for $ C^{*} $-dynamical systems via their square-integrable representations on Hilbert $ C^{*} $-modules. His method extends Marc Rieffel's construction of generalized fixed-point algebras from proper group actions in \cite{Rieffel2}.

This dissertation seeks to generalize Meyer's work to construct generalized fixed-point algebras for twisted $ C^{*} $-dynamical systems. To accomplish this, we must introduce some brand-new concepts, the foremost being that of a twisted Hilbert $ C^{*} $-module. A twisted Hilbert $ C^{*} $-module is basically a Hilbert $ C^{*} $-module equipped with a twisted group action that is compatible with the module's right $ C^{*} $-algebra action and its $ C^{*} $-algebra-valued inner product. Twisted Hilbert $ C^{*} $-modules form a category, where morphisms are twisted-equivariant adjointable operators, and we will establish that Meyer's bra-ket operators are morphisms between certain objects in this category.

A by-product of our work is a twisted-equivariant version of Kasparov's Stabilization Theorem, which states that every countably generated twisted Hilbert $ C^{*} $-module is isomorphic to an invariant orthogonal summand of the countable direct sum of a standard one if and only if the module is square-integrable.

Given a twisted $ C^{*} $-dynamical system, we provide a definition of a relatively continuous subspace of a twisted Hilbert $ C^{*} $-module (inspired by Ruy Exel's paper \cite{Exel}) and then prescribe a new method of constructing generalized fixed-point algebras that are Morita-Rieffel equivalent to an ideal of the corresponding reduced twisted crossed product. Our construction generalizes that of Meyer and, by extension, that of Rieffel in \cite{Rieffel2}.

Our main result is the description of a classifying category for the class of all Hilbert modules over a reduced twisted crossed product. This implies that every Hilbert module over a $ d $-dimensional non-commutative torus can be constructed from a Hilbert space endowed with a twisted $ \Z{d} $-action and a relatively continuous subspace.

\mbox{}



\newpage

\begin{center}
{\LARGE \textbf{Acknowledgments}}
\end{center}

Although this work was conceived in 2014 and started taking shape only within the past year, I began laying the foundation six years ago. Along the way, I received help and invaluable advice from several people, without whom I would not have finished my dissertation in a timely fashion.

I would like to thank Professor Ralf Meyer (Mathematisches Institut, Georg-August-Universit\"at G\"ottingen) for diverting precious time from his busy schedule to explain some of the subtle points in his seminal paper, \emph{Generalized Fixed Point Algebras and Square-Integrable Group Actions}. Although I have never met him in person, our email correspondences have been highly instrumental in helping me to develop a deep level of understanding of his paper, which forms the basis of my dissertation. I am therefore indebted to him for his generosity with his time.

I would also like to offer my thanks to Professor Judith Packer (University of Colorado Boulder) for her hospitality when I visited her in 2015 to seek her expertise on twisted $ C^{*} $-dynamical systems. Her research with Iain Raeburn on these objects has revolutionized the manner in which they are studied today, and it is to their work that I owe the discovery of many of my results herein. I am also extremely grateful for her very helpful suggestion of various references on the subject.

I wish to express my gratitude toward many of the professors in the Department of Mathematics: Professors Judith Roitman and Bill Fleissner, who reaffirmed my love for set theory; Professor Jack Porter, whose crystalline lectures on topology made it such a joy to learn; Professor David Lerner, whose courses on differential geometry and general relativity inspired me to adopt mathematical physics as a secondary research interest; Professor Rodolfo Torres, who taught me the deeper aspects of Fourier analysis; Professor Bozenna Pasik-Duncan, who shared with me her passion for applying mathematics to solve real-world problems; Professor Bangere Purnaprajna, who exposed me to the realm of algebraic geometry in a highly entertaining way; Professor Jeremy Martin, whose infusion of combinatorics into his algebraic topology class left a deep impression on me; and Professor Estela Gavosto, who freely offered encouragement and advice although I never took a class from her.

I owe my knowledge of functional analysis, $ C^{*} $-algebras and non-commutative geometry to the training of my indefatigable adviser, Professor Albert Sheu. His incisive mind, attention to detail, concise style and overarching view of mathematics have inspired me to be both a better researcher and a better expositor. I fondly remember my bi-weekly meetings with him in the summer of 2012, when he helped me to prepare for my preliminary exam. It was a very intense period, but I certainly emerged from it better equipped to undertake research.

Professor Sheu's office has become almost a second home for me over the years, where not only have we had discussions on mathematics, but also on family and life. The friendship that we have forged is one of the fondest memories that I will have of KU.

Lastly, I would like to dedicate my work to Magdalene. In the past ten years, her encouragement and support for me have been unwavering. There were many times when my frustrations in research could have gotten the better of me, but she was constantly there to bring the joy of mathematical discovery back to my heart. Her faith in my abilities is marrow-deep, and words alone do not suffice to describe my appreciation for her presence in my life.

Mathematical research is truly tough business. It is filled with long stretches of blissful learning, punctuated by moments of sheer terror when you see your proof turn to poof.


\tableofcontents


\section{Preliminaries} 

\pagenumbering{arabic}

In this section, we give a brief introduction to $ C^{*} $-algebras and some of their important properties. $ C^{*} $-algebras possess tightly intertwining algebraic and analytic structures, which give these algebras a remarkably wide range of applicability. They are thus important to other areas of mathematics, such as operator theory, harmonic analysis, algebraic topology and non-commutative geometry.

Some of the material here is taken from \cite{Raeburn|Williams,Wegge-Olsen,Williams}, which are standard references on the subject.

\subsection{$ C^{*} $-Algebras} 

A \emph{$ C^{*} $-algebra} is a complex Banach algebra $ A $ with an involution $ ^{*} $ that satisfies the \emph{$ C^{*} $-identity}:
$$
\forall a \in A: \quad
\Norm{a^{*} a}{A} = \Norm{a}{A}^{2}.
$$
It follows from the $ C^{*} $-identity that $ ^{*} $ is isometric, i.e., $ \Norm{a^{*}}{A} = \Norm{a}{A} $ for every $ a \in A $.

For a $ C^{*} $-algebra $ A $, we have the following standard terminology:
\begin{itemize}
\item
$ A $ is called \emph{unital} if and only if it has an identity element, i.e., an element $ 1_{A} $ such that
$$
\forall a \in A: \quad
1_{A} a = a 1_{A} = a.
$$

\item
$ a \in A $ is called \emph{self-adjoint} if and only if $ a^{*} = a $.

\item
$ a \in A $ is called \emph{normal} if and only if $ a^{*} a = a a^{*} $.

\item
$ a \in A $ is called \emph{unitary} if and only if $ a^{*} a = a a^{*} = 1_{A} $, assuming that $ A $ is unital.

\item
If $ A $ is unital, denote the set of unitary elements of $ A $ by $ \U{A} $. It is clear that $ \U{A} $ is a group with respect to multiplication in $ A $.
\end{itemize}

Let $ A $ and $ B $ be $ C^{*} $-algebras. Then a map $ \varphi: A \to B $ is called a \emph{$ * $-homomorphism} if and only if $ \varphi $ is a $ \C{} $-algebra homomorphism that satisfies $ \func{\varphi}{a^{*}} = \func{\varphi}{a}^{*} $ for every $ a \in A $. By spectral theory, $ * $-homomorphisms are bounded with norm $ \leq 1 $, and injective $ * $-homomorphisms are isometric.

An injective $ * $-homomorphism from a $ C^{*} $-algebra to another is sometimes called a \emph{$ * $-embedding}.

A bijective $ * $-homomorphism from a $ C^{*} $-algebra to another is called a \emph{$ * $-isomorphism}.

A $ * $-isomorphism from a $ C^{*} $-algebra to itself is called a \emph{$ * $-automorphism}, and we denote the set of $ * $-automorphisms on a $ C^{*} $-algebra $ A $ by $ \Aut{A} $, which is a group under composition.

Let $ X $ be a locally compact Hausdorff (l.c.H.) space, i.e., a Hausdorff space for which every point has a compact neighborhood. Let $ \Cont{0}{X} $ denote the set of all continuous functions $ f: X \to \C{} $ where for each $ \epsilon > 0 $, there exists a compact subset $ K $ of $ X $ such that $ \Abs{\func{f}{x}} < \epsilon $ for every $ x \in X \setminus K $. Then $ \Cont{0}{X} $ is a commutative $ C^{*} $-algebra when it is equipped with the usual pointwise operations (multiplication, scalar multiplication and conjugation) and the supremum norm.

By the famous Gelfand-Naimark Theorem, a commutative $ C^{*} $-algebra is $ * $-isomorphic to $ \Cont{0}{X} $ for some l.c.H. space $ X $, with $ X $ being compact if and only if the $ C^{*} $-algebra is unital.

The complex field $ \C{} $ is a unital $ C^{*} $-algebra (with complex conjugation serving as the involution), and we have the $ * $-isomorphism $ \C{} \cong \Cont{0}{\textnormal{pt}} $, where $ \textnormal{pt} $ denotes any one-point space.

Let $ \mathcal{H} $ be a Hilbert space, i.e., a $ \C{} $-vector space with a complete (conjugate-linear) inner product $ \Inner{\cdot}{\cdot}{\mathcal{H}} $. Let $ \B{\mathcal{H}} $ denote the set of all bounded operators on $ \mathcal{H} $. Then by the Riesz-Fr\'echet Theorem, every $ T \in \B{\mathcal{H}} $ has an adjoint, i.e., an operator $ T^{*} \in \B{\mathcal{H}} $ (necessarily unique) such that
$$
\forall v,w \in \mathcal{H}: \quad
\Inner{\func{T}{v}}{w}{\mathcal{H}} = \Inner{v}{\func{T^{*}}{w}}{\mathcal{H}}.
$$
Observe that $ \B{\mathcal{H}} $ is a $ C^{*} $-algebra, with composition of operators serving as the multiplication, the operator-adjoint as the involution, and the operator-norm as the $ C^{*} $-algebraic norm. If $ \mathcal{H} = \C{n} $ for some $ n \in \N{} $, this means that the $ n $-dimensional matrix algebra $ \func{M_{n}}{\C{}} $ is a $ C^{*} $-algebra.

The Gelfand-Naimark-Segal (GNS) Construction states that a $ C^{*} $-algebra is $ * $-isomorphic to an operator-norm-closed $ * $-subalgebra of $ \B{\mathcal{H}} $ for some Hilbert space $ \mathcal{H} $. Sometimes, $ C^{*} $-algebras are defined in this manner, but it makes more sense to reserve the name `operator algebras' for such concrete realizations of a $ C^{*} $-algebra.

\subsubsection{Positivity} 

Let $ A $ be a $ C^{*} $-algebra. Then $ a \in A $ is called \emph{positive} if and only if $ a = b^{*} b $ for some $ b \in A $. Given a Hilbert space $ \mathcal{H} $, this is consistent with calling an operator $ T \in \B{\mathcal{H}} $ \emph{positive} if and only if $ T = S^{*} \circ S $ for some $ S \in \B{\mathcal{H}} $. The set of all positive elements of $ A $, which we denote by $ A_{\geq} $, forms a positive cone, i.e.,
$$
A_{\geq} + A_{\geq} \subseteq A_{\geq} \qquad \text{and} \qquad
\RR{\geq 0}{} \cdot A_{\geq} \subseteq A_{\geq}.
$$
Hence, there is a partial order $ \leq_{A} $ on $ A_{\geq} $ given by $ a \leq_{A} b \iff b - a \in A_{\geq} $ for every $ a,b \in A_{\geq} $.

\subsubsection{Ideals} 

Let $ A $ be a $ C^{*} $-algebra. By an \emph{ideal} of $ A $, we mean a two-sided norm-closed algebraic ideal of $ A $.

Any ideal $ J $ of $ A $ is closed under involution (a non-trivial fact) and therefore a $ C^{*} $-algebra itself. The quotient algebra $ A / J $ is then a $ C^{*} $-algebra with the following properties:
\begin{itemize}
\item
The involution is defined by $ \Br{a + J}^{*} = a^{*} + J $ for every $ a \in A $ (well-defined because $ J^{*} = J $).

\item
The norm is defined by $ \D \Norm{a + J}{A / J} \df \inf_{x \in J} \Norm{a + x}{A} $ for every $ a \in A $.
\end{itemize}
We say that $ J $ is \emph{essential} if and only if $ a I = \SSet{0_{A}} = I a \iff a = 0_{A} $ for every $ a \in A $.

We say that $ A $ is \emph{simple} if and only if the only ideals of $ A $ are $ \SSet{0_{A}} $ and $ A $ itself. For any $ n \in \N{} $, the matrix algebra $ \func{M_{n}}{\C{}} $ is a simple $ C^{*} $-algebra.

Let $ \mathcal{H} $ be a Hilbert space. By a \emph{compact operator} on $ \mathcal{H} $, we mean an operator $ T $ on $ \mathcal{H} $ such that $ \Cl{\Im{T}{\Ball{\mathcal{H}}}}{\mathcal{H}} $ --- the closure of the $ T $-image of the open unit ball $ \Ball{\mathcal{H}} $ of $ \mathcal{H} $ --- is a compact subset of $ \mathcal{H} $. For every $ v,w \in \mathcal{H} $, we can define a compact operator $ \ket{v} \bra{w} $ on $ \mathcal{H} $ by
$$
\forall x \in \mathcal{H}: \quad
\func{\ket{v} \bra{w}}{x} \df \Inner{w}{x}{\mathcal{H}} v.
$$
These are called \emph{rank-1 operators} as the dimension of their range space is $ 1 $ (assuming $ v,w \neq 0_{\mathcal{H}} $).

If $ \Comp{\mathcal{H}} $ denotes the set of all compact operators on $ \mathcal{H} $, then $ \Comp{\mathcal{H}} $ is not just a subset but also an ideal of $ \B{\mathcal{H}} $. Furthermore, it can be shown that $ \Comp{\mathcal{H}} = \Cl{\Span{\Set{\ket{v} \bra{w}}{v,w \in \mathcal{H}}}}{\B{\mathcal{H}}} $.

\subsubsection{Approximate Identities} 

Let $ A $ be a $ C^{*} $-algebra. An \emph{approximate identity} for $ A $ is defined as a net $ \Seq{e_{i}}{i \in I} $ in $ A $ such that
$$
\forall a \in A: \quad
  \lim_{i \in I} e_{i} a
= \lim_{i \in I} a e_{i}
= a.
$$
If $ A $ is unital, then the sequence $ \Seq{1_{A}}{n \in \N{}} $ is an approximate identity. Even if $ A $ is not unital, it still has an approximate identity, which we may arrange to consist of positive elements norm-bounded by $ 1 $. We will assume that approximate identities have this special form, unless otherwise specified.

If $ A $ is separable (i.e., it has a countable dense subset), then it possesses an approximate identity that is not just a net but also a sequence.

\subsection{Hilbert $ C^{*} $-Modules} 

Hilbert $ C^{*} $-modules are generalizations of Hilbert spaces. They are extremely important for the structural analysis of $ C^{*} $-algebras, as they are used to define key $ C^{*} $-algebraic concepts such as Morita-Rieffel equivalence and operator $ KK $-theory.

\textbf{Throughout this subsection, $ A $, $ B $ and $ C $ denote arbitrary $ C^{*} $-algebras.}

\subsubsection{Right Hilbert $ C^{*} $-Modules} 

A \emph{right Hilbert $ A $-module} is a vector space $ \X $ endowed with a right $ A $-action $ \bullet: \X \times A \to \X $ and an $ A $-valued map $ \Inner{\cdot}{\cdot}{}: \X \times \X \to A $, called a \emph{right $ A $-inner product}, satisfying the following axioms:
\begin{itemize}
\item[(1)]
$ \Inner{x}{c y + z}{} = c \Inner{x}{y}{} + \Inner{x}{z}{} $ for every $ x,y,z \in \X $ and $ c \in \C{} $.

\item[(2)]
$ \Inner{y}{x}{} = \Inner{x}{y}{}^{*} $ for every $ x,y \in \X $.

\item[(3)]
$ \Inner{x}{y \bullet a}{} = \Inner{x}{y}{} a $ for every $ x,y \in \X $ and $ a \in A $.

\item[(4)]
$ \Inner{x}{x}{} \geq_{A} 0_{A} $ for every $ x \in X $.

\item[(5)]
$ \Inner{x}{x}{} = 0_{A} \iff x = 0_{\X} $ for every $ x \in \X $.

\item[(6)]
$ \X $ is complete with respect to the norm $ \Norm{\cdot}{\X} $ defined by $ \Norm{x}{\X} \df \Norm{\Inner{x}{x}{}}{A}^{\frac{1}{2}} $ for every $ x \in \X $.
\end{itemize}
From these axioms, it follows that
\begin{alignat*}{2}
\forall x,y,z \in \X, ~ \forall c \in \C{}: & \quad &
\Inner{c x + y}{z}{}     & = \Conj{c} \Inner{x}{z}{} + \Inner{y}{z}{} \\
\forall x,y \in \X, ~ \forall a \in A:      & \quad &
\Inner{x \bullet a}{y}{} & = a^{*} \Inner{x}{y}{}.
\end{alignat*}
In order to reflect the dependence of $ \Inner{\cdot}{\cdot}{} $ on $ \X $, we will write it as $ \Inner{\cdot}{\cdot}{\X} $.

If $ \X $ satisfies all of the conditions except (6), then we call it a \emph{right pre-Hilbert $ A $-module}.

It is easily shown that $ J \df \Cl{\Span{\Inner{\X}{\X}{\X}}}{A} $ is an ideal of $ A $. If $ J = A $, then we say that $ \X $ is \emph{full}.

Every Hilbert space is a right Hilbert $ \C{} $-module in the obvious manner.

Now, $ A $ itself is a right Hilbert $ A $-module, with the right $ A $-action being right-multiplication by elements of $ A $, and the right $ A $-inner product defined by $ \Inner{a}{b}{A} \df a^{*} b $ for every $ a,b \in A $. Whenever we want to express $ A $ as a right Hilbert $ A $-module, we will write it as $ A_{A} $.

\subsubsection{Left Hilbert $ C^{*} $-Modules} 

A \emph{left Hilbert $ A $-module} is a vector space $ \X $ endowed with a left $ A $-action $ \bullet: A \times \X \to \X $ and an $ A $-valued map $ \Inner{\cdot}{\cdot}{}: \X \times \X \to A $, called a \emph{left $ A $-inner product}, satisfying (2), (4), (5) and (6) above as well as the following ones:
\begin{enumerate}
\item[(2')]
$ \Inner{c x + y}{z}{} = c \Inner{x}{z}{} + \Inner{y}{z}{} $ for every $ x,y,z \in \X $ and $ c \in \C{} $.

\item[(3')]
$ \Inner{a \bullet x}{y}{} = a \Inner{x}{y}{} $ for every $ x,y \in \X $ and $ a \in A $.
\end{enumerate}
In order to reflect the dependence of $ \Inner{\cdot}{\cdot}{} $ on $ \X $, we will write it as $ \LeftInner{\cdot}{\cdot}{\X} $.

As before, $ J \df \Cl{\Span{\LeftInner{\X}{\X}{\X}}}{A} $ is an ideal of $ A $, and if $ J = A $, then we say that $ \X $ is \emph{full}.

Every Hilbert space $ \mathcal{H} $ is a left $ \Comp{\mathcal{H}} $-module with the following properties:
\begin{itemize}
\item
The left $ \Comp{\mathcal{H}} $-action is defined by $ T \bullet v \df \func{T}{v} $ for every $ T \in \Comp{\mathcal{H}} $ and $ v \in \mathcal{H} $.

\item
The left $ \Comp{\mathcal{H}} $-inner product is defined by $ \LeftInner{v}{w}{\mathcal{H}} \df \ket{v} \bra{w} $ for every $ v,w \in \mathcal{H} $.
\end{itemize}

\subsubsection{A Matter of Terminology} 

We will mostly be dealing with right Hilbert $ C^{*} $-modules, as is often the case throughout literature. When referring to right Hilbert $ C^{*} $-modules, it is common practice to omit the adjective \emph{right} unless such an omission would lead to confusion.

\subsubsection{Adjointable Operators} 

Let $ \X $ and $ \Y $ be Hilbert $ A $-modules. Then a map $ T: \X \to \Y $ is called \emph{adjointable} if and only if it has an adjoint, i.e., a map $ T^{*}: \Y \to \X $ (necessarily unique) such that
$$
\forall x \in \X, ~ \forall y \in \Y: \quad
\Inner{\func{T}{x}}{y}{\Y} = \Inner{x}{\func{T^{*}}{y}}{\X}.
$$
It follows readily from the definition of adjointability that $ T^{*} $ is adjointable with $ T^{**} = T $. An easy argument shows that $ T $ is linear and $ A $-linear (i.e., $ \func{T}{x \bullet a} = \func{T}{x} \bullet a $ for every $ x \in \X $ and $ a \in A $), and it is bounded as well by the Closed Graph Theorem. Hence, $ T $ is a bounded operator. However, a bounded operator from $ \X $ to $ \Y $ is not necessarily adjointable --- a counterexample was constructed by W. Paschke in \cite{Paschke}.

Denote the set of adjointable operators from $ \X $ to $ \Y $ by $ \AdjPair{\X}{\Y} $, and write $ \Adj{\X} $ for $ \AdjPair{\X}{\X} $. Note that $ \Adj{X} $ is a $ C^{*} $-algebra in much the same way that $ \B{\mathcal{H}} $ is one for any Hilbert space $ \mathcal{H} $.

We can also define an adjointable operator on a left Hilbert $ C^{*} $-module in an analogous fashion, but this notion is rarely used nowadays.

\subsubsection{Compact Operators} 

Let $ \X $ be a Hilbert $ A $-module. For every $ \Pair{x}{y} \in \X \times \X $, define an operator $ \ket{x} \bra{y} $ on $ \X $ by
$$
\forall z \in \X: \quad
\func{\ket{x} \bra{y}}{z} \df x \bullet \Inner{y}{z}{\X}.
$$
It is not hard to show that the set $ \Span{\Set{\ket{x} \bra{y}}{x,y \in \X}} $ is a two-sided algebraic ideal of $ \Adj{X} $. Its closure $ \Cl{\Span{\Set{\ket{x} \bra{y}}{x,y \in \X}}}{\Adj{X}} $ in $ \Adj{X} $ is thus an ideal of $ \Adj{X} $, which we denote by $ \Comp{\X} $. Call any element of $ \Comp{\X} $ a \emph{compact operator} on $ \X $. Note that $ \X $ is a full left Hilbert $ \Comp{\X} $-module, where the left $ \Comp{\X} $-inner product $ \LeftInner{\cdot}{\cdot}{\Comp{\X}} $ is defined by $ \LeftInner{x}{y}{\Comp{\X}} \df \ket{x} \bra{y} $ for every $ x,y \in \X $.

If $ \X $ is a Hilbert space, then the definition of a compact operator on $ \X $ given here coincides with the earlier definition of a Hilbert-space compact operator.

We can also define a compact operator on a left Hilbert $ C^{*} $-module in a similar manner.

In \cite{Rieffel2}, Rieffel used the term \emph{imprimitivity algebra} to refer to the algebra of compact operators on a left/right Hilbert $ C^{*} $-module.

\subsubsection{Multiplier Algebras} 

Define the \emph{multiplier algebra} of $ A $ to be $ \Adj{A_{A}} $, and denote it by $ \Mult{A} $. This is a unital $ C^{*} $-algebra whose identity element is $ \Id_{A} $ (the identity operator on $ A $).

There is an injective $ * $-homomorphism $ L: A \hookrightarrow \Mult{A} $ defined by
$$
\forall a \in A: \quad
L_{a} \df \Map{A}{A}{x}{a x}.
$$
(It is easily verified that $ L_{a} $ is adjointable with $ L_{a}^{*} = L_{a^{*}} $ for every $ a \in A $.) If $ A $ is already unital, then $ L $ is surjective. To see why, let $ T \in \Adj{A_{A}} $, so that $ \Inner{\func{T}{a}}{b}{A} = \Inner{a}{\func{T^{*}}{b}}{A} $ for every $ a,b \in A $. In particular,
$$
\forall a \in A: \quad
  \func{T}{a}^{*}
= \func{T}{a}^{*} 1_{A}
= \Inner{\func{T}{a}}{1_{A}}{A}
= \Inner{a}{\func{T^{*}}{1_{A}}}{A}
= a^{*} \func{T^{*}}{1_{A}}, \quad \text{so} \quad
  \func{T}{a} = \func{T^{*}}{1_{A}}^{*} a.
$$
Hence, $ T = L_{\func{T^{*}}{1_{A}}^{*}} $, and as $ T $ is arbitrary, $ \Adj{A_{A}} \subseteq \Range{L} $. However, $ \Range{L} \subseteq \Adj{A_{A}} $, which yields $ \Range{L} = \Adj{A_{A}} $.

Now, $ \Set{L_{a}}{a \in A} $ is an essential ideal of $ \Mult{A} $, and $ \Mult{A} $ has the following universal property: For every $ C^{*} $-algebra $ B $, if there is a $ * $-embedding $ j: A \hookrightarrow B $ where $ \Im{j}{A} $ is an essential ideal of $ B $, then there exists a unique $ * $-embedding $ \iota: B \hookrightarrow \Mult{A} $ satisfying $ \func{\iota}{\func{j}{a}} = L_{a} $ for every $ a \in A $.

By abuse of notation, we usually view $ A $ as a $ C^{*} $-subalgebra of $ \Mult{A} $ and write $ A \subseteq \Mult{A} $.

\subsection{Non-Degeneracy} 

A $ * $-homomorphism $ \varphi: A \to B $ is called \emph{non-degenerate} if and only if $ \Span{\Im{\varphi}{A} B} $ is dense in $ B $.

If $ \pi: A \to B $ is a non-degenerate $ * $-homomorphism, then there exists a unique $ * $-homomorphism $ \Conj{\varphi}: \Mult{A} \to \Mult{B} $ such that $ \Conj{\varphi}|_{A} = \varphi $.

Due to the existence of approximate identities, $ * $-isomorphisms are automatically non-degenerate.

A \emph{$ * $-representation} of $ A $ on a Hilbert $ B $-module $ \X $ is defined as a $ * $-homomorphism $ \pi: A \to \Adj{\X} $. We say that $ \pi $ is \emph{faithful} if and only if it is injective, and that it is \emph{non-degenerate} if and only if the set $ \Span{\Set{\FUNC{\func{\pi}{a}}{x}}{a \in A ~ \text{and} ~ x \in \X}} $ is dense in $ \X $.

We may define $ \Mult{A} $ via non-degenerate $ * $-representations. Let $ \X $ be a Hilbert $ B $-module, and suppose that $ \pi $ is a faithful and non-degenerate $ * $-representation of $ A $ on $ \X $. Then take $ \Mult{A} $ to be the idealizer of $ \Im{\pi}{A} $ in $ \Adj{\X} $:
$$
\Mult{A} \df \Set{T \in \Adj{\X}}{T \circ \Im{\pi}{A} \subseteq \Im{\pi}{A} ~ \text{and} ~ \Im{\pi}{A} \circ T \subseteq \Im{\pi}{A}}.
$$

\subsection{Morita-Rieffel Equivalence} 

The concept of Morita equivalence originates from ring theory. Two (unital) rings $ R $ and $ S $ are called \emph{Morita equivalent} if and only if there is an equivalence between the category of left $ R $-modules and the category of left $ S $-modules.

As the representations of a ring are given by the left modules over that ring, we can say that Morita-equivalent rings have the same representation theory.

While attempting to replace the unintuitive measure-theoretical foundation of George Mackey's theory of induced representations of locally compact Hausdorff (l.c.H.) groups by a more natural algebraic one, Rieffel was led to develop a specialized notion of Morita equivalence for $ C^{*} $-algebras, taking into account the presence of an involution and the fact that $ C^{*} $-algebras may be non-unital --- we say that $ A $ and $ B $ are \emph{Morita-Rieffel equivalent} (or \emph{strongly Morita equivalent}) if and only if there is an $ \Pair{A}{B} $-bimodule $ \X $ with the following properties:
\begin{itemize}
\item
$ \X $ is a full left Hilbert $ A $-module \emph{and} a full right Hilbert $ B $-module.

\item
$ \LeftInner{x}{y}{\X} \bullet z = x \bullet \Inner{y}{z}{\X} $ for every $ x,y,z \in \X $.
\end{itemize}
We then call $ \X $ an \emph{$ \Pair{A}{B} $-imprimitivity bimodule}.

A simple example of Morita-equivalent $ C^{*} $-algebras are $ \C{} $ and $ \Comp{\mathcal{H}} $, for any Hilbert space $ \mathcal{H} $. Hence, a non-commutative $ C^{*} $-algebra may be Morita-Rieffel equivalent to a commutative one.

Morita-Rieffel equivalence is a genuine equivalence relation on the class of $ C^{*} $-algebras:
\begin{itemize}
\item
Reflexivity comes from the fact that $ A $ itself is naturally an $ \Pair{A}{A} $-imprimitivity bimodule.

\item
If $ \X $ is an imprimitivity $ \Pair{A}{B} $-bimodule, then we can form an imprimitivity $ \Pair{B}{A} $-bimodule $ \widetilde{\X} $, called the \emph{dual module} of $ \X $, whose underlying vector space is the complex conjugate of that of $ \X $ and with operations involving $ A $ and $ B $ interchanged in a certain manner. This yields symmetry.

\item
If $ \X $ is an $ \Pair{A}{B} $-imprimitivity bimodule and $ \Y $ a $ \Pair{B}{C} $-imprimitivity bimodule, then by forming the $ B $-balanced algebraic tensor product $ \X \Odot{B} \Y $ and completing it with respect to a certain norm, we obtain an imprimitivity $ \Pair{A}{C} $-bimodule, denoted by $ \X \Otimes{B} \Y $. This gives us transitivity.
\end{itemize}

Morita-Rieffel equivalence is weaker than $ * $-isomorphism ($ \C{} $ and $ \Comp{\mathcal{H}} $ are not $ * $-isomorphic if $ \dim{\mathcal{H}} > 1 $), but it preserves many $ C^{*} $-algebraic properties. For example, Morita-Rieffel equivalent $ C^{*} $-algebras possess isomorphic ideal-lattice structures and isomorphic $ K $-theories.

The utility of the concept thus comes from the fact that if a non-commutative $ C^{*} $-algebra $ A $ is Morita-Rieffel equivalent to a commutative $ C^{*} $-algebra $ B $, then properties about $ A $ can be derived from $ B $, which is a more easily studied object.


\section{Topological Dynamical Systems and $ C^{*} $-Dynamical Systems} 

From now on, we will assume that every l.c.H. group is equipped with a choice of Haar measure $ \mu $, with respect to which integration can be carried out. The modular function corresponding to the group will be denoted by $ \Delta $.

\subsection{Topological Dynamical Systems} 


\begin{Def} \label{Topological Dynamical Systems}
A \emph{topological dynamical system} is a triple $ \Trip{G}{X}{\alpha} $, where:
\begin{itemize}
\item
$ G $ is an l.c.H. group.

\item
$ X $ is an l.c.H. space.

\item
$ \alpha: G \times X \to X $ is a jointly continuous $ G $-action on $ X $.
\end{itemize}
Denote the orbit space of $ \alpha $ by $ \Orb{G}{X}{\alpha} $, and equip it with the obvious quotient topology.
\end{Def}



\begin{Eg} \label{Classical Dynamical Systems}
The concept of a dynamical system originates from classical mechanics. In this context, a real dynamical system is a pair $ \Pair{M}{\Phi} $, where $ M $ is a smooth manifold, called the \emph{phase space}, and $ \Phi: \RR{}{} \times M \to M $ a smooth function with the following properties:
\begin{itemize}
\item
$ \func{\Phi}{0,x} = x $ for every $ x \in M $.

\item
$ \func{\Phi}{s + t,x} = \func{\Phi}{s,\func{\Phi}{t,x}} $ for every $ s,t \in \RR{}{} $ and $ x \in M $.
\end{itemize}
These properties imply that $ \Trip{\RR{}{}}{M}{\Phi} $ is a topological dynamical system.
\end{Eg}



\begin{Eg} \label{Translation in the Plane}
Fix a $ \theta \in \RR{}{} $, and define an $ \RR{}{} $-action $ \tau^{\theta} $ on $ \RR{}{2} $ by
$$
\forall \theta \in \RR{}{}, ~ \forall \Pair{x}{y} \in \RR{}{2}: \quad
\func{\tau^{\theta}}{r,\Pair{x}{y}} \df \Pair{x + r}{y + \theta r}.
$$
Then $ \Trip{\RR{}{}}{\RR{}{2}}{\tau^{\theta}} $ is a topological dynamical system, and $ \Orb{\RR{}{}}{\RR{}{2}}{\tau^{\theta}} $ is homeomorphic to $ \RR{}{} $.
\end{Eg}



\begin{Eg} \label{Translation in the 2-Torus}
Fix a $ \theta \in \RR{}{} $, and define an $ \RR{}{} $-action $ \tau^{\theta} $ on the torus $ \T{2} \df \RR{}{2} / \Z{2} $ by
$$
\forall r \in \RR{}{}, ~ \forall \Pair{x}{y} \in \RR{}{2}: \quad
\func{\tau^{\theta}}{r,\Pair{x}{y} + \Z{2}} \df \Pair{x + r}{y + \theta r} + \Z{2}.
$$
Then $ \Trip{\RR{}{}}{\T{2}}{\tau^{\theta}} $ is a topological dynamical system. If $ \theta \in \mathbb{Q} $, then $ \Orb{\RR{}{}}{\T{2}}{\tau^{\theta}} $ is homeomorphic to $ \T{} $, but if $ \theta \in \RR{}{} \setminus \mathbb{Q} $, then each orbit is dense in $ \T{2} $, so $ \Orb{\RR{}{}}{\T{2}}{\tau^{\theta}} $ is an indiscrete space.
\end{Eg}



\begin{Eg} \label{Self-Homeomorphisms}
Let $ X $ be an l.c.H. space and $ h: X \to X $ an arbitrary self-homeomorphism. Defining $ \sigma: \Z{} \times X \to X $ by
$$
\forall \Pair{n}{x} \in \Z{} \times X: \quad
\func{\sigma}{n,x} \df \func{h^{n}}{x},
$$
we find that $ \Trip{\Z{}}{X}{\sigma} $ is a topological dynamical system.
\end{Eg}



\begin{Def} \label{Proper Topological Dynamical Systems}
We say that a topological dynamical system $ \Trip{G}{X}{\alpha} $ is \emph{proper} if and only if the continuous function $ \Map{G \times X}{X \times X}{\Pair{r}{x}}{\Pair{x}{\func{\alpha}{r,x}}} $ is topologically proper, i.e., the pre-image of every compact subset of $ X \times X $ is a compact subset of $ G \times X $.
\end{Def}


The topological dynamical system in \autoref{Translation in the Plane} is easily checked to be proper for every $ \theta \in \RR{}{} $.

If $ \Trip{G}{X}{\alpha} $ is a proper topological dynamical system, then $ \Orb{G}{X}{\alpha} $ is an l.c.H. space, making $ \Cont{0}{\Orb{G}{X}{\alpha}} $ a commutative $ C^{*} $-algebra. In particular, the topological dynamical system in \autoref{Translation in the 2-Torus} is not proper for any $ \theta \in \RR{}{} $, although the associated orbit space is locally compact and Hausdorff when $ \theta \in \mathbb{Q} $.

\subsection{$ C^{*} $-Dynamical Systems and Crossed Products} 

\subsubsection{$ C^{*} $-Dynamical Systems} 


\begin{Def} \label{C*-Dynamical Systems}
A \emph{$ C^{*} $-dynamical system} is a triple $ \Trip{G}{A}{\alpha} $, where:
\begin{itemize}
\item
$ G $ is an l.c.H. group.

\item
$ A $ is a $ C^{*} $-algebra.

\item
$ \alpha $ is a strongly continuous group homomorphism from $ G $ to $ \Aut{A} $, i.e., the function
$$
\Map{G}{A}{r}{\alphArg{r}{a}}
$$
is continuous for each $ a \in A $.
\end{itemize}
We say that $ \Trip{G}{A}{\alpha} $ is \emph{commutative} if and only if $ A $ is commutative.
\end{Def}



\begin{Eg} \label{Transformation Groups}
Let $ \Trip{G}{X}{\alpha} $ be a topological dynamical system. Define $ \widetilde{\alpha}: G \to \Aut{\Cont{0}{X}} $ by
$$
\forall r \in G, ~ \forall f \in \Cont{0}{X}: \quad
\func{\widetilde{\alpha}_{r}}{f} \df \Map{G}{\C{}}{x}{\func{f}{\func{\alpha}{r^{-1},x}}}.
$$
Then $ \Trip{G}{\Cont{0}{X}}{\widetilde{\alpha}} $ is a commutative $ C^{*} $-dynamical system.
\end{Eg}



\begin{Eg} \label{Automorphisms of the Compact Operators}
Let $ G $ be an l.c.H. group and $ \mathcal{H} $ a Hilbert space, and suppose that $ U $ is a continuous unitary representation of $ G $ on $ \mathcal{H} $, i.e., $ U: G \to \U{\B{\mathcal{H}}} $ is a continuous group homomorphism. Define $ \Ad{U}: G \to \Aut{\Comp{\mathcal{H}}} $ by
$$
\forall r \in G, ~ \forall T \in \Comp{\mathcal{H}}: \quad
\FUNC{\Ad{U}_{r}}{T} \df U_{r} \circ T \circ U_{r}^{*}.
$$
Then $ \Trip{G}{\Comp{\mathcal{H}}}{\Ad{U}} $ is a commutative $ C^{*} $-dynamical system.
\end{Eg}


\subsubsection{Crossed Products} \label{Crossed Products} 

Let $ \Trip{G}{A}{\alpha} $ be a $ C^{*} $-dynamical system. From it, we can construct two kinds of $ C^{*} $-algebras, called the \emph{full crossed product} and the \emph{reduced crossed product}. To define these $ C^{*} $-algebras, let $ \Cc{G,A} $ denote the complex vector space of compactly supported continuous $ A $-valued functions on $ G $, and define two operations, $ \star: \Cc{G,A} \times \Cc{G,A} \to \Cc{G,A} $ and $ ^{*}: \Cc{G,A} \to \Cc{G,A} $, by
\begin{align*}
\forall f,g \in \Cc{G,A}: \quad
f \star g & \df \Map{G}{A}{x}{\Int{G}{\func{f}{y} ~ \alphArg{y}{\func{g}{y^{-1} x}}}{y}}, \qquad \Br{\text{Convolution}} \\
f^{*}     & \df \Map{G}{A}{x}{\Del{x}{-1} \alphArg{x}{\func{f}{x^{-1}}}^{*}}.             \qquad \Br{\text{Involution}}
\end{align*}
If $ \Norm{\cdot}{1} $ denotes the $ L^{1} $-norm on $ \Cc{G,A} $, i.e., $ \D \Norm{f}{1} = \Int{G}{\Norm{\func{f}{x}}{A}}{x} $ for every $ f \in \Cc{G,A} $, then $ \Quad{\Cc{G,A}}{\star}{^{*}}{\Norm{\cdot}{1}} $ is a normed $ * $-algebra. We denote the completion of $ \Cc{G,A} $ with respect to $ \Norm{\cdot}{1} $ by $ \LL{1}{G,A} $.

Now, a \emph{covariant representation} of $ \Trip{G}{A}{\alpha} $ is defined (see \cite{Williams}) as a triple $ \Trip{\X}{\pi}{U} $, where:
\begin{itemize}
\item
$ \X $ is a Hilbert $ B $-module for some $ C^{*} $-algebra $ B $.

\item
$ \pi $ is a $ * $-representation of $ A $ on $ \X $.

\item
$ U $ is a strongly continuous unitary representation of $ G $ on $ \X $, i.e., $ U $ is a group homomorphism from $ G $ to $ \U{\Adj{\X}} $ such that the function $ \Map{G}{\X}{r}{\func{U_{r}}{x}} $ is continuous for each $ x \in \X $.

\item
$ \func{\pi}{\alphArg{r}{a}} = U_{r} \circ \func{\pi}{a} \circ U_{r}^{*} $ for every $ r \in G $ and $ a \in A $.
\end{itemize}

For each covariant representation $ \Trip{\X}{\pi}{U} $ of $ \Trip{G}{A}{\alpha} $, we can define a $ * $-algebra homomorphism $ \rho_{\X,\pi,U}: \Trip{\Cc{G,A}}{\star}{^{*}} \to \Trip{\Adj{\X}}{\circ}{^{*}} $, called the \emph{integrated form} of $ \Trip{\X}{\pi}{U} $, by
$$
\forall f \in \Cc{G,A}: \quad
\func{\rho_{\X,\pi,U}}{f} \df \Map{\X}{\X}{x}{\Int{G}{\FUNC{\func{\pi}{\func{f}{y}}}{\func{U_{y}}{x}}}{y}},
$$
which yields the following norm inequality:
\begin{align*}
\forall f \in \Cc{G,A}: \quad
       \Norm{\func{\rho_{\X,\pi,U}}{f}}{\Adj{\X}}
& =    \sup_{\Norm{x}{\X} = 1} \Norm{\Int{G}{\FUNC{\func{\pi}{\func{f}{y}}}{\func{U_{y}}{x}}}{y}}{\X} \\
& \leq \sup_{\Norm{x}{\X} = 1} \Int{G}{\Norm{\FUNC{\func{\pi}{\func{f}{y}}}{\func{U_{y}}{x}}}{\X}}{y} \\
& \leq \sup_{\Norm{x}{\X} = 1} \Int{G}{\Norm{\func{\pi}{\func{f}{y}}}{\Adj{\X}} \Norm{\func{U_{y}}{x}}{\X}}{y} \\
& \leq \sup_{\Norm{x}{\X} = 1} \Int{G}{\Norm{\func{f}{y}}{A} \Norm{x}{\X}}{y} \\
& =    \Norm{f}{1}.
\end{align*}
Hence, if there exists a covariant representation $ \Trip{\X}{\pi}{U} $ of $ \Trip{G}{A}{\alpha} $ such that $ \rho_{\X,\pi,U} $ is injective, then we can define a new norm $ \Norm{\cdot}{u} $ on $ \Cc{G,A} $, called the \emph{universal norm}, by
$$
\forall f \in \Cc{G,A}: \quad
    \Norm{f}{u}
\df \func{\sup}{\Set{\Norm{\func{\rho_{\X,\pi,U}}{f}}{\Adj{\X}}}{\text{$ \Trip{\pi}{U}{\X} $ is a covariant rep. of $ \Trip{G}{A}{\alpha} $}}}.
$$

Let $ \X $ be a Hilbert $ B $-module for some $ C^{*} $-algebra $ B $, and let $ \LL{2}{G,\X} $ denote the completion of $ \Cc{G,\X} $ with respect to the norm $ \NNorm{\cdot} $ defined by
$$
\forall \phi \in \Cc{G,\X}: \quad
\NNorm{\phi} \df \Norm{\Int{G}{\Inner{\func{\phi}{x}}{\func{\phi}{x}}{\X}}{x}}{B}^{\frac{1}{2}}.
$$
Let $ q $ denote the canonical dense linear embedding of $ \Cc{G,\X} $ into $ \LL{2}{G,\X} $. Then $ \LL{2}{G,\X} $ is a Hilbert $ B $-module in the following manner:
\begin{itemize}
\item
For every $ b \in B $ and $ \phi \in \Cc{G,\X} $, we have
\begin{align*}
     & ~ \Norm{\func{q}{\Map{G}{\X}{x}{\func{\phi}{x} \bullet b}}}{\LL{2}{G,\X}} \\
=    & ~ \NNorm{\Map{G}{\X}{x}{\func{\phi}{x} \bullet b}} \\
=    & ~ \Norm{\Int{G}{\Inner{\func{\phi}{x} \bullet b}{\func{\phi}{x} \bullet b}{\X}}{x}}{B}^{\frac{1}{2}} \\
=    & ~ \Norm{\Int{G}{b^{*} \Inner{\func{\phi}{x}}{\func{\phi}{x}}{\X} b}{x}}{B}^{\frac{1}{2}} \\
=    & ~ \Norm{b^{*} \SqBr{\Int{G}{\Inner{\func{\phi}{x}}{\func{\phi}{x}}{\X}}{x}} b}{B}^{\frac{1}{2}}
         \qquad \Br{\text{As multiplication in $ B $ is continuous.}} \\
\leq & ~ \Norm{b^{*}}{B}^{\frac{1}{2}}
         \Norm{\Int{G}{\Inner{\func{\phi}{x}}{\func{\phi}{x}}{\X}}{x}}{B}^{\frac{1}{2}}
         \Norm{b}{B}^{\frac{1}{2}} \\
=    & ~ \Norm{b}{B}^{\frac{1}{2}} \Norm{\Int{G}{\Inner{\func{\phi}{x}}{\func{\phi}{x}}{\X}}{x}}{B}^{\frac{1}{2}} \Norm{b}{B}^{\frac{1}{2}} \\
=    & ~ \Norm{b}{B} \Norm{\Int{G}{\Inner{\func{\phi}{x}}{\func{\phi}{x}}{\X}}{x}}{B}^{\frac{1}{2}} \\
=    & ~ \Norm{b}{B} \NNorm{\phi} \\
=    & ~ \Norm{b}{B} \Norm{\func{q}{\phi}}{\LL{2}{G,\X}}.
\end{align*}
We can thus define a right $ B $-action $ \centerdot $ on $ \LL{2}{G,\X} $ by
$$
\forall b \in B, ~ \forall \Phi \in \LL{2}{G,\X}: \quad
\Phi \centerdot b \df \lim_{n \to \infty} ~ \func{q}{\Map{G}{\X}{x}{\func{\phi_{n}}{x} \bullet b}},
$$
where $ \Seq{\phi_{n}}{n \in \N{}} $ is \emph{any} sequence in $ \Cc{G,\X} $ with $ \D \lim_{n \to \infty} \func{q}{\phi_{n}} = \Phi $.

\item
The $ B $-valued bilinear map $ \SqBr{\cdot,\cdot}: \Cc{G,\X} \times \Cc{G,\X} \to B $ defined by
$$
\forall \phi,\psi \in \Cc{G,\X}: \quad
\SqBr{\phi,\psi} \df \Int{G}{\Inner{\func{\phi}{x}}{\func{\psi}{x}}{\X}}{x}
$$
is a $ B $-pre-inner product on $ \Cc{G,\X} $. Hence, by the Cauchy-Schwarz Inequality,
$$
\forall \phi,\psi \in \Cc{G,\X}: \quad
     \Norm{\SqBr{\phi,\psi}}{B}
\leq \Norm{\SqBr{\phi,\phi}}{B}^{\frac{1}{2}} \Norm{\SqBr{\psi,\psi}}{B}^{\frac{1}{2}}
=    \NNorm{\phi} \NNorm{\psi}
=    \Norm{\func{q}{\phi}}{\LL{2}{G,\X}} \Norm{\func{q}{\psi}}{\LL{2}{G,\X}}.
$$
We can thus define a $ B $-inner product $ \Inner{\cdot}{\cdot}{\LL{2}{G,\X}} $ on $ \LL{2}{G,\X} $ by
$$
\forall \Phi,\Psi \in \LL{2}{G,\X}: \quad
\Inner{\Phi}{\Psi}{\LL{2}{G,\X}} \df \lim_{n \to \infty} \Int{G}{\func{\phi_{n}}{x}^{*} \func{\psi_{n}}{x}}{x},
$$
where $ \Seq{\phi_{n}}{n \in \N{}},\Seq{\psi_{n}}{n \in \N{}} $ are \emph{any} sequences in $ \Cc{G,\X} $ with $ \D \lim_{n \to \infty} \func{q}{\phi_{n}} = \Phi $ and $ \D \lim_{n \to \infty} \func{q}{\psi_{n}} = \Psi $.
\end{itemize}
Now, let $ \pi $ be a faithful $ * $-representation of $ A $ on $ \X $. We can then define a $ * $-representation $ \widetilde{\pi} $ of $ A $ and a unitary representation $ \widetilde{U} $ of $ G $ on the Hilbert $ B $-module $ \LL{2}{G,\X} $ as follows:
\begin{itemize}
\item
Let $ a \in A $ and $ \Phi \in \LL{2}{G,\X} $. If $ \Seq{\phi_{n}}{n \in \N{}} $ is any sequence in $ \Cc{G,\X} $ with $ \D \lim_{n \to \infty} \func{q}{\phi_{n}} = \Phi $, then
$$
\FUNC{\func{\widetilde{\pi}}{a}}{\Phi} \df \lim_{n \to \infty} \func{q}{\Map{G}{\X}{x}{\FUNC{\func{\pi}{\alphArg{x}{a}}}{\func{\phi_{n}}{x}}}}.
$$

\item
Let $ r \in G $ and $ \Phi \in \LL{2}{G,\X} $. If $ \Seq{\phi_{n}}{n \in \N{}} $ is any sequence as above, then
$$
\func{\widetilde{U}_{r}}{\Phi} \df \lim_{n \to \infty} \func{q}{\Map{G}{\X}{x}{\Del{r}{\frac{1}{2}} \cdot \func{\phi_{n}}{x r}}}.
$$
\end{itemize}
Then $ \Trip{\LL{2}{G,\X}}{\widetilde{\pi}}{\widetilde{U}} $ is a covariant representation, called a \emph{right-regular representation}, and it is practically $ C^{*} $-folklore that $ \rho_{\LL{2}{G,\X},\widetilde{\pi},\widetilde{U}}: \Cc{G,A} \to \Adj{\LL{2}{G,\X}} $ is injective.

The \emph{reduced crossed product} for $ \Trip{G}{A}{\alpha} $, denoted by $ \RCP{G}{A}{\alpha} $, is now defined as
$$
\Cl{\Range{\rho_{\LL{2}{G,\X},\widetilde{\pi},\widetilde{U}}}}{\Adj{\LL{2}{G,\X}}}
$$
for any faithful $ * $-representation $ \pi $ of $ A $ on a Hilbert $ B $-module $ \X $, for an arbitrary $ C^{*} $-algebra $ B $. Different faithful $ * $-representations will yield the same $ \RCP{G}{A}{\alpha} $ up to $ * $-isomorphism (see \cite{Williams}).

Next, $ \Norm{\cdot}{u} $ satisfies the $ C^{*} $-identity, so the completion of $ \Trip{\Cc{G,A}}{\star}{^{*}} $ with respect to this norm yields a $ C^{*} $-algebra, which we call the \emph{full crossed product} for $ \Trip{G}{A}{\alpha} $ and denote by $ \FCP{G}{A}{\alpha} $.

Although it is rarely mentioned in the literature, our constructions of the crossed products are independent of the Haar measure used, so we do not require any structural data from $ G $ other than its group structure and group topology.

In general, $ \RCP{G}{A}{\alpha} \not\cong \FCP{G}{A}{\alpha} $, unless $ G $ is amenable.

\subsection{Rieffel-Properness} 

The most fundamental theorem on proper topological dynamical systems is perhaps the following, due to Philip Green.


\begin{Thm}[\cite{Green}] \label{The Fundamental Theorem of Proper Topological Dynamical Systems}
Let $ \Trip{G}{X}{\alpha} $ be a proper topological dynamical system. Then $ \Cont{0}{\Orb{G}{X}{\alpha}} $ is Morita-Rieffel equivalent to an ideal $ J $ of $ \CP{G}{\Cont{0}{X}}{\widetilde{\alpha}} $ (which is isomorphic to $ \RCP{G}{\Cont{0}{X}}{\widetilde{\alpha}} $ by Theorem 6.1 of \cite{Phillips}). Also, $ J = \RCP{G}{\Cont{0}{X}}{\widetilde{\alpha}} $ if and only if $ \alpha $ is a free $ G $-action on $ X $.
\end{Thm}


The $ C^{*} $-algebra $ \Cont{0}{\Orb{G}{X}{\alpha}} $ is called a \emph{fixed-point algebra} because each element of $ \Cont{0}{\Orb{G}{X}{\alpha}} $ can be identified with a function in $ \Cb{X} $ that is constant on the $ \alpha $-orbits of $ G $, thus fixed under the $ \alpha $-induced $ G $-action on $ \Cb{X} $.

As $ \Trip{G}{\Cont{0}{X}}{\widetilde{\alpha}} $ is a commutative $ C^{*} $-dynamical system, it makes sense to ask if there exists an analog of \autoref{The Fundamental Theorem of Proper Topological Dynamical Systems} for a non-commutative $ C^{*} $-dynamical system. This requires a definition of properness for an arbitrary $ C^{*} $-dynamical system --- the following one is due to Marc Rieffel.


\begin{Def}[\cite{Rieffel2}] \label{Rieffel-Properness}
We say that a $ C^{*} $-dynamical system $ \Trip{G}{A}{\alpha} $ is \emph{Rieffel-proper} if and only if there exists a dense $ \alpha $-invariant $ * $-subalgebra $ A_{0} $ of $ A $ satisfying the following conditions:
\begin{itemize}
\item
For every $ a,b \in A_{0} $, the following continuous $ A $-valued functions on $ G $ are integrable:
$$
\Map{G}{A}{x}{\Del{x}{- \frac{1}{2}} a ~ \alphArg{x}{b^{*}}} \qquad \text{and} \qquad
\Map{G}{A}{x}{a ~ \alphArg{x}{b^{*}}}.
$$

\item
For every $ a,b \in A_{0} $, there exists an $ m \in \Mult{A} $ --- necessarily unique --- such that:
\begin{itemize}
\item
$ m A_{0} \cup A_{0} m \subseteq A_{0} $, and $ \alphArgExt{x}{m} = m $ for every $ x \in G $, where $ \alphExt{x} $ denotes the extension of $ \alph{x} $ to an automorphism on $ \Mult{A} $.

\item
$ \D c m = \Int{G}{c ~ \alphArg{x}{a^{*} b}}{x} $ for every $ c \in A_{0} $.
\end{itemize}
\end{itemize}
(\textbf{Note:} Our use of the term \emph{Rieffel-proper} is adopted from \cite{Buss|Echterhoff}.)
\end{Def}


Rieffel showed in \cite{Rieffel3} that a topological dynamical system $ \Trip{G}{X}{\alpha} $ is proper if and only if $ \Trip{G}{\Cont{0}{X}}{\widetilde{\alpha}} $ is Rieffel-proper. In other words, Rieffel-properness is equivalent to \autoref{Proper Topological Dynamical Systems} in the case of a commutative $ C^{*} $-dynamical system, which is a clear indication of its success despite its unwieldy appearance.

The main success of \autoref{Rieffel-Properness}, however, lies in the fact that starting from a Rieffel-proper $ C^{*} $-dynamical system $ \Trip{G}{A}{\alpha} $, one can construct a \emph{generalized fixed-point algebra} that is Morita-Rieffel equivalent to an ideal of $ \RCP{G}{A}{\alpha} $. This yields a non-commutative analog of \autoref{The Fundamental Theorem of Proper Topological Dynamical Systems}. For completeness, we now give an outline of the construction.
\begin{itemize}
\item
Let $ \Trip{G}{A}{\alpha} $ be a Rieffel-proper $ C^{*} $-dynamical system.

\item
Let $ A_{0} $ be a dense $ \alpha $-invariant $ * $-subalgebra of $ A $ that satisfies the conditions in \autoref{Rieffel-Properness}.

\item
Define $ \Psi: A_{0} \times A_{0} \to \LL{1}{G,A} $ by $ \func{\Psi}{a,b} \df \Map{G}{A}{x}{\Del{x}{- \frac{1}{2}} a ~ \alphArg{x}{b^{*}}} $ for every $ a,b \in A_{0} $.

\item
Define $ \mu: A_{0} \times A_{0} \to \Mult{A} $ by $ \func{\mu}{a,b} \df m $, where $ m $ is the unique element of $ \Mult{A} $ satisfying the conditions in the second half of \autoref{Rieffel-Properness}.

\item
Let $ E_{0} \df \Span{\Set{\func{\Psi}{a,b}}{a,b \in A_{0}}} $, which is seen to be a $ * $-subalgebra of $ \Quad{\LL{1}{G,A}}{\star}{^{*}}{\Norm{\cdot}{1}} $. Then $ E_{0} $ can be embedded as a $ * $-subalgebra of $ \RCP{G}{A}{\alpha} $, so $ E \df \Cl{E_{0}}{\RCP{G}{A}{\alpha}} $ is well-defined.

\item
Define a left $ E_{0} $-action $ \diamond $ on $ A_{0} $ by $ \func{\Psi}{a,b} \diamond c \df a ~ \func{\mu}{b,c} $ for every $ a,b,c \in A_{0} $. With this action, $ A_{0} $ is a left pre-Hilbert $ E_{0} $-module for $ \func{\Psi}{\cdot,\cdot} $, which defines a left $ E_{0} $-pre-inner product on $ A_{0} $.

\item
Let $ \Norm{\cdot}{A_{0}} $ denote the norm on $ A_{0} $ induced by $ \func{\Psi}{\cdot,\cdot} $. Then $ \diamond $ extends to a left $ E $-action on $ \X $ --- the completion of $ A_{0} $ with respect to $ \Norm{\cdot}{A_{0}} $. Hence, $ \HModule{E}{\X}{} $ is a full left Hilbert $ E $-module for $ \Inner{\cdot}{\cdot}{E} $, which denotes the left $ E $-inner product on $ \X $ that continuously extends $ \func{\Psi}{\cdot,\cdot} $.

\item
Let $ \Norm{\cdot}{\X} $ denote the norm induced by $ \Inner{\cdot}{\cdot}{E} $. Obviously, $ \Norm{\cdot}{A_{0}} $ is the restriction of $ \Norm{\cdot}{\X} $ to $ A_{0} $.

\item
As $ E $ is an ideal of $ \RCP{G}{A}{\alpha} $, we obtain a left $ \RCP{G}{A}{\alpha} $-action on $ \X $.

\item
Next, let $ D_{0} \df \Span{\Set{\func{\mu}{a,b}}{a,b \in A_{0}}} $, which is a $ * $-subalgebra of $ \Mult{A} $, so $ D \df \Cl{D_{0}}{\Mult{A}} $ is a $ C^{*} $-subalgebra of $ \Mult{A} $.

\item
Define a right $ D_{0} $-action $ \bullet $ on $ A_{0} $ by right-multiplication, i.e.,
$$
\forall a \in A_{0}, ~ \forall d \in D_{0}: \quad
a \bullet d \df a d.
$$
Via $ \bullet $, each $ d \in D_{0} $ can be identified as a $ \Norm{\cdot}{A_{0}} $-bounded operator on $ A_{0} $ having norm $ \Norm{d}{\Mult{A}} $, which obviously extends to a $ \Norm{\cdot}{\X} $-bounded operator $ T_{d} $ on $ \X $ having norm $ \Norm{d}{\Mult{A}} $ also. Hence, for every $ d \in D_{0} $, it is true that $ T_{d} $ is an adjointable operator on $ \HModule{E}{\X}{} $ with adjoint $ T_{d^{*}} $.

\item
We thus have an isometric anti-homomorphism $ \Map{D_{0}}{\Adj{\HModule{E}{\X}{}}}{d}{T_{d}} $ that can be extended to an isometric anti-homomorphism from $ D $ to $ \Adj{\HModule{E}{\X}{}} $, which is actually an isometric anti-isomorphism from $ D $ to $ \Comp{\HModule{E}{\X}{}} $ because $ T_{\func{\mu}{a,b}} $ is a rank-$ 1 $ operator on $ \HModule{E}{\X}{} $ for every $ a,b \in A_{0} $:
$$
\forall c \in A_{0}: \quad
  \func{T_{\func{\mu}{a,b}}}{c}
= c ~ \func{\mu}{a,b}
= \func{\Psi}{a,b} \diamond b
= \Inner{c}{a}{E} \diamond b.
$$

\item
Therefore, $ D $ is an imprimitivity algebra of $ \HModule{E}{\X}{} $. This implies that there is a continuous extension of $ \func{\mu}{\cdot,\cdot} $ (which is a right $ D_{0} $-pre-inner product on $ A_{0} $) to a right $ D $-inner product $ \Inner{\cdot}{\cdot}{D} $ on $ \X $, with respect to which $ \HModule{}{\X}{D} $ is a full Hilbert $ D $-module.

\item
As $ \Inner{\cdot}{\cdot}{E} $ and $ \Inner{\cdot}{\cdot}{D} $ are compatible, $ \HModule{E}{\X}{D} $ is an imprimitivity $ \Pair{E}{D} $-bimodule. We then call $ D $ a \emph{generalized fixed-point algebra}, and it is Morita-Rieffel equivalent to the ideal $ E $ of $ \RCP{G}{A}{\alpha} $.
\end{itemize}

Evidently, $ D $ depends not only on $ \Trip{G}{A}{\alpha} $ but also on the dense $ * $-subalgebra $ A_{0} $. In order to fully reflect this dependence, we write $ \FixR{G}{A}{\alpha}{A_{0}} $ in place of $ D $. For precisely the same reasons, we write $ \ImpR{G}{A}{\alpha}{A_{0}} $ instead of $ \X $.

In \cite{Rieffel3}, Rieffel employed integrable group actions to provide yet another definition of properness, strictly weaker than \autoref{Rieffel-Properness} but still strong enough to build generalized fixed-point algebras. This definition had a drawback, for Ruy Exel showed in \cite{Exel} that the generalized fixed-point algebras, even in some cases where $ G $ is abelian, are too large to be equal to any ideal of $ \RCP{G}{A}{\alpha} $.

\subsection{Square-Integrable Representations of $ C^{*} $-Dynamical Systems} 

In \cite{Meyer2}, Ralf Meyer was able to construct generalized fixed-point algebras from minimal assumptions via square-integrable representations of $ C^{*} $-dynamical systems. To illustrate his idea, consider a \emph{Hilbert $ \Trip{G}{A}{\alpha} $-module}, i.e., a Hilbert $ A $-module $ \E $ endowed with a strongly continuous $ G $-action by linear isometries that is compatible with the right $ A $-action and the $ A $-inner product on $ \E $. If we denote the said $ G $-action on $ \E $ by $ \gam{\E} $, then what compatibility means is that
\begin{alignat*}{2}
\forall r \in G, ~ \forall a \in A, ~ \forall \zeta \in \E: & \quad &
\gammArg{\E}{r}{\zeta \bullet a}                            & = \gammArg{\E}{r}{\zeta} \bullet \alphArg{r}{a}, \\
\forall r \in G, ~ \forall \zeta,\eta \in \E:               & \quad &
\Inner{\gammArg{\E}{r}{\zeta}}{\gammArg{\E}{r}{\eta}}{\E}   & = \alphArg{r}{\Inner{\zeta}{\eta}{\E}}.
\end{alignat*}

Let $ \LL{2}{G,A} $ denote the Hilbert $ A $-module constructed in \autoref{Crossed Products} when $ \X = A_{A} $. The right $ A $-action $ \centerdot $ and the $ A $-valued inner product $ \Inner{\cdot}{\cdot}{\LL{2}{G,A}} $ thus satisfy the following properties:
\begin{itemize}
\item
$ \func{q}{\phi} \centerdot a = \func{q}{\Map{G}{A}{x}{\func{\phi}{x} ~ a}} $ for every $ a \in A $ and $ \phi \in \Cc{G,A} $.

\item
$ \D \Inner{\func{q}{\phi}}{\func{q}{\psi}}{\LL{2}{G,A}} = \Int{G}{\func{\phi}{x}^{*} \func{\psi}{x}}{x} $ for every $ \phi,\psi \in \Cc{G,A} $.
\end{itemize}
Define a strongly continuous $ G $-action $ \Gamma $ on $ \LL{2}{G,A} $ by linear isometries such that
$$
\forall r \in G, ~ \forall \phi \in \Cc{G,A}: \quad
\GammArg{r}{\func{q}{\phi}} = \func{q}{\Map{G}{A}{x}{\alphArg{x}{\func{\phi}{r^{-1} x}}}}.
$$
Then $ \LL{2}{G,A} $ is a Hilbert $ \Trip{G}{A}{\alpha} $-module.

For each $ \zeta \in \E $, define operators $ \Bra{\zeta}: \E \to \Cb{G,A} $ and $ \Ket{\zeta}: \Cc{G,A} \to \E $ by
\begin{alignat*}{2}
\forall \eta \in \E:       & \quad &
\BraArg{\zeta}{\eta}       & \df \Map{G}{A}{x}{\Inner{\gammArg{\E}{x}{\zeta}}{\eta}{\E}}, \\
\forall \phi \in \Cc{G,A}: & \quad &
\KetArg{\zeta}{\phi}       & \df \Int{G}{\gammArg{\E}{x}{\zeta} \bullet \func{\phi}{x}}{x},
\end{alignat*}
named the \emph{bra} and \emph{ket} of $ \zeta $ respectively. Then any $ \zeta \in \E $ is called \emph{square-integrable} if and only if for every $ \eta \in \E $ and every net $ \Seq{\varphi_{i}}{i \in I} $ in $ \Cc{G,\CC{0}{1}} $ converging uniformly to $ 1 $ on compact subsets of $ G $, the net $ \Seq{\func{q}{\varphi_{i} \BraArg{\zeta}{\eta}}}{i \in I} $ is Cauchy in $ \LL{2}{G,A} $, in which case the following are true:
\begin{itemize}
\item[(i)]
There is an operator $ \BBra{\zeta}: \E \to \LL{2}{G,A} $ defined by $ \D \BBraArg{\zeta}{\eta} \df \lim_{i \in I} \func{q}{\varphi_{i} \BraArg{\zeta}{\eta}} $ for every $ \eta \in \E $.

\item[(ii)]
There is an operator $ \KKet{\zeta}: \LL{2}{G,A} \to \E $ such that $ \KKetArg{\zeta}{\func{q}{\phi}} = \KetArg{\zeta}{\phi} $ for every $ \phi \in \Cc{G,A} $, whose adjoint is $ \BBra{\zeta} $.
\end{itemize}
Denote the set of square-integrable elements of $ \E $ by $ \Esi $, which is evidently a linear subspace of $ \E $. If $ \Esi $ is dense in $ \E $, then $ \E $ is called a \emph{square-integrable representation} of $ \Trip{G}{A}{\alpha} $.

Realizing $ \RCP{G}{A}{\alpha} $ as a $ C^{*} $-subalgebra of $ \AdjEq{\LL{2}{G,A}} $ (the set of equivariant adjointable operators on $ \E $), Meyer declared a linear subspace $ \R $ of $ \E $ to be \emph{relatively continuous} precisely when
$$
\R \subseteq \Esi \qquad \text{and} \qquad
\BBraKKet{\R}{\R} \subseteq \RCP{G}{A}{\alpha}.
$$
The concept of relative continuity was defined by Exel in \cite{Exel}, in the context of $ C^{*} $-dynamical systems where $ G $ is abelian. He defined it as a relation $ R $ on the set $ A_{\si} $ of square-integrable elements of $ A $ and showed that $ \Pair{a}{b} \in R \iff \BBraKKet{a}{b} \in \RCP{G}{A}{\alpha} $ for every $ a,b \in A_{\si} $.

From a relatively continuous subspace $ \R $, Meyer constructed a generalized fixed-point algebra as follows:
\begin{itemize}
\item
Let $ \Imp{\E}{\R} \df \Cl{\Span{\KKet{\R} \cup \Br{\KKet{\R} \circ \RCP{G}{A}{\alpha}}}}{\AdjEqPair{\LL{2}{G,A}}{\E}} $.

\item
Then $ \Imp{\E}{\R} $ is a Hilbert $ \RCP{G}{A}{\alpha} $-module, where the right $ \RCP{G}{A}{\alpha} $-action is defined by right-composition by elements of $ \RCP{G}{A}{\alpha} $, and the $ \RCP{G}{A}{\alpha} $-inner product $ \Inner{\cdot}{\cdot}{\Imp{\E}{\R}} $ by
$$
\forall P,Q \in \Imp{\E}{\R}: \quad
\Inner{P}{Q}{\Imp{\E}{\R}} \df P^{*} \circ Q.
$$

\item
$ \Imp{\E}{\R} $ is a \emph{full} Hilbert $ J $-module, where $ J \df \Cl{\Span{\Imp{\E}{\R}^{*} \circ \Imp{\E}{\R}}}{\RCP{G}{A}{\alpha}} $ is an ideal of $ \RCP{G}{A}{\alpha} $.

\item
The generalized fixed-point algebra is defined as $ \Fix{\E}{\R} \df \Cl{\Span{\Imp{\E}{\R} \circ \Imp{\E}{\R}^{*}}}{\AdjEq{\E}} $, which is isomorphic to $ \Comp{\Imp{\E}{\R}} $. Hence, $ \Fix{\E}{\R} $ is Morita-Rieffel equivalent to $ J $.
\end{itemize}

If we wish to fully reflect the dependence of $ \Imp{\E}{\R} $ and $ \Fix{\E}{\R} $ on $ \Trip{G}{A}{\alpha} $, we will utilize the notation $ \ImpFull{G}{A}{\alpha}{\E}{\R} $ and $ \FixFull{G}{A}{\alpha}{\E}{\R} $ respectively.

Meyer did not clarify the connection between his work and Rieffel's, but it is easily inferred. If $ \Trip{G}{A}{\alpha} $ is a $ C^{*} $-dynamical system, then $ A_{A} $ is a Hilbert $ \Trip{G}{A}{\alpha} $-module, where $ \gam{A} = \alpha $. If $ \Trip{G}{A}{\alpha} $ is Rieffel-proper, then any dense $ \alpha $-invariant $ * $-subalgebra $ A_{0} $ of $ A $ with the conditions in \autoref{Rieffel-Properness} is automatically relatively continuous (the fact that $ A_{0} \subseteq A_{\si} $ follows from Theorem 4.6 of \cite{Rieffel3}). Consequently, we can form $ \FixFull{G}{A}{\alpha}{A}{A_{0}} $, which is isomorphic to Rieffel's $ \FixR{G}{A}{\alpha}{A_{0}} $. One must be aware, however, that Rieffel's imprimitivity bimodule $ \ImpR{G}{A}{\alpha}{A_{0}} $ is the dual of $ \ImpFull{G}{A}{\alpha}{A}{A_{0}} $, because $ \FixR{G}{A}{\alpha}{A_{0}} $ acts on the left of $ \ImpR{G}{A}{\alpha}{A_{0}} $ whereas $ \FixFull{G}{A}{\alpha}{A}{A_{0}} $ acts on the right of $ \ImpFull{G}{A}{\alpha}{A}{A_{0}} $.

In \cite{Mingo|Phillips}, J. Mingo and W. Phillips proved for a countably generated Hilbert $ \Trip{G}{A}{\alpha} $-module $ \E $ that there exists an equivariant isomorphism $ \E \oplus \LL{2}{G,A}^{\infty} \cong \LL{2}{G,A}^{\infty} $ under certain conditions, such as when $ G $ is compact. This is the equivariant version of Kasparov's Stabilization Theorem. In \cite{Meyer1}, Meyer deduced the square-integrability of $ \E $ to be a necessary and sufficient condition; right at the core of his argument is the fact that his bra-ket operators are equivariant.

It seems natural to replace $ C^{*} $-dynamical systems and Hilbert $ C^{*} $-modules in Meyer's work by twisted ones and see which results can be generalized. Twisted $ C^{*} $-dynamical systems have been studied extensively (\cite{Busby|Smith,Packer|Raeburn}), but twisted Hilbert $ C^{*} $-modules appear to be a new concept. As such, one of the challenges that I faced was to propose a correct definition of a twisted Hilbert $ C^{*} $-module that would allow Meyer's ideas to be generalized to twisted $ C^{*} $-dynamical systems.


\section{Twisted $ C^{*} $-Dynamical Systems and Twisted Hilbert $ C^{*} $-Modules} 

\subsection{Twisted $ C^{*} $-Dynamical Systems} 


\begin{Def}[\cite{Busby|Smith,Packer|Raeburn}] \label{Twisted C*-Dynamical Systems}
A \emph{twisted $ C^{*} $-dynamical system} is a quadruple $ \Quad{G}{A}{\alpha}{\omega} $, where:
\begin{enumerate}
\item[(1)]
$ G $ is an l.c.H. group and $ A $ a $ C^{*} $-algebra.

\item[(2)]
$ \alpha: G \to \Aut{A} $ is a strongly continuous map, i.e., the function $ \Map{G}{A}{x}{\alphArg{x}{a}} $ is continuous for every $ a \in A $.

\item[(3)]
$ \omega: G \times G \to \UMult{A} $ is a strictly continuous map, i.e., the functions
$$
\Map{G \times G}{A}{\Pair{x}{y}}{\om{x}{y} ~ a} \qquad \text{and} \qquad
\Map{G \times G}{A}{\Pair{x}{y}}{a ~ \om{x}{y}}
$$
are continuous for every $ a \in A $. We call $ \omega $ an \emph{$ A $-multiplier on $ G $}.

\item[(4)]
$ \alph{e} = \Id_{A} $, and $ \om{e}{r} = 1_{\Mult{A}} = \om{r}{e} $ for every $ r \in G $.

\item[(5)]
$ \alphExt{r} \circ \alphExt{s} = \Ad{\om{r}{s}} \circ \alphExt{r s} $ for every $ r,s \in G $, i.e.,
$$
\forall m \in \Mult{A}: \quad
\alphArgExt{r}{\alphArgExt{s}{m}} = \om{r}{s} ~ \alphArgExt{r s}{m} ~ \om{r}{s}^{*}.
$$

\item[(6)]
$ \alphArgExt{r}{\om{s}{t}} ~ \om{r}{s t} = \om{r}{s} ~ \om{r s}{t} $ for every $ r,s,t \in G $.
\end{enumerate}
\end{Def}


\textbf{From now on, $ \Quad{G}{A}{\alpha}{\omega} $ denotes an arbitrary twisted $ C^{*} $-dynamical system.}


\begin{Eg} \label{C*-Dynamical Systems Are Twisted C*-Dynamical Systems with a Trivial Multiplier}
Let $ \Trip{H}{B}{\beta} $ be a $ C^{*} $-dynamical system. If $ \omega $ is the trivial $ B $-multiplier on $ H $, i.e., $ \om{r}{s} = 1_{\Mult{B}} $ for every $ r,s \in G $, then $ \Quad{H}{B}{\beta}{\omega} $ is a twisted $ C^{*} $-dynamical system.
\end{Eg}


\begin{Eg} \label{Twisted C*-Dynamical Systems Associated with Higher-Dimensional Non-Commutative Tori}
Let $ d \in \N{} $. Let $ \Theta $ be any skew-symmetric bilinear form on $ \RR{}{d} $. Then a well-known example of a twisted $ C^{*} $-dynamical system (used to define $ d $-dimensional non-commutative tori) is $ \Quad{\Z{d}}{\C{}}{\tr}{\omega_{\Theta}} $, where $ \tr $ denotes the trivial action of $ \Z{d} $ on $ \C{} $, and $ \omega_{\Theta}: \Z{d} \times \Z{d} \to \T{} $ the normalized $ 2 $-cocycle corresponding to $ \Theta $, i.e., $ \func{\omega_{\Theta}}{\mathbf{m},\mathbf{n}} = e^{\pi i \func{\Theta}{\mathbf{m},\mathbf{n}}} $ for every $ \Pair{\mathbf{m}}{\mathbf{n}} \in \Z{d} \times \Z{d} $.
\end{Eg}


\subsection{Twisted Hilbert $ C^{*} $-Modules} 


\begin{Def} \label{Twisted Hilbert C*-Modules}
A \emph{Hilbert $ \Quad{G}{A}{\alpha}{\omega} $-module} is a Hilbert $ A $-module $ \E $ with a \emph{strongly continuous} map $ \gamma: G \to \Isom{\E} $ (the set of linear isometries on $ \E $), called a \emph{twisted action}, having the following properties:
\begin{enumerate}
\item[(1)]
$ \gamma_{e} = \Id_{\E} $.

\item[(2)]
$ \func{\gamma_{r}}{\zeta \bullet a} = \func{\gamma_{r}}{\zeta} \bullet \alphArg{r}{a} $ for every $ r \in G $, $ a \in A $ and $ \zeta \in \E $.

\item[(3)]
$ \Inner{\func{\gamma_{r}}{\zeta}}{\func{\gamma_{r}}{\eta}}{\E} = \alphArg{r}{\Inner{\zeta}{\eta}{\E}} $ for every $ r \in G $ and $ \zeta,\eta \in \E $.

\item[(4)]
$ \func{\gamma_{r}}{\func{\gamma_{s}}{\zeta}} = \func{\gamma_{r s}}{\zeta} \bullet \om{r}{s}^{*} $ for every $ r,s \in G $ and $ \zeta \in \E $.
\end{enumerate}
For better clarity, we denote the twisted action on $ \E $ by $ \gam{\E} $. If $ \Quad{G}{A}{\alpha}{\omega} $ is clear from the context, then we simply call $ \E $ a \emph{twisted Hilbert $ C^{*} $-module}.
\end{Def}


\begin{Eg} \label{A Basic Twisted Hilbert C*-Module}
Recall the Hilbert $ A $-module $ \LL{2}{G,A} $ defined earlier. To construct a twisted $ G $-action on $ \LL{2}{G,A} $, first observe for every $ r \in G $ and $ \phi \in \Cc{G,A} $ that
\begin{align*}
  & ~ \Norm{\func{q}{\Map{G}{A}{x}{\om{r}{r^{-1} x}^{*} \alphArg{r}{\func{\phi}{r^{-1} x}}}}}{\LL{2}{G,A}} \\
= & ~ \NNorm{\Map{G}{A}{x}{\om{r}{r^{-1} x}^{*} \alphArg{r}{\func{\phi}{r^{-1} x}}}} \\
= & ~ \Norm{
           \Int{G}{
                  \SqBr{\om{r}{r^{-1} x}^{*} \alphArg{r}{\func{\phi}{r^{-1} x}}}^{*}
                  \SqBr{\om{r}{r^{-1} x}^{*} \alphArg{r}{\func{\phi}{r^{-1} x}}}
                  }{x}
           }{A}^{\frac{1}{2}} \\
= & ~ \Norm{
           \Int{G}{\alphArg{r}{\func{\phi}{r^{-1} x}}^{*} \om{r}{r^{-1} x} ~ \om{r}{r^{-1} x}^{*} \alphArg{r}{\func{\phi}{r^{-1} x}}}{x}
           }{A}^{\frac{1}{2}} \\
= & ~ \Norm{\Int{G}{\alphArg{r}{\func{\phi}{r^{-1} x}}^{*} \alphArg{r}{\func{\phi}{r^{-1} x}}}{x}}{A}^{\frac{1}{2}} \\
= & ~ \Norm{\Int{G}{\alphArg{r}{\func{\phi}{r^{-1} x}^{*}} ~ \alphArg{r}{\func{\phi}{r^{-1} x}}}{x}}{A}^{\frac{1}{2}} \\
= & ~ \Norm{\Int{G}{\alphArg{r}{\func{\phi}{r^{-1} x}^{*} \func{\phi}{r^{-1} x}}}{x}}{A}^{\frac{1}{2}} \\
= & ~ \Norm{\alphArg{r}{\Int{G}{\func{\phi}{r^{-1} x}^{*} \func{\phi}{r^{-1} x}}{x}}}{A}^{\frac{1}{2}} \qquad
      \Br{\text{As $ \alph{r} $ is continuous.}} \\
= & ~ \Norm{\Int{G}{\func{\phi}{r^{-1} x}^{*} \func{\phi}{r^{-1} x}}{x}}{A}^{\frac{1}{2}} \qquad \Br{\text{As $ \alph{r} $ is isometric.}} \\
= & ~ \Norm{\Int{G}{\func{\phi}{x}^{*} \func{\phi}{x}}{x}}{A}^{\frac{1}{2}} \qquad \Br{\text{By the change of variables $ x \mapsto r x $.}} \\
= & ~ \NNorm{\phi} \\
= & ~ \Norm{\func{q}{\phi}}{\LL{2}{G,A}}.
\end{align*}
We can thus define a map $ \Gamma: G \to \Isom{\LL{2}{G,A}} $ by
$$
\forall r \in G, ~ \forall \Phi \in \LL{2}{G,A}: \quad
\GammArg{r}{\Phi} \df \lim_{n \to \infty} \func{q}{\Map{G}{A}{x}{\om{r}{r^{-1} x}^{*} \alphArg{r}{\func{\phi_{n}}{r^{-1} x}}}},
$$
where $ \Seq{\phi_{n}}{n \in \N{}} $ is \emph{any} sequence in $ \Cc{G,A} $ with $ \D \lim_{n \to \infty} \func{q}{\phi_{n}} = \Phi $.

We now check the strong continuity of $ \Gamma $. Let $ \epsilon > 0 $, $ r \in G $, $ \phi \in \Cc{G,A} \setminus \SSet{\bm{0}} $ and $ S \df \Supp{\phi} $. Fix a compact subset $ K $ of $ G $ containing $ r $ in its interior. Our aim then is to obtain the limit
$$
\lim_{s \to r} \Norm{\GammArg{s}{\func{q}{\phi}} - \GammArg{r}{\func{q}{\phi}}}{\LL{2}{G,A}} = 0.
$$
As the function
$$
\Map{G \times G}{A}{\Pair{y}{s}}{\om{s}{s^{-1} y}^{*} \alphArg{s}{\func{\phi}{s^{-1} y}}}
$$
is continuous, we can find $ K S $-indexed sequences $ \Seq{V_{x}}{x \in K S} $ and $ \Seq{W_{x}}{x \in K S} $ of subsets of $ G $ having the following properties for every $ x \in K S $:
\begin{itemize}
\item
$ V_{x} $ is the intersection of $ K S $ with an open neighborhood of $ x $.

\item
$ W_{x} $ is the intersection of $ K^{\circ} $ with an open neighborhood of $ r $.

\item
$ \Norm{\om{s}{s^{-1} y}^{*} \alphArg{s}{\func{\phi}{s^{-1} y}} - \om{r}{r^{-1} x}^{*} \alphArg{r}{\func{\phi}{r^{-1} x}}}{A} < \dfrac{\epsilon}{2 \sqrt{\func{\mu}{K S}}} $ for every $ \Pair{y}{s} \in V_{x} \times W_{x} $, whence
$$
\forall \Pair{y}{s} \in V_{x} \times W_{x}: \quad
  \Norm{\om{s}{s^{-1} y}^{*} \alphArg{s}{\func{\phi}{s^{-1} y}} - \om{r}{r^{-1} y}^{*} \alphArg{r}{\func{\phi}{r^{-1} y}}}{A}
< \dfrac{\epsilon}{\sqrt{\func{\mu}{K S}}}.
$$
\end{itemize}
By the compactness of $ K S $, there exist points $ x_{1},\ldots,x_{n} $ that satisfy $ \D K S = \bigcup_{k = 1}^{n} V_{x_{k}} $. Pick any open neighborhood $ N $ of $ r $ contained within $ \D \bigcap_{k = 1}^{n} W_{x_{k}} $, and let $ \Pair{x}{s} \in K S \times N $. Find a $ k \in \SSet{1,\ldots,n} $ such that $ x \in V_{x_{k}} $. As $ s \in W_{x_{k}} $, we have
$$
  \Norm{\om{s}{s^{-1} x}^{*} \alphArg{s}{\func{\phi}{s^{-1} x}} - \om{r}{r^{-1} x}^{*} \alphArg{r}{\func{\phi}{r^{-1} x}}}{A}
< \frac{\epsilon}{\sqrt{\func{\mu}{K S}}}.
$$
We chose $ x $ arbitrarily, so
\begin{align*}
     & ~ \Norm{\GammArg{s}{\func{q}{\phi}} - \GammArg{r}{\func{q}{\phi}}}{\LL{2}{G,A}} \\
=    & ~ \NNorm{
               \Map{G}{A}{x}{
                            \om{s}{s^{-1} x}^{*} \alphArg{s}{\func{\phi}{s^{-1} x}} -
                            \om{r}{r^{-1} x}^{*} \alphArg{r}{\func{\phi}{r^{-1} x}}
                            }
               } \\
\leq & ~ \Norm{
              \Map{G}{A}{x}{
                           \om{s}{s^{-1} x}^{*} \alphArg{s}{\func{\phi}{s^{-1} x}} -
                           \om{r}{r^{-1} x}^{*} \alphArg{r}{\func{\phi}{r^{-1} x}}
                           }
              }{2} \\
=    & ~ \SqBr{
              \Int{G}{
                     \Norm{
                          \om{s}{s^{-1} x}^{*} \alphArg{s}{\func{\phi}{s^{-1} x}} - \om{r}{r^{-1} x}^{*} \alphArg{r}{\func{\phi}{r^{-1} x}}
                          }{A}^{2}
                     }{x}
              }^{\frac{1}{2}} \\
=    & ~ \SqBr{
              \Int{K S}{
                       \Norm{
                            \om{s}{s^{-1} x}^{*} \alphArg{s}{\func{\phi}{s^{-1} x}} - \om{r}{r^{-1} x}^{*} \alphArg{r}{\func{\phi}{r^{-1} x}}
                            }{A}^{2}
                       }{x}
              }^{\frac{1}{2}} \\
     & ~ \Br{\text{As the integrand vanishes outside of $ K S $.}} \\
<    & ~ \SqBr{\Int{K S}{\frac{\epsilon^{2}}{\func{\mu}{K S}}}{x}}^{\frac{1}{2}} \\
=    & ~ \SqBr{\frac{\epsilon^{2}}{\func{\mu}{K S}} \cdot \func{\mu}{K S}}^{\frac{1}{2}} \\
=    & ~ \epsilon.
\end{align*}
As $ \epsilon $ and $ \phi $ are arbitrary, we get $ \D \lim_{s \to r} \Norm{\GammArg{s}{\func{q}{\phi}} - \GammArg{r}{\func{q}{\phi}}}{\LL{2}{G,A}} = 0 $ for every $ \phi \in \Cc{G,A} $.

Let $ \Phi \in \LL{2}{G,A} $. Let $ \epsilon > 0 $ once more, and pick a $ \phi \in \Cc{G,A} $ so that $ \Norm{\Phi - \func{q}{\phi}}{\LL{2}{G,A}} < \dfrac{\epsilon}{3} $. By the argument above, there is an open neighborhood $ N $ of $ r $ such that
$$
\forall s \in N: \quad
\Norm{\GammArg{s}{\func{q}{\phi}} - \GammArg{r}{\func{q}{\phi}}}{\LL{2}{G,A}} < \frac{\epsilon}{3}, \quad \text{from which it follows that}
$$
\begin{align*}
     & ~ \Norm{\GammArg{s}{\Phi} - \GammArg{r}{\Phi}}{\LL{2}{G,A}} \\
\leq & ~ \Norm{\GammArg{s}{\Phi} - \GammArg{s}{\func{q}{\phi}}}{\LL{2}{G,A}} ~ + \\
     & ~ \Norm{\GammArg{s}{\func{q}{\phi}} - \GammArg{r}{\func{q}{\phi}}}{\LL{2}{G,A}} ~ + \\
     & ~ \Norm{\GammArg{r}{\func{q}{\phi}} - \GammArg{r}{\Phi}}{\LL{2}{G,A}} \\
=    & ~ \Norm{\Phi - \func{q}{\phi}}{\LL{2}{G,A}} +
         \Norm{\GammArg{s}{\func{q}{\phi}} - \GammArg{r}{\func{q}{\phi}}}{\LL{2}{G,A}} +
         \Norm{\func{q}{\phi} - \Phi}{\LL{2}{G,A}} \\
<    & ~ \frac{\epsilon}{3} + \frac{\epsilon}{3} + \frac{\epsilon}{3} \\
=    & ~ \epsilon.
\end{align*}
As $ \epsilon $ is arbitrary, we get $ \D \lim_{s \to r} \Norm{\GammArg{s}{\Phi} - \GammArg{r}{\Phi}}{\LL{2}{G,A}} = 0 $. Then as $ r $ and $ \Phi $ are arbitrary, we conclude that $ \Gamma $ is strongly continuous.

This type of compactness argument will be a recurring theme throughout this work.

To show that $ \Gamma $ is a twisted action, the four conditions in \autoref{Twisted Hilbert C*-Modules} must be verified:
\begin{enumerate}
\item[(1)]
Trivial.

\item[(2)]
For every $ r \in G $, $ a \in A $ and $ \phi \in \Cc{G,A} $, we have
\begin{align*}
    \GammArg{r}{\func{q}{\phi} \centerdot a}
& = \GammArg{r}{\func{q}{\Map{G}{A}{x}{\func{\phi}{x} ~ a}}} \\
& = \func{q}{\Map{G}{A}{x}{\om{r}{r^{-1} x}^{*} \alphArg{r}{\func{\phi}{r^{-1} x} ~ a}}} \\
& = \func{q}{\Map{G}{A}{x}{\om{r}{r^{-1} x}^{*} \alphArg{r}{\func{\phi}{r^{-1} x}} ~ \alphArg{r}{a}}} \\
& = \func{q}{\Map{G}{A}{x}{\om{r}{r^{-1} x}^{*} \alphArg{r}{\func{\phi}{r^{-1} x}}}} \centerdot \alphArg{r}{a} \\
& = \GammArg{r}{\func{q}{\phi}} \centerdot \alphArg{r}{a},
\end{align*}
so by continuity, $ \GammArg{r}{\Phi \centerdot a} = \GammArg{r}{\Phi} \centerdot \alphArg{r}{a} $ for every $ \Phi \in \LL{2}{G,A} $.

\item[(3)]
For every $ r \in G $ and $ \phi,\psi \in \Cc{G,A} $, we have
\begin{align*}
    \Inner{\GammArg{r}{\func{q}{\phi}}}{\GammArg{r}{\func{q}{\psi}}}{\LL{2}{G,A}}
& = \Int{G}{
           \SqBr{\om{r}{r^{-1} x}^{*} \alphArg{r}{\func{\phi}{r^{-1} x}}}^{*} \SqBr{\om{r}{r^{-1} x}^{*} \alphArg{r}{\func{\psi}{r^{-1} x}}}
           }{x} \\
& = \Int{G}{\alphArg{r}{\func{\phi}{r^{-1} x}}^{*} \om{r}{r^{-1} x} ~ \om{r}{r^{-1} x}^{*} \alphArg{r}{\func{\psi}{r^{-1} x}}}{x} \\
& = \Int{G}{\alphArg{r}{\func{\phi}{r^{-1} x}}^{*} \alphArg{r}{\func{\psi}{r^{-1} x}}}{x} \\
& = \Int{G}{\alphArg{r}{\func{\phi}{r^{-1} x}^{*}} ~ \alphArg{r}{\func{\psi}{r^{-1} x}}}{x} \\
& = \Int{G}{\alphArg{r}{\func{\phi}{r^{-1} x}^{*} \func{\psi}{r^{-1} x}}}{x} \\
& = \alphArg{r}{\Int{G}{\func{\phi}{r^{-1} x}^{*} \func{\psi}{r^{-1} x}}{x}} \\
& = \alphArg{r}{\Int{G}{\func{\phi}{x}^{*} \func{\psi}{x}}{x}} \qquad \Br{\text{By the change of variables $ x \mapsto r x $.}} \\
& = \alphArg{r}{\Inner{\func{q}{\phi}}{\func{q}{\psi}}{\LL{2}{G,A}}},
\end{align*}
so by continuity, $ \Inner{\GammArg{r}{\Phi}}{\GammArg{r}{\Psi}}{\LL{2}{G,A}} = \alphArg{r}{\Inner{\Phi}{\Psi}{\LL{2}{G,A}}} $ for every $ \Phi,\Psi \in \LL{2}{G,A} $.

\item[(4)]
Finally, for every $ r,s \in G $ and $ \phi \in \Cc{G,A} $, we have
\begin{align*}
  & ~ \GammArg{r}{\GammArg{s}{\func{q}{\phi}}} \\
= & ~ \GammArg{r}{\func{q}{\Map{G}{A}{x}{\om{s}{s^{-1} x}^{*} \alphArg{s}{\func{\phi}{s^{-1} x}}}}} \\
= & ~ \func{q}{\Map{G}{A}{x}{\om{r}{r^{-1} x}^{*} \alphArg{r}{\om{s}{s^{-1} r^{-1} x}^{*} \alphArg{s}{\func{\phi}{s^{-1} r^{-1} x}}}}} \\
= & ~ \func{q}{
              \Map{G}{A}{x}{
                           \om{r}{r^{-1} x}^{*}
                           \alphArgExt{r}{\om{s}{s^{-1} r^{-1} x}^{*}} ~
                           \alphArg{r}{\alphArg{s}{\func{\phi}{s^{-1} r^{-1} x}}}
                           }
              } \\
= & ~ \func{q}{
              \Map{G}{A}{x}{
                           \om{r}{r^{-1} x}^{*}
                           \alphArgExt{r}{\om{s}{s^{-1} r^{-1} x}}^{*}
                           \alphArg{r}{\alphArg{s}{\func{\phi}{s^{-1} r^{-1} x}}}
                           }
              } \\
= & ~ \func{q}{
              \Map{G}{A}{x}{
                           \om{r}{r^{-1} x}^{*}
                           \alphArgExt{r}{\om{s}{s^{-1} r^{-1} x}}^{*}
                           \om{r}{s} ~
                           \alphArg{r s}{\func{\phi}{s^{-1} r^{-1} x}} ~
                           \om{r}{s}^{*}
                           }
              } \\
= & ~ \func{q}{
              \Map{G}{A}{x}{
                           \SqBr{\alphArgExt{r}{\om{s}{s^{-1} r^{-1} x}} ~ \om{r}{r^{-1} x}}^{*}
                           \om{r}{s} ~
                           \alphArg{r s}{\func{\phi}{s^{-1} r^{-1} x}} ~
                           \om{r}{s}^{*}
                           }
              } \\
= & ~ \func{q}{
              \Map{G}{A}{x}{
                           \SqBr{\om{r}{s} ~ \om{r s}{s^{-1} r^{-1} x}}^{*}
                           \om{r}{s} ~
                           \alphArg{r s}{\func{\phi}{s^{-1} r^{-1} x}} ~
                           \om{r}{s}^{*}
                           }
              } \\
= & ~ \func{q}{
              \Map{G}{A}{x}{
                           \om{r s}{s^{-1} r^{-1} x}^{*} \om{r}{s}^{*} \om{r}{s} ~ \alphArg{r s}{\func{\phi}{s^{-1} r^{-1} x}} ~ \om{r}{s}^{*}
                           }
              } \\
= & ~ \func{q}{\Map{G}{A}{x}{\om{r s}{s^{-1} r^{-1} x}^{*} \alphArg{r s}{\func{\phi}{s^{-1} r^{-1} x}} ~ \om{r}{s}^{*}}} \\
= & ~ \func{q}{\Map{G}{A}{x}{\om{r s}{s^{-1} r^{-1} x}^{*} \alphArg{r s}{\func{\phi}{s^{-1} r^{-1} x}}}} \centerdot \om{r}{s}^{*} \\
= & ~ \func{q}{\Map{G}{A}{x}{\om{r s}{\Br{r s}^{-1} x}^{*} \alphArg{r s}{\func{\phi}{\Br{r s}^{-1} x}}}} \centerdot \om{r}{s}^{*} \\
= & ~ \GammArg{r s}{\func{q}{\phi}} \centerdot \om{r}{s}^{*},
\end{align*}
so by continuity, $ \GammArg{r}{\GammArg{s}{\Phi}} = \GammArg{r s}{\Phi} \centerdot \om{r}{s}^{*} $ for every $ \Phi \in \LL{2}{G,A} $.
\end{enumerate}
\end{Eg}



\begin{Rmk} \label{Non-Equivalent Norms}
Observe that $ \NNorm{\cdot} \leq \Norm{\cdot}{2} $, where $ \Norm{\cdot}{2} $ denotes the $ L^{2} $-norm on $ \Cc{G,A} $, i.e.,
$$
\forall f \in \Cc{G,A}: \quad
\Norm{f}{2} = \Br{\Int{G}{\Norm{\func{f}{x}}{A}^{2}}{x}}^{\frac{1}{2}}.
$$
Despite this, unless $ A = \C{} $ and/or $ G $ is finite, $ \NNorm{\cdot} $ and $ \Norm{\cdot}{2} $ are generally not equivalent. To see why, let $ G = \Z{} $ and $ A = \Cont{0}{\Z{}} $. Let $ \Seq{f_{n}}{n \in \N{}} $ be a sequence in $ A $ whose members have disjoint support, have sup-norm $ 1 $, and are non-negative. Next, define a sequence $ \Seq{F_{n}}{n \in \N{}} $ in $ \Cc{G,A} $ by
$$
\forall \Pair{n}{k} \in \N{} \times \Z{}: \quad
\func{F_{n}}{k} \df
\begin{cases}
f_{k}  & \text{if $ 1 \leq k \leq n $}; \\
\bm{0} & \text{otherwise}.
\end{cases}
$$
Then
$$
\forall n \in \N{}: \quad
  \frac{\NNorm{F_{n}}}{\Norm{F_{n}}{2}}
= \frac{\D \Norm{\sum_{k = 1}^{n} f_{k}^{2}}{\infty}^{\frac{1}{2}}}{\D \Br{\sum_{k = 1}^{n} \Norm{f_{k}}{\infty}^{2}}^{\frac{1}{2}}}
= \frac{1}{\sqrt{n}},
$$
which results in $ \D \lim_{n \to \infty} \frac{\NNorm{F_{n}}}{\Norm{F_{n}}{2}} = 0 $. Therefore, $ \NNorm{\cdot} $ and $ \Norm{\cdot}{2} $ are not equivalent. This tells us that we \emph{should not} identify $ \LL{2}{G,A} $ with the Banach space of all (equivalence classes of) square-integrable strongly measurable $ A $-valued functions on $ G $.
\end{Rmk}



\begin{Eg} \label{The Direct Sum of Twisted Hilbert C*-Modules}
If $ \Seq{\E_{n}}{n \in \N{}{}} $ is a sequence of Hilbert $ \Quad{G}{A}{\alpha}{\omega} $-modules, we can define a \emph{direct-sum Hilbert $ \Quad{G}{A}{\alpha}{\omega} $-module} $ \D \bigoplus_{n = 1}^{\infty} \E_{n} $ in the following manner:
\begin{itemize}
\item
The underlying vector space is the set of sequences $ \D \Seq{\zeta_{n}}{n \in \N{}{}} \in \prod_{n \in \N{}{}} \E_{n} $ such that $ \sum_{n = 1}^{\infty} \Inner{\zeta_{n}}{\zeta_{n}}{\E_{n}} $ converges. Such sequences are either denoted by a bold Greek letter or written as a formal sum. For example, $ \Seq{\zeta_{n}}{n \in \N{}{}} $ is either denoted by $ \bm{\zeta} $ or expanded as $ \D \sum_{n = 1}^{\infty} \zeta_{n} \cdot \mathbf{e}_{n} $.

\item
Define the right $ A $-action by $ \D \bm{\zeta} \bullet a \df \sum_{n = 1}^{\infty} \Br{\zeta_{n} \bullet a} \cdot \mathbf{e}_{n} $ for every $ a \in A $ and $ \D \bm{\zeta} \in \bigoplus_{n = 1}^{\infty} \E_{n} $.

\item
Define the $ A $-valued inner product by $ \D \Inner{\bm{\zeta}}{\bm{\eta}}{\bigoplus_{n = 1}^{\infty} \E_{n}} \df \sum_{n = 1}^{\infty} \Inner{\zeta_{n}}{\eta_{n}}{\E_{n}} $ for every $ \D \bm{\zeta},\bm{\eta} \in \bigoplus_{n = 1}^{\infty} \E_{n} $.

\item
Define the twisted $ G $-action by $ \D \func{\Br{\gam{\bigoplus_{n = 1}^{\infty} \E_{n}}}_{r}}{\bm{\zeta}} \df \sum_{n = 1}^{\infty} \gammArg{\E_{n}}{r}{\zeta_{n}} \cdot \mathbf{e}_{n} $ for every $ r \in G $ and $ \D \bm{\zeta} \in \bigoplus_{n = 1}^{\infty} \E_{n} $.
\end{itemize}
For every Hilbert $ \Quad{G}{A}{\alpha}{\omega} $-module $ \E $, let $ \D \E^{\infty} \df \bigoplus_{n = 1}^{\infty} \E $. Also, let $ \Gamma^{\infty} $ denote the twisted $ G $-action on $ \LL{2}{G,A}^{\infty} $.
\end{Eg}


\subsection{The Category of Twisted Hilbert $ C^{*} $-Modules} 


\begin{Def} \label{The Category of Twisted Hilbert C*-Modules}
Write $ \Hilb{G}{A}{\alpha}{\omega} $ for the category of Hilbert $ \Quad{G}{A}{\alpha}{\omega} $-modules. If $ \E $ and $ \F $ are Hilbert $ \Quad{G}{A}{\alpha}{\omega} $-modules, then a morphism from $ \E $ to $ \F $ is an adjointable operator $ T: \E \to \F $ such that $ \func{T}{\gammArg{\E}{r}{\zeta}} = \gammArg{\F}{r}{\func{T}{\zeta}} $ for every $ r \in G $ and $ \zeta \in \E $ (we say that $ T $ is \emph{twisted-equivariant}). Denote the set of morphisms from $ \E $ to $ \F $ by $ \AdjEqPair{\E}{\F} $, and write $ \AdjEq{\E} $ for $ \AdjEqPair{\E}{\E} $.
\end{Def}


It is important to know that $ \Hilb{G}{A}{\alpha}{\omega} $-morphisms are closed under the operator-adjoint.


\begin{Lem} \label{Morphisms of Twisted Hilbert C*-Modules Are Closed Under the Operator-Adjoint}
If $ T: \E \to \F $ is a $ \Hilb{G}{A}{\alpha}{\omega} $-morphism, then so is $ T^{*}: \F \to \E $.
\end{Lem}

\begin{proof}
Let $ r \in G $ and $ \zeta \in \F $. Then for every $ \eta \in \E $,
\begin{align*}
    \Inner{\func{T^{*}}{\gammArg{\F}{r}{\zeta}}}{\eta}{\E}
& = \Inner{\gammArg{\F}{r}{\zeta}}{\func{T}{\eta}}{\F} \\
& = \Inner{\gammArg{\F}{r}{\zeta}}{\gammArg{\F}{r}{\gammArg{\F}{r^{-1}}{\func{T}{\eta}}} \bullet \om{r}{r^{-1}}}{\F} \qquad
    \Br{\text{By (4) of \autoref{Twisted Hilbert C*-Modules}.}} \\
& = \Inner{\gammArg{\F}{r}{\zeta}}{\gammArg{\F}{r}{\gammArg{\F}{r^{-1}}{\func{T}{\eta}}}}{\F} ~ \om{r}{r^{-1}} \\
& = \alphArg{r}{\Inner{\zeta}{\gammArg{\F}{r^{-1}}{\func{T}{\eta}}}{\F}} ~ \om{r}{r^{-1}} \qquad
    \Br{\text{By (3) of \autoref{Twisted Hilbert C*-Modules}.}} \\
& = \alphArg{r}{\Inner{\zeta}{\func{T}{\gammArg{\E}{r^{-1}}{\eta}}}{\F}} ~ \om{r}{r^{-1}} \qquad
    \Br{\text{As $ T $ is twisted-equivariant.}} \\
& = \alphArg{r}{\Inner{\func{T^{*}}{\zeta}}{\gammArg{\E}{r^{-1}}{\eta}}{\E}} ~ \om{r}{r^{-1}} \\
& = \Inner{\gammArg{\E}{r}{\func{T^{*}}{\zeta}}}{\gammArg{\E}{r}{\gammArg{\E}{r^{-1}}{\eta}}}{\E} ~ \om{r}{r^{-1}} \qquad
    \Br{\text{By (3) of \autoref{Twisted Hilbert C*-Modules} again.}} \\
& = \Inner{\gammArg{\E}{r}{\func{T^{*}}{\zeta}}}{\eta \bullet \om{r}{r^{-1}}^{*}}{\E} ~ \om{r}{r^{-1}} \qquad
    \Br{\text{By (4) of \autoref{Twisted Hilbert C*-Modules} again.}} \\
& = \Inner{\gammArg{\E}{r}{\func{T^{*}}{\zeta}}}{\eta \bullet \om{r}{r^{-1}}^{*} \om{r}{r^{-1}}}{\E} \\
& = \Inner{\gammArg{\E}{r}{\func{T^{*}}{\zeta}}}{\eta}{\E}.
\end{align*}
Therefore, $ \func{T^{*}}{\gammArg{\F}{r}{\zeta}} = \gammArg{\E}{r}{\func{T^{*}}{\zeta}} $, and as $ r $ and $ \zeta $ are arbitrary, we are done.
\end{proof}



\begin{Rmk} \label{The Set of Twisted-Equivariant Adjointable Operators on a Twisted Hilbert C*-Module Is a C*-Algebra}
By \autoref{Morphisms of Twisted Hilbert C*-Modules Are Closed Under the Operator-Adjoint}, $ \AdjEq{\E} $ is a $ C^{*} $-algebra for every Hilbert $ \Quad{G}{A}{\alpha}{\omega} $-module $ \E $. We will later define the reduced twisted crossed product for $ \Quad{G}{A}{\alpha}{\omega} $ as a $ C^{*} $-subalgebra of $ \AdjEq{\LL{2}{G,A}} $.
\end{Rmk}



\begin{Eg} \label{Invariant Orthogonal Summands of a Twisted Hilbert C*-Module}
Let $ \E $ be a Hilbert $ \Quad{G}{A}{\alpha}{\omega} $-module and $ \F $ a $ \gam{\E} $-invariant orthogonal summand of $ \E $, i.e., $ \gammArg{\E}{r}{\zeta} \in \F $ for every $ r \in G $ and $ \zeta \in \F $. Then $ \F $ is a Hilbert $ \Quad{G}{A}{\alpha}{\omega} $-submodule of $ \E $. To see that the orthogonal complement $ \F^{\perp} $ of $ \F $ is $ \gam{\E} $-invariant, suppose that $ \eta \in \F^{\perp} $. Then
\begin{align*}
\forall r \in G, ~ \forall \zeta \in \F: \quad
      \Inner{\zeta}{\gammArg{\E}{r}{\eta}}{\E}
& =   \Inner{\gammArg{\E}{r}{\gammArg{\E}{r^{-1}}{\zeta}} \bullet \om{r}{r^{-1}}}{\gammArg{\E}{r}{\eta}}{\E} \\
& =   \om{r}{r^{-1}}^{*} \Inner{\gammArg{\E}{r}{\gammArg{\E}{r^{-1}}{\zeta}}}{\gammArg{\E}{r}{\eta}}{\E} \\
& =   \om{r}{r^{-1}}^{*} \alphArg{r}{\Inner{\gammArg{\E}{r^{-1}}{\zeta}}{\eta}{\E}} \\
& =   0, \quad \text{so} \\
      \gammArg{\E}{r}{\eta}
& \in \F^{\perp}.
\end{align*}
Let $ P: \E \to \F $ denote the projection map from $ \E $ onto $ \F $, which is an adjointable operator whose adjoint is the inclusion map $ \iota: \F \hookrightarrow \E $. By the foregoing argument, we have
\begin{align*}
\forall r \in G, ~ \forall \zeta \in \E: \quad
    \func{P}{\gammArg{\E}{r}{\zeta}}
& = \func{P}{\gammArg{\E}{r}{\func{P}{\zeta} \oplus \Func{\Id - P}{\zeta}}} \\
& = \func{P}{\gammArg{\E}{r}{\func{P}{\zeta}} \oplus \gammArg{\E}{r}{\Func{\Id - P}{\zeta}}} \\
& = \gammArg{\E}{r}{\func{P}{\zeta}}.
\end{align*}
Therefore, $ P \in \AdjEqPair{\E}{\F} $.
\end{Eg}


We now define a $ * $-representation of $ A $ by $ \Hilb{G}{A}{\alpha}{\omega} $-endomorphisms on $ \LL{2}{G,A} $ that is both faithful and non-degenerate.


\begin{Eg} \label{Covariant Representation Part 1}
For every $ a \in A $, $ \phi \in \Cc{G,A} $ and $ x \in G $, we have
$$
         \SqBr{\alphArg{x}{a} ~ \func{\phi}{x}}^{*} \SqBr{\alphArg{x}{a} ~ \func{\phi}{x}}
=        \func{\phi}{x}^{*} \alphArg{x}{a}^{*} \alphArg{x}{a} ~ \func{\phi}{x}
\leq_{A} \Norm{\alphArg{x}{a}}{A}^{2} \func{\phi}{x}^{*} \func{\phi}{x}
=        \Norm{a}{A}^{2} \func{\phi}{x}^{*} \func{\phi}{x}, \quad \text{so}
$$
\begin{align*}
       \Norm{\func{q}{\Map{G}{A}{x}{\alphArg{x}{a} ~ \func{\phi}{x}}}}{\LL{2}{G,A}}
& =    \NNorm{\Map{G}{A}{x}{\alphArg{x}{a} ~ \func{\phi}{x}}} \\
& =    \Norm{\Int{G}{\SqBr{\alphArg{x}{a} ~ \func{\phi}{x}}^{*} \SqBr{\alphArg{x}{a} ~ \func{\phi}{x}}}{x}}{A}^{\frac{1}{2}} \\
& \leq \Norm{\Norm{a}{A}^{2} \Int{G}{\func{\phi}{x}^{*} \func{\phi}{x}}{x}}{A}^{\frac{1}{2}} \\
& =    \Norm{a}{A} \Norm{\Int{G}{\func{\phi}{x}^{*} \func{\phi}{x}}{x}}{A}^{\frac{1}{2}} \\
& =    \Norm{a}{A} \NNorm{\phi} \\
& =    \Norm{a}{A} \Norm{\func{q}{\phi}}{\LL{2}{G,A}}.
\end{align*}
We can thus define a map $ \pi $ from $ A $ to the set of bounded operators on $ \LL{2}{G,A} $ by
$$
\forall a \in A, ~ \forall \Phi \in \LL{2}{G,A}: \quad
\FUNC{\func{\pi}{a}}{\Phi} \df \lim_{n \to \infty} \func{q}{\Map{G}{A}{x}{\alphArg{x}{a} ~ \func{\phi_{n}}{x}}},
$$
where $ \Seq{\phi_{n}}{n \in \N{}} $ is \emph{any} sequence in $ \Cc{G,A} $ with $ \D \lim_{n \to \infty} \func{q}{\phi_{n}} = \Phi $.

For every $ a,b \in A $, it is easy to check the following:
\begin{itemize}
\item
$ \func{\pi}{a b} = \func{\pi}{a} \circ \func{\pi}{b} $.

\item
$ \func{\pi}{a} $ is adjointable with $ \func{\pi}{a}^{*} = \func{\pi}{a^{*}} $.
\end{itemize}
Hence, $ \pi $ is a $ * $-representation of $ A $ on $ \LL{2}{G,A} $, and we also have that $ \pi: A \to \AdjEq{\LL{2}{G,A}} $ --- observe for every $ r \in G $, $ a \in A $ and $ \phi \in \Cc{G,A} $ that
\begin{align*}
    \GammArg{r}{\FUNC{\func{\pi}{a}}{\func{q}{\phi}}}
& = \GammArg{r}{\func{q}{\Map{G}{A}{x}{\alphArg{x}{a} ~ \func{\phi}{x}}}} \\
& = \func{q}{\Map{G}{A}{x}{\om{r}{r^{-1} x}^{*} \alphArg{r}{\alphArg{r^{-1} x}{a} ~ \func{\phi}{r^{-1} x}}}} \\
& = \func{q}{\Map{G}{A}{x}{\om{r}{r^{-1} x}^{*} \alphArg{r}{\alphArg{r^{-1} x}{a}} ~ \alphArg{r}{\func{\phi}{r^{-1} x}}}} \\
& = \func{q}{
            \Map{G}{A}{x}{\om{r}{r^{-1} x}^{*} \om{r}{r^{-1} x} ~ \alphArg{x}{a} ~ \om{r}{r^{-1} x}^{*} \alphArg{r}{\func{\phi}{r^{-1} x}}}
            } \\
& = \func{q}{\Map{G}{A}{x}{\alphArg{x}{a} ~ \om{r}{r^{-1} x}^{*} \alphArg{r}{\func{\phi}{r^{-1} x}}}} \\
& = \FUNC{\func{\pi}{a}}{\func{q}{\Map{G}{A}{x}{\om{r}{r^{-1} x}^{*} \alphArg{r}{\func{\phi}{r^{-1} x}}}}} \\
& = \FUNC{\func{\pi}{a}}{\GammArg{r}{\func{q}{\phi}}},
\end{align*}
so by continuity, $ \GammArg{r}{\FUNC{\func{\pi}{a}}{\Phi}} = \FUNC{\func{\pi}{a}}{\GammArg{r}{\Phi}} $ for every $ \Phi \in \LL{2}{G,A} $.

Now, we prove that $ \pi $ is faithful. Suppose that $ a \in A $ and $ \func{\pi}{a} = 0_{\AdjEq{\LL{2}{G,A}}} $. Pick a non-zero $ \phi \in \Cc{G,\RR{\geq 0}{}} $ so that $ \func{\phi}{e} > 0 $. Then $ \Map{G}{A}{x}{\func{\phi}{x} \cdot \alphArg{x}{a^{*}}} \in \Cc{G,A} $, and
$$
  0_{\LL{2}{G,A}}
= \FUNC{\func{\pi}{a}}{\func{q}{\Map{G}{A}{x}{\func{\phi}{x} \cdot \alphArg{x}{a^{*}}}}}
= \func{q}{\Map{G}{A}{x}{\func{\phi}{x} \cdot \alphArg{x}{a a^{*}}}}.
$$
It follows that $ \func{\phi}{e} \cdot a a^{*} = 0_{A} $, or equivalently, $ a = 0_{A} $ because $ \func{\phi}{e} > 0 $. Therefore, $ \pi $ is faithful.

Next, we prove that $ \pi $ is non-degenerate. Let $ \phi \in \Cc{G,A} $. By (5) of \autoref{Twisted C*-Dynamical Systems}, we have
$$
\forall r \in G, ~ \forall a \in A: \quad
\func{\alph{r}^{-1}}{a} = \om{r^{-1}}{r}^{*} \alphArg{r^{-1}}{a} ~ \om{r^{-1}}{r}.
$$
Hence, $ \Map{G}{A}{x}{\func{\alph{x}^{-1}}{\func{\phi}{x}}} $ is a continuous function, making $ K \df \Set{\func{\alph{x}^{-1}}{\func{\phi}{x}}}{x \in \Supp{\phi}} $ a compact subset of $ A $. The $ C^{*} $-subalgebra $ B $ of $ A $ generated by $ K $ is thus separable, so there is a sequential approximate identity $ \Seq{e_{n}}{n \in \N{}} $ for $ B $. Then
\begin{align*}
\forall n \in \N{}: \quad
       \Norm{\Func{\func{\pi}{e_{n}}}{\func{q}{\phi}} - \func{q}{\phi}}{\LL{2}{G,A}}
& =    \NNorm{\Map{G}{A}{x}{\alphArg{x}{e_{n}} ~ \func{\phi}{x} - \func{\phi}{x}}} \\
& =    \Norm{
            \Int{G}{
                   \SqBr{\alphArg{x}{e_{n}} ~ \func{\phi}{x} - \func{\phi}{x}}^{*} \SqBr{\alphArg{x}{e_{n}} ~ \func{\phi}{x} - \func{\phi}{x}}
                   }{x}
            }{A}^{\frac{1}{2}} \\
& \leq \SqBr{\Int{G}{\Norm{\alphArg{x}{e_{n}} ~ \func{\phi}{x} - \func{\phi}{x}}{A}^{2}}{x}}^{\frac{1}{2}} \\
& =    \SqBr{\Int{G}{\Norm{\alphArg{x}{e_{n} ~ \func{\alph{x}^{-1}}{\func{\phi}{x}}} - \func{\phi}{x}}{A}^{2}}{x}}^{\frac{1}{2}} \\
& =    \SqBr{\Int{\Supp{\phi}}{\Norm{\alphArg{x}{e_{n} ~ \func{\alph{x}^{-1}}{\func{\phi}{x}}} - \func{\phi}{x}}{A}^{2}}{x}}^{\frac{1}{2}}.
\end{align*}
The integrand in the last line is dominated by the integrable function $ \Map{G}{\RR{\geq 0}{}}{x}{4 \Norm{\func{\phi}{x}}{A}^{2}} $. As
$$
\forall x \in \Supp{\phi}: \quad
\lim_{n \to \infty} \Norm{\alphArg{x}{e_{n} ~ \func{\alph{x}^{-1}}{\func{\phi}{x}}} - \func{\phi}{x}}{A}^{2} = 0,
$$
the Lebesgue Dominated Convergence Theorem then implies that
$$
  \lim_{n \to \infty}
  \SqBr{\Int{\Supp{\phi}}{\Norm{\alphArg{x}{e_{n} ~ \func{\alph{x}^{-1}}{\func{\phi}{x}}} - \func{\phi}{x}}{A}^{2}}{x}}^{\frac{1}{2}}
= 0.
$$
By the Squeeze Theorem, $ \D \lim_{n \to \infty} \Norm{\FUNC{\func{\pi}{e_{n}}}{\func{q}{\phi}} - \func{q}{\phi}}{\LL{2}{G,A}} = 0 $. As $ \phi $ is arbitrary and $ \Im{q}{\Cc{G,A}} $ is dense in $ \LL{2}{G,A} $, we conclude that $ \pi $ is non-degenerate.
\end{Eg}


We also have a multiplier representation of $ G $ by $ \Hilb{G}{A}{\alpha}{\omega} $-endomorphisms on $ \LL{2}{G,A} $.


\begin{Eg} \label{Covariant Representation Part 2}
Observe for every $ r \in G $ and $ \phi \in \Cc{G,A} $ that
\begin{align*}
  & ~ \Norm{\func{q}{\Map{G}{A}{x}{\Del{r}{\frac{1}{2}} \om{x}{r} ~ \func{\phi}{x r}}}}{\LL{2}{G,A}} \\
= & ~ \NNorm{\Map{G}{A}{x}{\Del{r}{\frac{1}{2}} \om{x}{r} ~ \func{\phi}{x r}}} \\
= & ~ \Norm{
           \Int{G}{\SqBr{\Del{r}{\frac{1}{2}} \om{x}{r} ~ \func{\phi}{x r}}^{*} \SqBr{\Del{r}{\frac{1}{2}} \om{x}{r} ~ \func{\phi}{x r}}}{x}
           }{A}^{\frac{1}{2}} \\
= & ~ \Norm{\Int{G}{\Del{r}{} ~ \func{\phi}{x r}^{*} \om{x}{r}^{*} \om{x}{r} ~ \func{\phi}{x r}}{x}}{A}^{\frac{1}{2}} \\
= & ~ \Norm{\Int{G}{\Del{r}{} ~ \func{\phi}{x r}^{*} \func{\phi}{x r}}{x}}{A}^{\frac{1}{2}} \\
= & ~ \Norm{\Int{G}{\func{\phi}{x}^{*} \func{\phi}{x}}{x}}{A}^{\frac{1}{2}} \\
= & ~ \NNorm{\phi} \\
= & ~ \Norm{\func{q}{\phi}}{\LL{2}{G,A}}.
\end{align*}
We can thus define a map $ \lambda: G \to \Isom{\LL{2}{G,A}} $ by
$$
\forall r \in G, ~ \forall \Phi \in \LL{2}{G,A}: \quad
\FUNC{\func{\lambda}{r}}{\Phi} \df \lim_{n \to \infty} \func{q}{\Map{G}{A}{x}{\Del{r}{\frac{1}{2}} \om{x}{r} ~ \func{\phi_{n}}{x r}}},
$$
where $ \Seq{\phi_{n}}{n \in \N{}} $ is \emph{any} sequence in $ \Cc{G,A} $ with $ \D \lim_{n \to \infty} \func{q}{\phi_{n}} = \Phi $. It is then not difficult to check the following:
\begin{itemize}
\item
$ \func{\lambda}{r} $ is unitary for every $ r \in G $.

\item
For every $ r \in G $ and $ \Phi \in \LL{2}{G,A} $,
$$
  \FUNC{\func{\lambda}{r}^{*}}{\Phi}
= \lim_{n \to \infty} \func{q}{\Map{G}{A}{x}{\Del{r}{- \frac{1}{2}} \om{x r^{-1}}{r}^{*} \func{\phi_{n}}{x r^{-1}}}},
$$
where $ \Seq{\phi_{n}}{n \in \N{}} $ is \emph{any} sequence in $ \Cc{G,A} $ with $ \D \lim_{n \to \infty} \func{q}{\phi_{n}} = \Phi $.

\item
$ \func{\lambda}{r} \circ \func{\lambda}{s} = \func{\Conj{\pi}}{\om{r}{s}} \circ \func{\lambda}{r s} $ for every $ r,s \in G $, where $ \pi: A \to \AdjEq{\LL{2}{G,A}} $ is as defined in \autoref{Covariant Representation Part 1}. Hence, $ \lambda $ is not a unitary representation of $ G $ on $ \LL{2}{G,A} $ unless $ \omega $ is trivial.
\end{itemize}
To obtain $ \lambda: G \to \U{\AdjEq{\LL{2}{G,A}}} $, observe for every $ r,s \in G $ and $ \phi \in \Cc{G,A} $ that
\begin{align*}
    \GammArg{r}{\FUNC{\func{\lambda}{s}}{\func{q}{\phi}}}
& = \GammArg{r}{\func{q}{\Map{G}{A}{x}{\Del{s}{\frac{1}{2}} \om{x}{s} ~ \func{\phi}{x s}}}} \\
& = \func{q}{\Map{G}{A}{x}{\om{r}{r^{-1} x}^{*} \alphArg{r}{\Del{s}{\frac{1}{2}} \om{r^{-1} x}{s} ~ \func{\phi}{r^{-1} x s}}}} \\
& = \func{q}{\Map{G}{A}{x}{\Del{s}{\frac{1}{2}} \om{r}{r^{-1} x}^{*} \alphArg{r}{\om{r^{-1} x}{s} ~ \func{\phi}{r^{-1} x s}}}} \\
& = \func{q}{
            \Map{G}{A}{x}{\Del{s}{\frac{1}{2}} \om{r}{r^{-1} x}^{*} \alphArgExt{r}{\om{r^{-1} x}{s}} ~ \alphArg{r}{\func{\phi}{r^{-1} x s}}}
            } \\
& = \func{q}{\Map{G}{A}{x}{\Del{s}{\frac{1}{2}} \om{x}{s} ~ \om{r}{r^{-1} x s}^{*} \alphArg{r}{\func{\phi}{r^{-1} x s}}}} \\
& = \FUNC{\func{\lambda}{s}}{\func{q}{\Map{G}{A}{x}{\om{r}{r^{-1} x}^{*} \alphArg{r}{\func{\phi}{r^{-1} x}}}}} \\
& = \FUNC{\func{\lambda}{s}}{\GammArg{r}{\func{q}{\phi}}},
\end{align*}
so by continuity, $ \GammArg{r}{\FUNC{\func{\lambda}{s}}{\Phi}} = \FUNC{\func{\lambda}{s}}{\GammArg{r}{\Phi}} $ for every $ \Phi \in \LL{2}{G,A} $.

A point of observation is that $ \lambda $ is strongly continuous. Let $ \epsilon > 0 $, $ r \in G $, $ \phi \in \Cc{G,A} \setminus \SSet{\bm{0}} $ and $ S \df \Supp{\phi} $. Fix a compact subset $ K $ of $ G $ containing $ r $ in its interior. Then by continuity, we can find $ S K^{-1} $-indexed sequences $ \Seq{V_{x}}{x \in S K^{-1}} $ and $ \Seq{W_{x}}{x \in S K^{-1}} $ of subsets of $ G $ with the following properties for every $ x \in S K^{-1} $:
\begin{itemize}
\item
$ V_{x} $ is the intersection of $ S K^{-1} $ with an open neighborhood of $ x $.

\item
$ W_{x} $ is the intersection of $ K^{\circ} $ with an open neighborhood of $ r $.

\item
$ \D \forall \Pair{y}{s} \in V_{x} \times W_{x}: \quad \Norm{\Del{s}{- \frac{1}{2}} \om{y}{s} ~ \func{\phi}{y s} - \Del{r}{- \frac{1}{2}} \om{x}{r} ~ \func{\phi}{x r}}{A} < \dfrac{\epsilon}{2 \sqrt{\func{\mu}{S K^{-1}}}} $, whence
$$
\forall \Pair{y}{s} \in V_{x} \times W_{x}: \quad
  \Norm{\Del{s}{- \frac{1}{2}} \om{y}{s} ~ \func{\phi}{y s} - \Del{r}{- \frac{1}{2}} \om{y}{r} ~ \func{\phi}{y r}}{A}
< \dfrac{\epsilon}{\sqrt{\func{\mu}{S K^{-1}}}}.
$$
\end{itemize}
By the compactness of $ S K^{-1} $, there exist points $ x_{1},\ldots,x_{n} \in S K^{-1} $ that satisfy $ \D S K^{-1} = \bigcup_{k = 1}^{n} V_{x_{k}} $. Pick any open neighborhood $ N $ of $ r $ contained within $ \D \bigcap_{k = 1}^{n} W_{x_{k}} $, and let $ \Pair{x}{s} \in S K^{-1} \times N $. Find a $ k \in \SSet{1,\ldots,n} $ such that $ x \in V_{x_{k}} $. As $ s \in W_{x_{k}} $, we have
$$
  \Norm{\Del{s}{- \frac{1}{2}} \om{x}{s} ~ \func{\phi}{x s} - \Del{r}{- \frac{1}{2}} \om{x}{r} ~ \func{\phi}{x r}}{A}
< \dfrac{\epsilon}{\sqrt{\func{\mu}{S K^{-1}}}}.
$$
We chose $ x $ arbitrarily, so
\begin{align*}
     & ~ \Norm{\FUNC{\func{\lambda}{s}}{\func{q}{\phi}} - \FUNC{\func{\lambda}{r}}{\func{q}{\phi}}}{\LL{2}{G,A}} \\
=    & ~ \NNorm{\Map{G}{A}{x}{\Del{s}{- \frac{1}{2}} \om{x}{s} ~ \func{\phi}{x s} - \Del{r}{- \frac{1}{2}} \om{x}{r} ~ \func{\phi}{x r}}} \\
\leq & ~ \Norm{\Map{G}{A}{x}{\Del{s}{- \frac{1}{2}} \om{x}{s} ~ \func{\phi}{x s} - \Del{r}{- \frac{1}{2}} \om{x}{r} ~ \func{\phi}{x r}}}{2} \\
=    & ~ \SqBr{
              \Int{G}{
                     \Norm{\Del{s}{- \frac{1}{2}} \om{x}{s} ~ \func{\phi}{x s} - \Del{r}{- \frac{1}{2}} \om{x}{r} ~ \func{\phi}{x r}}{A}^{2}
                     }{x}
              }^{\frac{1}{2}} \\
=    & ~ \SqBr{
              \Int{S K^{-1}}{
                            \Norm{
                                 \Del{s}{- \frac{1}{2}} \om{x}{s} ~ \func{\phi}{x s} - \Del{r}{- \frac{1}{2}} \om{x}{r} ~ \func{\phi}{x r}
                                 }{A}^{2}
                            }{x}
              }^{\frac{1}{2}} \\
     & ~ \Br{\text{As the integrand vanishes outside of $ S K^{-1} $.}} \\
<    & ~ \SqBr{\Int{S K^{-1}}{\frac{\epsilon^{2}}{\func{\mu}{S K^{-1}}}}{x}}^{\frac{1}{2}} \\
=    & ~ \SqBr{\frac{\epsilon^{2}}{\func{\mu}{S K^{-1}}} \cdot \func{\mu}{S K^{-1}}}^{\frac{1}{2}} \\
=    & ~ \epsilon.
\end{align*}
As $ \epsilon $ and $ \phi $ are arbitrary, we get $ \D \lim_{s \to r} \Norm{\FUNC{\func{\lambda}{s}}{\func{q}{\phi}} - \FUNC{\func{\lambda}{r}}{\func{q}{\phi}}}{\LL{2}{G,A}} = 0 $ for every $ \phi \in \Cc{G,A} $.

Let $ \Phi \in \LL{2}{G,A} $. Let $ \epsilon > 0 $ once more, and pick a $ \phi \in \Cc{G,A} $ so that $ \Norm{\Phi - \func{q}{\phi}}{\LL{2}{G,A}} < \dfrac{\epsilon}{3} $. By the argument above, there is an open neighborhood $ N $ of $ r $ such that
$$
\forall s \in N: \quad
\Norm{\FUNC{\func{\lambda}{s}}{\func{q}{\phi}} - \FUNC{\func{\lambda}{r}}{\func{q}{\phi}}}{\LL{2}{G,A}} < \frac{\epsilon}{3},
$$
from which it follows that
\begin{align*}
\forall s \in N: \qquad
     & ~ \Norm{\FUNC{\func{\lambda}{s}}{\Phi} - \FUNC{\func{\lambda}{r}}{\Phi}}{\LL{2}{G,A}} \\
\leq & ~ \Norm{\FUNC{\func{\lambda}{s}}{\Phi} - \FUNC{\func{\lambda}{s}}{\func{q}{\phi}}}{\LL{2}{G,A}} + \\
     & ~ \Norm{\FUNC{\func{\lambda}{s}}{\func{q}{\phi}} - \FUNC{\func{\lambda}{r}}{\func{q}{\phi}}}{\LL{2}{G,A}} + \\
     & ~ \Norm{\FUNC{\func{\lambda}{r}}{\func{q}{\phi}} - \FUNC{\func{\lambda}{r}}{\Phi}}{\LL{2}{G,A}} \\
=    & ~ \Norm{\Phi - \func{q}{\phi}}{\LL{2}{G,A}} +
         \Norm{\FUNC{\func{\lambda}{s}}{\func{q}{\phi}} - \FUNC{\func{\lambda}{r}}{\func{q}{\phi}}}{\LL{2}{G,A}} +
         \Norm{\func{q}{\phi} - \Phi}{\LL{2}{G,A}} \\
<    & ~ \frac{\epsilon}{3} + \frac{\epsilon}{3} + \frac{\epsilon}{3} \\
=    & ~ \epsilon.
\end{align*}
As $ \epsilon $ is arbitrary, we obtain $ \D \lim_{s \to r} \Norm{\FUNC{\func{\lambda}{s}}{\Phi} - \FUNC{\func{\lambda}{r}}{\Phi}}{\LL{2}{G,A}} = 0 $. Then as $ r $ and $ \Phi $ are arbitrary, we conclude that $ \lambda $ is strongly continuous.
\end{Eg}



\section{Meyer's Bra-Ket Operators} 

\textbf{In this section, $ \E $ is a Hilbert $ \Quad{G}{A}{\alpha}{\omega} $-module.}

\subsection{Square-Integrability} 


\begin{Def} \label{Meyer's Bra-Ket Operators}
For each $ \zeta \in \E $, we can define operators $ \Bra{\zeta}: \E \to \Cb{G,A} $ and $ \Ket{\zeta}: \Cc{G,A} \to \E $, called the \emph{bra} and \emph{ket} of $ \zeta $ respectively, by
\begin{alignat*}{2}
\forall \eta \in \E:       & \quad &
\BraArg{\zeta}{\eta}       & \df \Map{G}{A}{x}{\Inner{\gammArg{\E}{x}{\zeta}}{\eta}{\E}}, \\
\forall \phi \in \Cc{G,A}: & \quad &
\KetArg{\zeta}{\phi}       & \df \Int{G}{\gammArg{\E}{x}{\zeta} \bullet \func{\phi}{x}}{x}.
\end{alignat*}
\end{Def}



\begin{Lem} \label{The Basic Norm Inequalities}
For every $ \zeta,\eta \in \E $ and $ \phi \in \Cc{G,A} $, the following norm inequalities hold:
\begin{align}
\BigNorm{\BraArg{\zeta}{\eta}}{\infty} & \leq \Norm{\zeta}{\E} \Norm{\eta}{\E}, \label{Basic Norm Inequality 1} \\
\BigNorm{\KetArg{\zeta}{\phi}}{\E}     & \leq \Norm{\zeta}{\E} \Norm{\phi}{1}.  \label{Basic Norm Inequality 2}
\end{align}
\end{Lem}

\begin{proof}
For every $ \zeta,\eta \in \E $ and $ \phi \in \Cc{G,A} $,
\begin{align*}
       \Norm{\BraArg{\zeta}{\eta}}{\infty}
& =    \sup_{x \in G} \Norm{\Inner{\gammArg{\E}{x}{\zeta}}{\eta}{\E}}{A} \\
& \leq \sup_{x \in G} \Norm{\gammArg{\E}{x}{\zeta}}{\E} \Norm{\eta}{\E} \qquad \Br{\text{By the Cauchy-Schwarz Inequality.}} \\
& =    \sup_{x \in G} \Norm{\zeta}{\E} \Norm{\eta}{\E} \\
& =    \Norm{\zeta}{\E} \Norm{\eta}{\E}, \\
       \Norm{\KetArg{\zeta}{\phi}}{\E}
& =    \Norm{\Int{G}{\gammArg{\E}{x}{\zeta} \bullet \func{\phi}{x}}{x}}{\E} \\
& \leq \Int{G}{\Norm{\gammArg{\E}{x}{\zeta} \bullet \func{\phi}{x}}{\E}}{x} \\
& \leq \Int{G}{\Norm{\gammArg{\E}{x}{\zeta}}{\E} \Norm{\func{\phi}{x}}{A}}{x} \\
& =    \Int{G}{\Norm{\zeta}{\E} \Norm{\func{\phi}{x}}{A}}{x} \\
& =    \Norm{\zeta}{\E} \Int{G}{\Norm{\func{\phi}{x}}{A}}{x} \\
& =    \Norm{\zeta}{\E} \Norm{\phi}{1}.
\end{align*}
This concludes the proof.
\end{proof}


We now state a result to the effect that every element of $ \E $ yields a unique ket operator.


\begin{Prop} \label{Every Vector of a Twisted Hilbert C*-Module Gives a Unique Ket Operator}
Let $ \zeta,\eta \in \E $. If $ \Ket{\zeta} = \Ket{\eta} $, then $ \zeta = \eta $.
\end{Prop}

\begin{proof}
Suppose that $ \Ket{\zeta} = \Ket{\eta} $. Let $ \mathcal{N} $ denote the open-neighborhood base of $ e $ directed by inclusion. Let $ \Seq{f_{N}}{N \in \mathcal{N}} $ be a net in $ \Cc{G,\RR{\geq 0}{}} $ so that $ \Supp{f_{N}} \subseteq N $ and $ \D \Int{G}{\func{f_{N}}{x}}{x} = 1 $ for every $ N \in \mathcal{N} $. Then $ \Seq{f_{N}}{N \in \mathcal{N}} $ is an approximating delta at $ e $, and if $ A' $ denotes the dual space of $ A $, we have
\begin{align*}
\forall N \in \mathcal{N}, ~ \forall a \in A, ~ \forall \varphi \in A': \quad
    \Int{G}{\func{f_{N}}{x} ~ \func{\varphi}{\gammArg{\E}{x}{\zeta} \bullet a}}{x}
& = \func{\varphi}{\Int{G}{\gammArg{\E}{x}{\zeta} \bullet \SqBr{\func{f_{N}}{x} ~ a}}{x}} \\
& = \func{\varphi}{\KetArg{\zeta}{f_{N} a}} \\
& = \func{\varphi}{\KetArg{\eta}{f_{N} a}} \\
& = \func{\varphi}{\Int{G}{\gammArg{\E}{x}{\eta} \bullet \SqBr{\func{f_{N}}{x} ~ a}}{x}} \\
& = \Int{G}{\func{f_{N}}{x} ~ \func{\varphi}{\gammArg{\E}{x}{\eta} \bullet a}}{x}, \quad \text{so} \\
    \func{\varphi}{\zeta \bullet a}
& = \func{\varphi}{\gammArg{\E}{e}{\zeta} \bullet a} \\
& = \lim_{N \in \mathcal{N}} \Int{G}{\func{f_{N}}{x} ~ \func{\varphi}{\gammArg{\E}{x}{\zeta} \bullet a}}{x} \\
& = \lim_{N \in \mathcal{N}} \Int{G}{\func{f_{N}}{x} ~ \func{\varphi}{\gammArg{\E}{x}{\eta} \bullet a}}{x} \\
& = \func{\varphi}{\gammArg{\E}{e}{\eta} \bullet a} \\
& = \func{\varphi}{\eta \bullet a}.
\end{align*}
By the Hahn-Banach Theorem, $ \zeta \bullet a = \eta \bullet a $ for every $ a \in A $. Therefore, $ \zeta = \eta $.
\end{proof}



\begin{Def} \label{Square-Integrability}
We say that $ \zeta \in \E $ is \emph{square-integrable} if and only if for every $ \eta \in \E $ and every net $ \Seq{\varphi_{i}}{i \in I} $ in $ \Cc{G,\CC{0}{1}} $ converging uniformly to $ 1 $ on compact subsets of $ G $, the net $ \Seq{\func{q}{\varphi_{i} \BraArg{\zeta}{\eta}}}{i \in I} $ is Cauchy in $ \LL{2}{G,A} $, in which case we can define an operator $ \BBra{\zeta}: \E \to \LL{2}{G,A} $ by
$$
\forall \eta \in \E: \quad
\BBraArg{\zeta} \eta \df \lim_{i \in I} \func{q}{\varphi_{i} \BraArg{\zeta}{\eta}}.
$$
The definition of $ \BBra{\zeta} $ is independent of our choice of $ \Seq{\varphi_{i}}{i \in I} $, which we will establish in a moment. The set of square-integrable elements of $ \E $ is denoted by $ \Esi $ --- it is clearly a linear subspace of $ \E $. We call $ \E $ a \emph{square-integrable representation} of $ \Quad{G}{A}{\alpha}{\omega} $ if and only if $ \Esi $ is dense in $ \E $.
\end{Def}



\begin{Lem} \label{Multiplication Operators}
For every $ \varphi \in \Cb{G} $, there exists a unique $ M_{\varphi} \in \Adj{\LL{2}{G,A}} $, with norm $ \leq \Norm{\varphi}{\infty} $, such that $ \func{M_{\varphi}}{\func{q}{\phi}} = \func{q}{\Map{G}{A}{x}{\func{\varphi}{x} ~ \func{\phi}{x}}} $ for every $ \phi \in \Cc{G,A} $. Furthermore, $ M_{\varphi}^{*} = M_{\Conj{\varphi}} $ for every $ \varphi \in \Cc{G} $.
\end{Lem}

\begin{proof}
Let $ \varphi \in \Cb{G} $. Then for every $ \phi \in \Cc{G,A} $ and $ x \in G $, we have
$$
         0_{A}
\leq_{A} \SqBr{\func{\varphi}{x} ~ \func{\phi}{x}}^{*} \SqBr{\func{\varphi}{x} ~ \func{\phi}{x}}
=        \Abs{\func{\varphi}{x}}^{2} \func{\phi}{x}^{*} \func{\phi}{x}
\leq_{A} \Norm{\varphi}{\infty}^{2} \func{\phi}{x}^{*} \func{\phi}{x}, \quad \text{so}
$$
$$
         \Int{G}{\SqBr{\func{\varphi}{x} ~ \func{\phi}{x}}^{*} \SqBr{\func{\varphi}{x} ~ \func{\phi}{x}}}{x}
\leq_{A} \Norm{\varphi}{\infty}^{2} \Int{G}{\func{\phi}{x}^{*} \func{\phi}{x}}{x}, \quad \text{which yields}
$$
\begin{align*}
       \Norm{\func{q}{\Map{G}{A}{x}{\func{\varphi}{x} ~ \func{\phi}{x}}}}{\LL{2}{G,A}}
& =    \NNorm{\Map{G}{A}{x}{\func{\varphi}{x} ~ \func{\phi}{x}}} \\
& =    \Norm{\Int{G}{\SqBr{\func{\varphi}{x} ~ \func{\phi}{x}}^{*} \SqBr{\func{\varphi}{x} ~ \func{\phi}{x}}}{x}}{A}^{\frac{1}{2}} \\
& \leq \Norm{\Norm{\varphi}{\infty}^{2} \Int{G}{\func{\phi}{x}^{*} \func{\phi}{x}}{x}}{A}^{\frac{1}{2}} \\
& =    \Norm{\varphi}{\infty} \Norm{\Int{G}{\func{\phi}{x}^{*} \func{\phi}{x}}{x}}{A}^{\frac{1}{2}} \\
& =    \Norm{\varphi}{\infty} \NNorm{\phi} \\
& =    \Norm{\varphi}{\infty} \Norm{\func{q}{\phi}}{\LL{2}{G,A}}.
\end{align*}
Hence, there is a unique bounded operator $ M_{\varphi} $ on $ \LL{2}{G,A} $, with norm $ \leq \Norm{\varphi}{\infty} $, such that
$$
\forall \phi \in \Cc{G,A}: \quad
\func{M_{\varphi}}{\func{q}{\phi}} = \func{q}{\Map{G}{A}{x}{\func{\varphi}{x} ~ \func{\phi}{x}}}.
$$
Next, observe that
\begin{align*}
\forall \phi,\psi \in \Cc{G,A}: \quad
    \Inner{\func{M_{\varphi}}{\func{q}{\phi}}}{\func{q}{\psi}}{\LL{2}{G,A}}
& = \Int{G}{\SqBr{\func{\varphi}{x} ~ \func{\phi}{x}}^{*} \func{\psi}{x}}{x} \\
& = \Int{G}{\func{\phi}{x}^{*} \SqBr{\func{\Conj{\varphi}}{x} ~ \func{\psi}{x}}}{x} \\
& = \Inner{\func{q}{\phi}}{\func{M_{\Conj{\varphi}}}{\func{q}{\psi}}}{\LL{2}{G,A}},
\end{align*}
so by continuity, $ \Inner{\func{M_{\varphi}}{\Phi}}{\Psi}{\LL{2}{G,A}} = \Inner{\Phi}{\func{M_{\Conj{\varphi}}}{\Psi}}{\LL{2}{G,A}} $ for every $ \Phi,\Psi \in \LL{2}{G,A} $. Therefore, $ M_{\varphi} $ is adjointable with $ M_{\Conj{\varphi}} $ as its adjoint, and as $ \varphi $ is arbitrary, we are done.
\end{proof}



\begin{Lem} \label{A Convergence Lemma for Multiplication Operators}
Let $ \Seq{\varphi_{i}}{i \in I} $ be a net in $ \Cc{G,\CC{0}{1}} $ converging uniformly to $ 1 $ on compact subsets of $ G $. Then $ \D \lim_{i \in I} \func{M_{\varphi_{i}}}{\Phi} = \Phi $ for every $ \Phi \in \LL{2}{G,A} $.
\end{Lem}

\begin{proof}
Let $ \Phi \in \LL{2}{G,A} $ and $ \epsilon > 0 $. Pick a $ \phi \in \Cc{G,A} $ so that $ \Norm{\Phi - \func{q}{\phi}}{\LL{2}{G,A}} < \dfrac{\epsilon}{3} $. Then
\begin{align*}
\forall i \in I: \quad
       \Norm{\func{M_{\varphi_{i}}}{\Phi} - \func{M_{\varphi_{i}}}{\func{q}{\phi}}}{\LL{2}{G,A}}
& =    \Norm{\func{M_{\varphi_{i}}}{\Phi - \func{q}{\phi}}}{\LL{2}{G,A}} \\
& \leq \Norm{\varphi_{i}}{\infty} \Norm{\Phi - \func{q}{\phi}}{\LL{2}{G,A}} \quad \Br{\text{By \autoref{Multiplication Operators}.}} \\
& <    1 \cdot \frac{\epsilon}{3} \\
& =    \frac{\epsilon}{3}.
\end{align*}
Furthermore, for every $ i \in I $ and $ x \in \Supp{\phi} $, we have
\begin{align*}
           \Abs{\func{\varphi_{i}}{x} - 1}^{2} \func{\phi}{x}^{*} \func{\phi}{x}
& \leq_{A} \SqBr{\max_{x \in \Supp{\phi}} \Abs{\func{\varphi_{i}}{x} - 1}}^{2} \func{\phi}{x}^{*} \func{\phi}{x}, \quad \text{whence} \\
           \Int{\Supp{\phi}}{\Abs{\func{\varphi_{i}}{x} - 1}^{2} \func{\phi}{x}^{*} \func{\phi}{x}}{x}
& \leq_{A} \SqBr{\max_{\Supp{\phi}} \Abs{\func{\varphi_{i}}{x} - 1}}^{2} \Int{\Supp{\phi}}{\func{\phi}{x}^{*} \func{\phi}{x}}{x}.
\end{align*}
Pick an $ i_{0} \in I $ so that for every $ i \in I_{\geq i_{0}} $,
$$
  \max_{x \in \Supp{\phi}} \Abs{\func{\varphi_{i}}{x} - 1}
< \frac{\epsilon}{3 \Br{1 + \Norm{\func{q}{\phi}}{\LL{2}{G,A}}}}, \quad \text{in which case}
$$
\begin{align*}
     & ~ \Norm{\func{M_{\varphi_{i}}}{\Phi} - \Phi}{\LL{2}{G,A}} \\
\leq & ~ \Norm{\func{M_{\varphi_{i}}}{\Phi} - \func{M_{\varphi_{i}}}{\func{q}{\phi}}}{\LL{2}{G,A}} +
         \Norm{\func{M_{\varphi_{i}}}{\func{q}{\phi}} - \func{q}{\phi}}{\LL{2}{G,A}} +
         \Norm{\func{q}{\phi} - \Phi}{\LL{2}{G,A}} \\
<    & ~ \frac{2 \epsilon}{3} + \Norm{\func{M_{\varphi_{i}}}{\func{q}{\phi}} - \func{q}{\phi}}{\LL{2}{G,A}} \\
=    & ~ \frac{2 \epsilon}{3} + \Norm{\func{q}{\varphi_{i} \phi} - \func{q}{\phi}}{\LL{2}{G,A}} \\
=    & ~ \frac{2 \epsilon}{3} + \Norm{\func{q}{\Br{\varphi_{i} - \bm{1}} \phi}}{\LL{2}{G,A}} \\
=    & ~ \frac{2 \epsilon}{3} + \NNorm{\Br{\varphi_{i} - \bm{1}} \phi} \\
=    & ~ \frac{2 \epsilon}{3} +
         \Norm{\Int{G}{\FUNC{\Br{\varphi_{i} - \bm{1}} \phi}{x}^{*} \FUNC{\Br{\varphi_{i} - \bm{1}} \phi}{x}}{x}}{A}^{\frac{1}{2}} \\
=    & ~ \frac{2 \epsilon}{3} + \Norm{\Int{G}{\Abs{\func{\varphi_{i}}{x} - 1}^{2} \func{\phi}{x}^{*} \func{\phi}{x}}{x}}{A}^{\frac{1}{2}} \\
=    & ~ \frac{2 \epsilon}{3} +
         \Norm{\Int{\Supp{\phi}}{\Abs{\func{\varphi_{i} - 1}{x}}^{2} \func{\phi}{x}^{*} \func{\phi}{x}}{x}}{A}^{\frac{1}{2}} \\
\leq & ~ \frac{2 \epsilon}{3} +
         \Norm{
              \SqBr{\max_{x \in \Supp{\phi}} \Abs{\func{\varphi_{i}}{x} - 1}}^{2} \Int{\Supp{\phi}}{\func{\phi}{x}^{*} \func{\phi}{x}}{x}
              }{A}^{\frac{1}{2}} \\
=    & ~ \frac{2 \epsilon}{3} +
         \SqBr{\max_{x \in \Supp{\phi}} \Abs{\func{\varphi_{i}}{x} - 1}}
         \Norm{\Int{\Supp{\phi}}{\func{\phi}{x}^{*} \func{\phi}{x}}{x}}{A}^{\frac{1}{2}} \\
=    & ~ \frac{2 \epsilon}{3} +
         \SqBr{\max_{x \in \Supp{\phi}} \Abs{\func{\varphi_{i} - 1}{x}}} \Norm{\Int{G}{\func{\phi}{x}^{*} \func{\phi}{x}}{x}}{A}^{\frac{1}{2}}
         \\
\leq & ~ \frac{2 \epsilon}{3} + \frac{\epsilon}{3 \Br{1 + \Norm{\func{q}{\phi}}{\LL{2}{G,A}}}} \NNorm{\phi} \\
=    & ~ \frac{2 \epsilon}{3} + \frac{\epsilon}{3 \Br{1 + \Norm{\func{q}{\phi}}{\LL{2}{G,A}}}} \Norm{\func{q}{\phi}}{\LL{2}{G,A}} \\
<    & ~ \epsilon.
\end{align*}
As $ \epsilon $ is arbitrary, we obtain $ \D \lim_{i \in I} \func{M_{\varphi_{i}}}{\Phi} = \Phi $, and as $ \Phi $ is arbitrary, we are done.
\end{proof}


We now return to the unproven assertion that for $ \zeta \in \Esi $, the definition of $ \BBra{\zeta} $ in \autoref{Square-Integrability} does not depend on our choice of a net $ \Seq{\varphi_{i}}{i \in I} $ having the properties listed there. Let $ \Seq{\psi_{j}}{j \in J} $ be another net with the same properties. Then the continuity of $ M_{\varphi_{i}} $ implies that
$$
\forall i \in I, ~ \forall \eta \in \E: \quad
  \lim_{j \in J} \func{q}{\varphi_{i} \psi_{j} \BraArg{\zeta}{\eta}}
= \lim_{j \in J} \func{M_{\varphi_{i}}}{\func{q}{\psi_{j} \BraArg{\zeta}{\eta}}}
= \func{M_{\varphi_{i}}}{\lim_{j \in J} \func{q}{\psi_{j} \BraArg{\zeta}{\eta}}},
$$
whereas \autoref{A Convergence Lemma for Multiplication Operators} implies that
$$
\forall i \in I, ~ \forall \eta \in \E: \quad
  \lim_{j \in J} \func{q}{\varphi_{i} \psi_{j} \BraArg{\zeta}{\eta}}
= \lim_{j \in J} \func{M_{\psi_{j}}}{\func{q}{\varphi_{i} \BraArg{\zeta}{\eta}}}
= \func{q}{\varphi_{i} \BraArg{\zeta}{\eta}}.
$$
Hence, by another application of \autoref{A Convergence Lemma for Multiplication Operators},
$$
\forall \eta \in \E: \quad
  \lim_{j \in J} \func{q}{\psi_{j} \BraArg{\zeta}{\eta}}
= \lim_{i \in I} \func{M_{\varphi_{i}}}{\lim_{j \in J} \func{q}{\psi_{j} \BraArg{\zeta}{\eta}}}
= \lim_{i \in I} \func{q}{\varphi_{i} \BraArg{\zeta}{\eta}}.
$$
The definition of $ \BBra{\zeta} $ is therefore consistent as it stands.


\begin{Prop} \label{A Direct Sum of Square-Integrable Representations Is Square-Integrable}
Let $ \Seq{\E_{n}}{n \in \N{}{}} $ be a sequence of square-integrable representations of $ \Quad{G}{A}{\alpha}{\omega} $. Then the direct sum $ \D \bigoplus_{n = 1}^{\infty} \E_{n} $ is also a square-integrable representation of $ \Quad{G}{A}{\alpha}{\omega} $.
\end{Prop}

\begin{proof}
Let $ \D\bm{\zeta} \in \bigoplus_{n = 1}^{\infty} \E_{n} $. Suppose that $ \bm{\zeta} $ has finitely many non-zero components so that $ \D \bm{\zeta} = \sum_{n \in N} \zeta_{n} \cdot \mathbf{e}_{n} $ for some finite subset $ N $ of $ \N{}{} $. Then
$$
\forall \bm{\eta} \in \bigoplus_{n = 1}^{\infty} \E_{n}: \quad
  \BraArg{\bm{\zeta}}{\bm{\eta}}
= \sum_{n \in N} \BraArg{\zeta_{n} \cdot \mathbf{e}_{n}}{\bm{\eta}}
= \sum_{n \in N} \Map{G}{A}{x}{\Inner{\gammArg{\E_{n}}{x}{\zeta_{n}}}{\eta_{n}}{\E_{n}}}
= \sum_{n \in N} \BraArg{\zeta_{n}}{\eta_{n}}.
$$
Pick a net $ \Seq{\varphi_{i}}{i \in I} $ in $ \Cc{G,\CC{0}{1}} $ converging uniformly to $ 1 $ on compact subsets of $ G $. By assumption, $ \Seq{\varphi_{i} \BraArg{\zeta_{n}}{\eta_{n}}}{i \in I} $ is Cauchy in $ \LL{2}{G,A} $ for every $ n \in N $, so the same applies to $ \Seq{\varphi_{i} \BraArg{\bm{\zeta}}{\bm{\eta}}}{i \in I} $. Hence, $ \bm{\zeta} $ is square-integrable. As the set of all elements of $ \D \bigoplus_{n = 1}^{\infty} \E_{n} $ with finitely many non-zero components is dense, we conclude that $ \D \bigoplus_{n = 1}^{\infty} \E_{n} $ is a square-integrable representation of $ \Quad{G}{A}{\alpha}{\omega} $.
\end{proof}



\begin{Lem} \label{The Fundamental Lemma of Square-Integrability}
If $ \zeta \in \Esi $, then $ \BigInner{\func{q}{\phi}}{\BBraArg{\zeta}{\eta}}{\LL{2}{G,A}} = \BigInner{\KetArg{\zeta}{\phi}}{\eta}{\E} $ for every $ \eta \in \E $ and $ \phi \in \Cc{G,A} $.
\end{Lem}

\begin{proof}
Let $ \phi \in \Cc{G,A} $, $ \eta \in \E $ and $ K \df \Supp{\phi} $. If $ \Seq{\varphi_{i}}{i \in I} $ is a net in $ \Cc{G,\CC{0}{1}} $ converging uniformly to $ 1 $ on compact subsets of $ G $, then
\begin{align*}
    \BigInner{\func{q}{\phi}}{\BBraArg{\zeta}{\eta}}{\LL{2}{G,A}}
& = \Inner{\func{q}{\phi}}{\lim_{i \in I} \func{q}{\varphi_{i} \BraArg{\zeta}{\eta}}}{\LL{2}{G,A}} \\
& = \lim_{i \in I} ~ \BigInner{\func{q}{\phi}}{\func{q}{\varphi_{i} \BraArg{\zeta}{\eta}}}{\LL{2}{G,A}} \\
& = \lim_{i \in I} \Int{G}{\func{\phi}{x}^{*} \SqBr{\func{\varphi_{i}}{x} ~ \Inner{\gammArg{\E}{x}{\zeta}}{\eta}{\E}}}{x} \\
& = \lim_{i \in I} \Int{G}{\func{\varphi_{i}}{x} ~ \func{\phi}{x}^{*} \Inner{\gammArg{\E}{x}{\zeta}}{\eta}{\E}}{x} \\
& = \lim_{i \in I} \Int{K}{\func{\varphi_{i}}{x} ~ \func{\phi}{x}^{*} \Inner{\gammArg{\E}{x}{\zeta}}{\eta}{\E}}{x} \qquad
    \Br{\text{As $ \Supp{\phi} = K $.}} \\
& = \Int{K}{\func{\phi}{x}^{*} \Inner{\gammArg{\E}{x}{\zeta}}{\eta}{\E}}{x} \qquad
    \Br{\text{As $ \varphi_{i} \rightrightarrows 1 $ on $ K $.}} \\
& = \Int{G}{\func{\phi}{x}^{*} \Inner{\gammArg{\E}{x}{\zeta}}{\eta}{\E}}{x} \\
& = \Int{G}{\Inner{\gammArg{\E}{x}{\zeta} \bullet \func{\phi}{x}}{\eta}{\E}}{x} \\
& = \Inner{\Int{G}{\gammArg{\E}{x}{\zeta} \bullet \func{\phi}{x}}{x}}{\eta}{\E} \\
& = \BigInner{\KetArg{\zeta}{\phi}}{\eta}{\E}.
\end{align*}
As $ \phi $ and $ \eta $ are arbitrary, we are finished.
\end{proof}



\begin{Prop} \label{The BBra Operator of a Square-Integrable Vector Is Bounded}
Let $ \zeta \in \Esi $. Then $ \BBra{\zeta}: \E \to \LL{2}{G,A} $ is a bounded $ A $-linear operator.
\end{Prop}

\begin{proof}
The proof of the $ A $-linearity of $ \BBra{\zeta} $ is trivial, so we omit it.

Let $ \Seq{\eta_{n}}{n \in \N{}} $ be a sequence in $ \E $ where $ \Seq{\eta_{n},\BBraArg{\zeta}{\eta_{n}}}{n \in \N{}} $ converges to some $ \Pair{\eta}{\Phi} \in \E \times \LL{2}{G,A} $. Then
\begin{align*}
\forall \phi \in \Cc{G,A}: \quad
    \Inner{\func{q}{\phi}}{\Phi}{\LL{2}{G,A}}
& = \lim_{n \to \infty} \BigInner{\func{q}{\phi}}{\BBraArg{\zeta}{\eta_{n}}}{\LL{2}{G,A}} \\
& = \lim_{n \to \infty} \BigInner{\KetArg{\zeta}{\phi}}{\eta_{n}}{\E} \qquad
    \Br{\text{By \autoref{The Fundamental Lemma of Square-Integrability}.}} \\
& = \BigInner{\KetArg{\zeta}{\phi}}{\eta}{\E} \\
& = \BigInner{\func{q}{\phi}}{\BBraArg{\zeta}{\eta}}{\LL{2}{G,A}}. \qquad
    \Br{\text{By \autoref{The Fundamental Lemma of Square-Integrability} again.}}
\end{align*}
As $ \Im{q}{\Cc{G,A}} $ is dense in $ \LL{2}{G,A} $, we obtain $ \Pair{\eta}{\Phi} \in \Graph{\BBra{\zeta}} $. Therefore, $ \BBra{\zeta} $ is  bounded by the Closed Graph Theorem.
\end{proof}



\begin{Prop} \label{The Square-Integrability of a Vector Implies the Extendibility of Its Ket Operator}
If $ \zeta \in \Esi $, then there exists a unique adjointable operator $ \KKet{\zeta}: \LL{2}{G,A} \to \E $, with adjoint $ \BBra{\zeta} $, such that $ \KKetArg{\zeta}{\func{q}{\phi}} = \KetArg{\zeta}{\phi} $ for every $ \phi \in \Cc{G,A} $.
\end{Prop}

\begin{proof}
By \autoref{The BBra Operator of a Square-Integrable Vector Is Bounded}, there exists a $ C > 0 $ such that $ \Norm{\BBraArg{\zeta}{\eta}}{\LL{2}{G,A}} \leq C \Norm{\eta}{\E} $ for every $ \eta \in \E $. Hence, by the Cauchy-Schwarz Inequality,
\begin{align*}
\forall \phi,\psi \in \Cc{G,A}: \quad
       \BigNorm{\KetArg{\zeta}{\phi - \psi}}{\E}^{2}
& =    \Norm{\BigInner{\KetArg{\zeta}{\phi - \psi}}{\KetArg{\zeta}{\phi - \psi}}{\E}}{A} \\
& =    \Norm{\BigInner{\func{q}{\phi - \psi}}{\BBraArg{\zeta}{\KetArg{\zeta}{\phi - \psi}}}{\LL{2}{G,A}}}{A} \qquad
       \Br{\text{By \autoref{The Fundamental Lemma of Square-Integrability}.}} \\
& \leq \Norm{\func{q}{\phi - \psi}}{\LL{2}{G,A}} \BigNorm{\BBraArg{\zeta}{\KetArg{\zeta}{\phi - \psi}}}{\LL{2}{G,A}} \\
& =    \Norm{\func{q}{\phi} - \func{q}{\psi}}{\LL{2}{G,A}} \BigNorm{\BBraArg{\zeta}{\KetArg{\zeta}{\phi - \psi}}}{\LL{2}{G,A}} \\
& \leq \Norm{\func{q}{\phi} - \func{q}{\psi}}{\LL{2}{G,A}} \cdot C \BigNorm{\KetArg{\zeta}{\phi - \psi}}{\E}, \quad \text{so} \\
       \BigNorm{\KetArg{\zeta}{\phi - \psi}}{\E}
& \leq C \Norm{\func{q}{\phi} - \func{q}{\psi}}{\LL{2}{G,A}}. \numberthis \label{Inequality}
\end{align*}

Let $ \Phi \in \LL{2}{G,A} $. Pick a sequence $ \Seq{\phi_{n}}{n \in \N{}} $ in $ \Cc{G,A} $ so that $ \D \lim_{n \to \infty} \func{q}{\phi_{n}} = \Phi $. By \Inequality{Inequality},
$$
\forall m,n \in \N{}: \quad
\BigNorm{\KetArg{\zeta}{\phi_{m} - \phi_{n}}}{\E} \leq C \Norm{\func{q}{\phi_{m}} - \func{q}{\phi_{n}}}{\LL{2}{G,A}}.
$$
As $ \Seq{\func{q}{\phi_{n}}}{n \in \N{}} $ is Cauchy in $ \LL{2}{G,A} $, we find that $ \Seq{\KetArg{\zeta}{\phi_{n}}}{n \in \N{}} $ is Cauchy in $ \E $. The completeness of $ \E $ implies that $ \D \lim_{n \to \infty} \KetArg{\zeta}{\phi_{n}} $ exists, and by \Inequality{Inequality} again, this limit depends only on $ \Phi $ and not on any particular choice of the sequence $ \Seq{\phi_{n}}{n \in \N{}} $, so we denote it by $ \theta_{\Phi} $.

Observe that $ \theta_{\func{q}{\phi}} = \KetArg{\zeta}{\phi} $ for every $ \phi \in \Cc{G,A} $ and that
\begin{align*}
\forall \eta \in \E: \quad
    \Inner{\theta_{\Phi}}{\eta}{\E}
& = \lim_{n \to \infty} \BigInner{\KetArg{\zeta}{\phi_{n}}}{\eta}{\E} \\
& = \lim_{n \to \infty} \BigInner{\func{q}{\phi_{n}}}{\BBraArg{\zeta}{\eta}}{\LL{2}{G,A}} \qquad
    \Br{\text{By \autoref{The Fundamental Lemma of Square-Integrability}.}} \\
& = \BigInner{\Phi}{\BBraArg{\zeta}{\eta}}{\LL{2}{G,A}}.
\end{align*}
Defining $ \KKet{\zeta}: \LL{2}{G,A} \to \E $ by $ \KKetArg{\zeta}{\Phi} \df \theta_{\Phi} $ for every $ \Phi \in \LL{2}{G,A} $, we see that $ \KKet{\zeta} $ is adjoint to $ \BBra{\zeta} $. Furthermore, $ \KKetArg{\zeta}{\func{q}{\phi}} = \theta_{\func{q}{\phi}} = \KetArg{\zeta}{\phi} $ for every $ \phi \in \Cc{G,A} $, concluding the proof.
\end{proof}


The converse of \autoref{The Square-Integrability of a Vector Implies the Extendibility of Its Ket Operator} is also true, as the next proposition shows.


\begin{Prop} \label{The Extendibility of a Ket Operator Implies the Square-Integrability of the Underlying Vector}
Let $ \zeta \in \E $. If there exists a $ T \in \AdjPair{\LL{2}{G,A}}{\E} $ so that $ \func{T}{\func{q}{\phi}} = \KetArg{\zeta}{\phi} $ for every $ \phi \in \Cc{G,A} $, then $ \zeta \in \Esi $.
\end{Prop}

\begin{proof}
Suppose that there is a $ T \in \AdjPair{\LL{2}{G,A}}{\E} $ with $ \func{T}{\func{q}{\phi}} = \KetArg{\zeta}{\phi} $ for every $ \phi \in \Cc{G,A} $. Let $ \Seq{\varphi_{i}}{i \in I} $ be a net in $ \Cc{G,\CC{0}{1}} $ converging uniformly to $ 1 $ on compact subsets of $ G $. Then for every $ \eta \in \E $, $ \phi \in \Cc{G,A} $ and $ i \in I $, we have
\begin{align*}
    \BigInner{\func{q}{\phi}}{\func{q}{\varphi_{i} \BraArg{\zeta}{\eta}}}{\LL{2}{G,A}}
& = \Int{G}{\func{\phi}{x}^{*} \SqBr{\func{\varphi_{i}}{x} ~ \Inner{\gammArg{\E}{x}{\zeta}}{\eta}{\E}}}{x} \\
& = \Int{G}{\func{\varphi_{i}}{x} ~ \func{\phi}{x}^{*} \Inner{\gammArg{\E}{x}{\zeta}}{\eta}{\E}}{x} \\
& = \Inner{\Int{G}{\gammArg{\E}{x}{\zeta} \bullet \Func{\varphi_{i} \phi}{x}}{x}}{\eta}{\E} \\
& = \BigInner{\KetArg{\zeta}{\varphi_{i} \phi}}{\eta}{\E} \\
& = \Inner{\func{T}{\func{q}{\varphi_{i} \phi}}}{\eta}{\E} \qquad \Br{\text{By assumption.}} \\
& = \Inner{\func{q}{\varphi_{i} \phi}}{\func{T^{*}}{\eta}}{\LL{2}{G,A}} \\
& = \Inner{\func{M_{\varphi_{i}}}{\func{q}{\phi}}}{\func{T^{*}}{\eta}}{\LL{2}{G,A}} \\
& = \Inner{\func{q}{\phi}}{\func{M_{\Conj{\varphi_{i}}}}{\func{T^{*}}{\eta}}}{\LL{2}{G,A}} \qquad
    \Br{\text{By \autoref{Multiplication Operators}.}} \\
& = \Inner{\func{q}{\phi}}{\func{M_{\varphi_{i}}}{\func{T^{*}}{\eta}}}{\LL{2}{G,A}}. \qquad
    \Br{\text{As $ \varphi_{i} $ is real-valued.}}
\end{align*}
As $ \Im{q}{\Cc{G,A}} $ is dense in $ \LL{2}{G,A} $, it follows that $ \func{q}{\varphi_{i} \BraArg{\zeta}{\eta}} = \func{M_{\varphi_{i}}}{\func{T^{*}}{\eta}} $ for every $ i \in I $. By \autoref{A Convergence Lemma for Multiplication Operators}, $ \D \lim_{i \in I} \func{M_{\varphi_{i}}}{\func{T^{*}}{\eta}} = \func{T^{*}}{\eta} $, so $ \Seq{\func{q}{\varphi_{i} \BraArg{\zeta}{\eta}}}{i \in I} $ is Cauchy in $ \LL{2}{G,A} $, giving $ \zeta \in \Esi $.
\end{proof}


The following proposition makes it possible to expand Meyer's framework to accommodate twisted $ C^{*} $-dynamical systems and twisted Hilbert $ C^{*} $-modules.


\begin{Prop} \label{Meyer's BBra-KKet Operators Are Morphisms of Twisted Hilbert C*-Modules}
Let $ \zeta \in \Esi $. Then $ \BBra{\zeta}: \E \to \LL{2}{G,A} $ and $ \KKet{\zeta}: \LL{2}{G,A} \to \E $ are morphisms of Hilbert $ \Quad{G}{A}{\alpha}{\omega} $-modules.
\end{Prop}

There are two ways to prove \autoref{Meyer's BBra-KKet Operators Are Morphisms of Twisted Hilbert C*-Modules}: (1) The obvious direct approach. (2) Show that either $ \BBra{\zeta} $ or $ \KKet{\zeta} $ is a $ \Hilb{G}{A}{\alpha}{\omega} $-morphism, and then apply \autoref{Morphisms of Twisted Hilbert C*-Modules Are Closed Under the Operator-Adjoint}. We prefer (2).

\begin{proof}
We already know that $ \KKet{\zeta} $ is adjointable, so it remains to prove its twisted-equivariance. Note for every $ r \in G $ and $ \phi \in \Cc{G,A} $ that
\begin{align*}
    \KKetArg{\zeta}{\GammArg{r}{\func{q}{\phi}}}
& = \KKetArg{\zeta}{\func{q}{\Map{G}{A}{x}{\om{r}{r^{-1} x}^{*} \alphArg{r}{\func{\phi}{r^{-1} x}}}}} \\
& = \KetArg{\zeta}{\Map{G}{A}{x}{\om{r}{r^{-1} x}^{*} \alphArg{r}{\func{\phi}{r^{-1} x}}}} \\
& = \Int{G}{\gammArg{\E}{x}{\zeta} \bullet \om{r}{r^{-1} x}^{*} \alphArg{r}{\func{\phi}{r^{-1} x}}}{x} \\
& = \Int{G}{\gammArg{\E}{r x}{\zeta} \bullet \om{r}{x}^{*} \alphArg{r}{\func{\phi}{x}}}{x} \qquad
    \Br{\text{By the change of variables $ x \mapsto r x $.}} \\
& = \Int{G}{\gammArg{\E}{r}{\gammArg{\E}{x}{\zeta}} \bullet \om{r}{x} ~ \om{r}{x}^{*} \alphArg{r}{\func{\phi}{x}}}{x} \\
& = \Int{G}{\gammArg{\E}{r}{\gammArg{\E}{x}{\zeta}} \bullet \alphArg{r}{\func{\phi}{x}}}{x} \\
& = \Int{G}{\gammArg{\E}{r}{\gammArg{\E}{x}{\zeta} \bullet \func{\phi}{x}}}{x} \\
& = \gammArg{\E}{r}{\Int{G}{\gammArg{\E}{x}{\zeta} \bullet \func{\phi}{x}}{x}} \qquad \Br{\text{As $ \gamm{\E}{r} $ is continuous.}} \\
& = \gammArg{\E}{r}{\KetArg{\zeta}{\phi}} \\
& = \gammArg{\E}{r}{\KKetArg{\zeta}{\func{q}{\phi}}},
\end{align*}
so by continuity, $ \KKetArg{\zeta}{\GammArg{r}{\Phi}} = \gammArg{\E}{r}{\KKetArg{\zeta}{\Phi}} $ for every $ \Phi \in \LL{2}{G,A} $. The twisted-equivariance of $ \KKet{\zeta} $ is therefore established, and by \autoref{Morphisms of Twisted Hilbert C*-Modules Are Closed Under the Operator-Adjoint}, the proof is complete.
\end{proof}


\subsection{A Complete Norm on $ \Esi $} 

We can equip $ \Esi $ with a special norm $ \Norm{\cdot}{\E,\si} $ defined by
\begin{equation} \label{S.i.-Norm}
\forall \zeta \in \Esi: \quad
\Norm{\zeta}{\E,\si} \df \Norm{\zeta}{\E} + \BigNorm{\KKet{\zeta}}{\AdjEqPair{\LL{2}{G,A}}{\E}}.
\end{equation}
As $ \BBra{\zeta}^{*} = \KKet{\zeta} $ for every $ \zeta \in \Esi $, it follows from straightforward norm arguments that
$$
\Norm{\zeta}{\E,\si} = \Norm{\zeta}{\E} + \BigNorm{\BBraKKet{\zeta}{\zeta}}{\AdjEq{\LL{2}{G,A}}}^{\frac{1}{2}}.
$$


\begin{Prop} \label{The S.i.-Norm Is Complete}
$ \Pair{\Esi}{\Norm{\cdot}{\E,\si}} $ is a Banach space.
\end{Prop}

\begin{proof}
Let $ \Seq{\zeta_{n}}{n \in \N{}} $ be a Cauchy sequence in $ \Pair{\Esi}{\Norm{\cdot}{\E,\si}} $, so that it is Cauchy in $ \E $ and $ \Seq{\KKet{\zeta_{n}}}{n \in \N{}} $ is Cauchy in $ \AdjEqPair{\LL{2}{G,A}}{\E} $. Let $ \D \zeta \df \lim_{n \to \infty} \zeta_{n} $ and $ \D T \df \lim_{n \to \infty} \KKet{\zeta_{n}} $. Then for every $ \eta \in E $ and $ \phi \in \Cc{G,A} $, we have
\begin{align*}
    \BigInner{\KetArg{\zeta}{\phi}}{\eta}{\E}
& = \lim_{n \to \infty} \BigInner{\KetArg{\zeta_{n}}{\phi}}{\eta}{\E} \qquad
    \Br{\text{By \Inequality{Basic Norm Inequality 2} and the Cauchy-Schwarz Inequality.}} \\
& = \lim_{n \to \infty} \BigInner{\KKetArg{\zeta_{n}}{\func{q}{\phi}}}{\eta}{\E} \\
& = \Inner{\func{T}{\func{q}{\phi}}}{\eta}{\E}, \quad \text{resulting in} \\
    \KetArg{\zeta}{\phi}
& = \func{T}{\func{q}{\phi}}.
\end{align*}
By \autoref{The Extendibility of a Ket Operator Implies the Square-Integrability of the Underlying Vector}, $ \zeta \in \Esi $, which implies that $ T = \KKet{\zeta} $ and $ \D \lim_{n \to \infty} \Norm{\zeta_{n} - \zeta}{\E,\si} = 0 $. Therefore, $ \Seq{\zeta_{n}}{n \in \N{}} $ has a limit in $ \Pair{\Esi}{\Norm{\cdot}{\E,\si}} $, and we are done.
\end{proof}


When the group $ G $ in $ \Quad{G}{A}{\alpha}{\omega} $ is compact, every element of $ \E $ is square-integrable and there is really no topological difference between $ \Norm{\cdot}{\E} $ and $ \Norm{\cdot}{\E,\si} $.


\begin{Prop} \label{The Compactness of G Implies the Equivalence of the Hilbert-Module Norm and the S.i.-Norm}
Suppose that the group $ G $ in $ \Quad{G}{A}{\alpha}{\omega} $ is compact. Then $ \Esi = \E $ and
$$
     \Norm{\cdot}{\E}
\leq \Norm{\cdot}{\E,\si}
\leq \SqBr{1 + \sqrt{\func{\mu}{G}}} \Norm{\cdot}{\E}.
$$
In other words, $ \Norm{\cdot}{\E} $ and $ \Norm{\cdot}{\E,\si} $ are equivalent norms on $ \E $.
\end{Prop}

\begin{proof}
Let $ \zeta \in \E $. Pick an $ \eta \in \E $ and a net $ \Seq{\varphi_{i}}{i \in I} $ in $ \Cc{G,\CC{0}{1}} $ that converges uniformly to $ 1 $ on compact subsets of $ G $. As $ G $ is compact, $ \Seq{\varphi_{i}}{i \in I} $ converges uniformly to the constant function $ \bm{1} $. Hence, as $ \BraArg{\zeta}{\eta} $ is bounded, $ \Seq{\varphi_{i} \BraArg{\zeta}{\eta}}{i \in I} $ converges uniformly to $ \BraArg{\zeta}{\eta} $.

Observe that $ \BraArg{\zeta}{\eta} \in \Cc{G,A} $ and that $ \Seq{\func{q}{\varphi_{i} \BraArg{\zeta}{\eta}}}{i \in I} $ converges in $ \LL{2}{G,A} $ to $ \func{q}{\BraArg{\zeta}{\eta}} $, which makes it a Cauchy net in $ \LL{2}{G,A} $. The first assertion is clear. To prove the latter, let $ \epsilon > 0 $. For every $ i \in I $,
\begin{align*}
           \Abs{\func{\varphi_{i}}{x} - 1}^{2} \FUNC{\BraArg{\zeta}{\eta}}{x}^{*} \FUNC{\BraArg{\zeta}{\eta}}{x}
& \leq_{A} \Norm{\varphi_{i} - \bm{1}}{\infty}^{2} \FUNC{\BraArg{\zeta}{\eta}}{x}^{*} \FUNC{\BraArg{\zeta}{\eta}}{x},
           \quad \text{whence} \\
           \Int{G}{\Abs{\func{\varphi_{i}}{x} - 1}^{2} \FUNC{\BraArg{\zeta}{\eta}}{x}^{*} \FUNC{\BraArg{\zeta}{\eta}}{x}}{x}
& \leq_{A} \Norm{\varphi_{i} - \bm{1}}{\infty}^{2} \Int{G}{\FUNC{\BraArg{\zeta}{\eta}}{x}^{*} \FUNC{\BraArg{\zeta}{\eta}}{x}}{x}.
\end{align*}
Pick an $ i_{0} \in I $ so that for every $ i \in I_{\geq i_{0}} $,
\begin{align*}
       \Norm{\bm{1} - \varphi_{i}}{\infty}
& \leq \dfrac{\epsilon}{1 + \Norm{\func{q}{\BraArg{\zeta}{\eta}}}{\LL{2}{G,A}}}, \quad \text{in which case} \\
       \Norm{\func{q}{\varphi_{i} \BraArg{\zeta}{\eta}} - \func{q}{\BraArg{\zeta}{\eta}}}{\LL{2}{G,A}}
& =    \Norm{\func{q}{\Br{\varphi_{i} - \bm{1}} \BraArg{\zeta}{\eta}}}{\LL{2}{G,A}} \\
& =    \NNorm{\Map{G}{A}{x}{\FUNC{\Br{\varphi_{i} - \bm{1}} \BraArg{\zeta}{\eta}}{x}}} \\
& =    \Norm{
            \Int{G}{\FUNC{\Br{\varphi_{i} - \bm{1}} \BraArg{\zeta}{\eta}}{x}^{*} \FUNC{\Br{\varphi_{i} - \bm{1}} \BraArg{\zeta}{\eta}}{x}}{x}
            }{A}^{\frac{1}{2}} \\
& =    \Norm{
            \Int{G}{\Abs{\func{\varphi_{i}}{x} - 1}^{2} \FUNC{\BraArg{\zeta}{\eta}}{x}^{*} \FUNC{\BraArg{\zeta}{\eta}}{x}}{x}
            }{A}^{\frac{1}{2}} \\
& \leq \Norm{
            \Norm{\varphi_{i} - \bm{1}}{\infty}^{2} \Int{G}{\FUNC{\BraArg{\zeta}{\eta}}{x}^{*} \FUNC{\BraArg{\zeta}{\eta}}{x}}{x}
            }{A}^{\frac{1}{2}} \\
& =    \Norm{\varphi_{i} - \bm{1}}{\infty}
       \Norm{\Int{G}{\FUNC{\BraArg{\zeta}{\eta}}{x}^{*} \FUNC{\BraArg{\zeta}{\eta}}{x}}{x}}{A}^{\frac{1}{2}} \\
& =    \Norm{\varphi_{i} - \bm{1}}{\infty} \NNorm{\BraArg{\zeta}{\eta}} \\
& =    \Norm{\varphi_{i} - \bm{1}}{\infty} \Norm{\func{q}{\BraArg{\zeta}{\eta}}}{\LL{2}{G,A}} \\
& \leq \frac{\epsilon}{1 + \Norm{\func{q}{\BraArg{\zeta}{\eta}}}{\LL{2}{G,A}}} \Norm{\func{q}{\BraArg{\zeta}{\eta}}}{\LL{2}{G,A}} \\
& <    \epsilon.
\end{align*}
As $ \epsilon $ is arbitrary, we get $ \D \lim_{i \in I} \func{q}{\varphi_{i} \BraArg{\zeta}{\eta}} = \func{q}{\BraArg{\zeta}{\eta}} $. Hence, $ \zeta \in \Esi $, and as $ \zeta $ is arbitrary, $ \Esi = \E $.

By the definition of $ \Norm{\cdot}{\E,\si} $, we have $ \Norm{\cdot}{\E} \leq \Norm{\cdot}{\E,\si} $, so only the second half of the inequality is non-trivial. Let $ \zeta \in \E $ as before and $ \phi \in \Cc{G,A} $. Then for every $ \eta \in \E $,
\begin{align*}
     & ~ \Norm{\BigInner{\KKetArg{\zeta}{\func{q}{\phi}}}{\eta}{\E}}{A} \\
=    & ~ \Norm{\BigInner{\KetArg{\zeta}{\phi}}{\eta}{\E}}{A} \\
=    & ~ \Norm{\Inner{\Int{G}{\gammArg{\E}{x}{\zeta} \bullet \func{\phi}{x}}{x}}{\eta}{\E}}{A} \\
=    & ~ \Norm{\Int{G}{\Inner{\gammArg{\E}{x}{\zeta} \bullet \func{\phi}{x}}{\eta}{\E}}{x}}{A} \\
=    & ~ \Norm{\Int{G}{\func{\phi}{x}^{*} \Inner{\gammArg{\E}{x}{\zeta}}{\eta}{\E}}{x}}{A} \\
=    & ~ \Norm{\Inner{\func{q}{\phi}}{\func{q}{\BraArg{\zeta}{\eta}}}{\LL{2}{G,A}}}{A} \qquad
         \Br{\text{Note that $ \BraArg{\zeta}{\eta} \in \Cc{G,A} $.}} \\
\leq & ~ \Norm{\func{q}{\phi}}{\LL{2}{G,A}} \Norm{\func{q}{\BraArg{\zeta}{\eta}}}{\LL{2}{G,A}} \qquad
         \Br{\text{By the Cauchy-Schwarz Inequality.}} \\
=    & ~ \Norm{\func{q}{\phi}}{\LL{2}{G,A}}
         \Norm{\Int{G}{\FUNC{\BraArg{\zeta}{\eta}}{x}^{*} \FUNC{\BraArg{\zeta}{\eta}}{x}}{x}}{A}^{\frac{1}{2}} \\
\leq & ~ \Norm{\func{q}{\phi}}{\LL{2}{G,A}} \SqBr{\Int{G}{\Norm{\FUNC{\BraArg{\zeta}{\eta}}{x}}{A}^{2}}{x}}^{\frac{1}{2}} \\
=    & ~ \Norm{\func{q}{\phi}}{\LL{2}{G,A}} \SqBr{\Int{G}{\Norm{\Inner{\gammArg{\E}{x}{\zeta}}{\eta}{\E}}{A}^{2}}{x}}^{\frac{1}{2}} \\
\leq & ~ \Norm{\func{q}{\phi}}{\LL{2}{G,A}} \SqBr{\Int{G}{\Norm{\gammArg{\E}{x}{\zeta}}{\E}^{2} \Norm{\eta}{\E}^{2}}{x}}^{\frac{1}{2}} \qquad
         \Br{\text{By the Cauchy-Schwarz Inequality again.}} \\
=    & ~ \Norm{\func{q}{\phi}}{\LL{2}{G,A}} \SqBr{\Int{G}{\Norm{\zeta}{\E}^{2} \Norm{\eta}{\E}^{2}}{x}}^{\frac{1}{2}} \qquad
         \Br{\text{As $ \gamm{\E}{x} \in \Isom{\E} $ for every $ x \in G $.}} \\
=    & ~ \Norm{\func{q}{\phi}}{\LL{2}{G,A}} \SqBr{\Norm{\zeta}{\E}^{2} \Norm{\eta}{\E}^{2} \func{\mu}{G}}^{\frac{1}{2}} \\
=    & ~ \Norm{\func{q}{\phi}}{\LL{2}{G,A}} \Norm{\zeta}{\E} \Norm{\eta}{\E} \sqrt{\func{\mu}{G}}.
\end{align*}
Letting $ \eta = \KKetArg{\zeta}{\func{q}{\phi}} $ in the foregoing derivation gives us
$$
     \BigNorm{\KKetArg{\zeta}{\func{q}{\phi}}}{\E}^{2}
=    \Norm{\BigInner{\KKetArg{\zeta}{\func{q}{\phi}}}{\KKetArg{\zeta}{\func{q}{\phi}}}{\E}}{A}
\leq \Norm{\func{q}{\phi}}{\LL{2}{G,A}} \Norm{\zeta}{\E} \BigNorm{\KKetArg{\zeta}{\func{q}{\phi}}}{\E} \sqrt{\func{\mu}{G}}.
$$
Hence, $ \BigNorm{\KKetArg{\zeta}{\func{q}{\phi}}}{\E} \leq \Norm{\func{q}{\phi}}{\LL{2}{G,A}} \Norm{\zeta}{\E} \sqrt{\func{\mu}{G}} $, and as $ \phi $ is arbitrary, we obtain
$$
\forall \Phi \in \LL{2}{G,A}: \quad
\BigNorm{\KKetArg{\zeta}{\Phi}}{\E} \leq \Norm{\Phi}{\LL{2}{G,A}} \Norm{\zeta}{\E} \sqrt{\func{\mu}{G}}.
$$
Therefore, $ \BigNorm{\KKet{\zeta}}{\AdjEqPair{\LL{2}{G,A}}{\E}} \leq \Norm{\zeta}{\E} \sqrt{\func{\mu}{G}} $, so
$$
     \Norm{\zeta}{\E,\si}
\df  \Norm{\zeta}{\E} + \BigNorm{\KKet{\zeta}}{\AdjEqPair{\LL{2}{G,A}}{\E}}
\leq \Norm{\zeta}{\E} + \Norm{\zeta}{\E} \sqrt{\func{\mu}{G}}
=    \SqBr{1 + \sqrt{\func{\mu}{G}}} \Norm{\zeta}{\E}.
$$
As $ \zeta $ is arbitrary, we are finished.
\end{proof}



\section{Reduced Twisted Crossed Products} 

Here, we present the reduced twisted crossed product for $ \Quad{G}{A}{\alpha}{\omega} $ as a certain $ C^{*} $-subalgebra of $ \AdjEq{\LL{2}{G,A}} $. Some of the material has been sourced from the papers \cite{Busby|Smith,Packer|Raeburn}.

\textbf{In this section, $ \E $ and $ \F $ are Hilbert $ \Quad{G}{A}{\alpha}{\omega} $-modules.}

\subsection{Covariant Representations of a Twisted $ C^{*} $-Dynamical System} 

Define operations $ \star: \Cc{G,A} \times \Cc{G,A} \to \Cc{G,A} $ and $ ^{*}: \Cc{G,A} \to \Cc{G,A} $ as follows:
\begin{align*}
\forall f,g \in \Cc{G,A}: \quad
f \star g & \df \Map{G}{A}{x}{\Int{G}{\func{f}{y} ~ \alphArg{y}{\func{g}{y^{-1} x}} ~ \om{y}{y^{-1} x}}{y}}, \qquad \Br{\text{Convolution}} \\
f^{*}     & \df \Map{G}{A}{x}{\Del{x}{-1} \om{x}{x^{-1}}^{*} \alphArg{x}{\func{f}{x^{-1}}}^{*}}.             \qquad \Br{\text{Involution}}
\end{align*}
The quadruple $ \Quad{\Cc{G,A}}{\star}{^{*}}{\Norm{\cdot}{1}} $ is then a normed $ * $-algebra, with $ ^{*} $ as an isometric involution.

Next, a \emph{covariant representation} of $ \Quad{G}{A}{\alpha}{\omega} $ is a triple $ \Trip{\X}{\pi}{U} $, where:
\begin{itemize}
\item
$ \X $ is a Hilbert $ B $-module for some $ C^{*} $-algebra $ B $.

\item
$ \pi $ is a $ * $-representation of $ A $ on $ \X $.

\item
$ U $ is a strongly continuous map from $ G $ to $ \U{\Adj{\X}} $.

\item
$ \func{\pi}{\alphArg{r}{a}} = \func{U}{r} \circ \func{\pi}{a} \circ \func{U}{r}^{*} $ for every $ r,s \in G $.

\item
$ \func{U}{r} \circ \func{U}{s} = \func{\Conj{\pi}}{\om{r}{s}} \circ \func{U}{r s} $ for every $ r,s \in G $.
\end{itemize}
For each covariant representation $ \Trip{\X}{\pi}{U} $ of $ \Quad{G}{A}{\alpha}{\omega} $, we can define a $ * $-algebra homomorphism $ \rho_{\X,\pi,U}: \Trip{\Cc{G,A}}{\star}{^{*}} \to \Trip{\Adj{\X}}{\circ}{^{*}} $ by
$$
\forall f \in \Cc{G,A}: \quad
\func{\rho_{\X,\pi,U}}{f} \df \Map{\X}{\X}{x}{\Int{G}{\FUNC{\func{\pi}{\func{f}{y}} \circ \func{U}{y}}{x}}{y}},
$$
and it is not hard to show that $ \Norm{\func{\rho_{\X,\pi,U}}{f}}{\Adj{\X}} \leq \Norm{f}{1} $ for every $ f \in \Cc{G,A} $.

If $ \pi: A \to \AdjEq{\LL{2}{G,A}} $ and $ \lambda: G \to \U{\AdjEq{\LL{2}{G,A}}} $ are the maps given in \autoref{Covariant Representation Part 1} and \autoref{Covariant Representation Part 2} respectively, then $ \Trip{\LL{2}{G,A}}{\pi}{\lambda} $ is clearly a covariant representation of $ \Quad{G}{A}{\alpha}{\omega} $. If $ \rho \df \rho_{\LL{2}{G,A},\pi,\lambda} $, then the following statements are true:
\begin{itemize}
\item
$ \func{\rho}{f} \in \AdjEq{\LL{2}{G,A}} $ for every $ f \in \Cc{G,A} $, as $ \func{\pi}{\func{f}{y}} \circ \func{\lambda}{y} \in \AdjEq{\LL{2}{G,A}} $ for every $ y \in G $.

\item
$ \Norm{\func{\rho}{f}}{\AdjEq{\LL{2}{G,A}}} \leq \Norm{f}{1} $ for every $ f \in \Cc{G,A} $.

\item
For every $ f \in \Cc{G,A} $, we have
\begin{equation}
\forall \phi \in \Cc{G,A}: \quad
  \FUNC{\func{\rho}{f}}{\func{q}{\phi}}
= \func{q}{\Map{G}{A}{x}{\Int{G}{\Del{y}{\frac{1}{2}} \alphArg{x}{\func{f}{y}} ~ \om{x}{y} ~ \func{\phi}{x y}}{y}}},
  \label{An Explicit Formula for the Integrated Form of the Covariant Representation of G and A}
\end{equation}
so $ \Im{q}{\Cc{G,A}} $ is invariant under $ \func{\rho}{f} $.

\item
The set $ \Span{\Set{\FUNC{\func{\rho}{f}}{\func{q}{\phi}}}{f,\phi \in \Cc{G,A}}} $ is dense in $ \LL{2}{G,A} $, which makes $ \rho $ non-degenerate. (See Proposition 2.23 of \cite{Williams}.)

\item
$ \rho $ is injective, as $ \pi $ is injective. (See Lemma 2.26 of \cite{Williams}.)
\end{itemize}


\begin{Def} \label{Reduced Twisted Crossed Products}
Define the \emph{reduced twisted crossed product} for the twisted $ C^{*} $-dynamical system $ \Quad{G}{A}{\alpha}{\omega} $ as the $ C^{*} $-algebra $ \Cl{\Range{\rho}}{\AdjEq{\LL{2}{G,A}}} $, and denote it by $ \RTCP{G}{A}{\alpha}{\omega} $.
\end{Def}



\begin{Rmk} \label{Some Remarks on Reduced Twisted Crossed Products}
When $ \omega $ is trivial, this agrees with the earlier definition of a reduced crossed product. Although we will not work with it here, we can define the full twisted crossed product for $ \Quad{G}{A}{\alpha}{\omega} $ as the completion of $ \Trip{\Cc{G,A}}{\star}{^{*}} $ with respect to the norm $ \Norm{\cdot}{u} $ defined by
$$
\forall f \in \Cc{G,A}: \quad
    \Norm{f}{u}
\df \func{\sup}
    {\Set{\Norm{\func{\rho_{\X,\pi,U}}{f}}{\Adj{\X}}}{\text{$ \Trip{\X}{\pi}{U} $ is a covariant rep. of $ \Quad{G}{A}{\alpha}{\omega} $}}}.
$$
That $ \Norm{\cdot}{u} $ is not merely a semi-norm is due to the injectivity of $ \rho $.
\end{Rmk}


Let us now insert a lemma that relates Meyer's bra-ket operators to the maps $ \pi $ and $ \lambda $.


\begin{Lem} \label{Some Ket Identities}
For every $ T \in \AdjEqPair{\E}{\F} $, $ r \in G $, $ a \in A $ and $ \zeta \in \Esi $, the following ket identities hold:
\begin{align}
\Ket{\func{T}{\zeta}}        & = T \circ \Ket{\zeta} \circ q. \label{Ket Identity 1} \\
\Ket{\zeta \bullet a}        & = \KKet{\zeta} \circ \func{\pi}{a} \circ q. \label{Ket Identity 2} \\
\Ket{\gammArg{\E}{r}{\zeta}} & = \Del{r}{- \frac{1}{2}} \SqBr{\KKet{\zeta} \circ \func{\lambda}{r}^{*} \circ q}.
                                 \label{Ket Identity 3}
\end{align}
\end{Lem}

\begin{proof}
For every $ T \in \AdjEqPair{\E}{\F} $, $ \zeta \in \Esi $ and $ \phi \in \Cc{G,A} $, we have
\begin{align*}
    \KetArg{\func{T}{\zeta}}{\phi}
& = \Int{G}{\gammArg{\F}{x}{\func{T}{\zeta}} \bullet \func{\phi}{x}}{x} \\
& = \Int{G}{\func{T}{\gammArg{\E}{x}{\zeta}} \bullet \func{\phi}{x}}{x} \qquad \Br{\text{As $ T $ is twisted-equivariant.}} \\
& = \Int{G}{\func{T}{\gammArg{\E}{x}{\zeta} \bullet \func{\phi}{x}}}{x} \qquad \Br{\text{As $ T $ is $ A $-linear.}} \\
& = \func{T}{\Int{G}{\gammArg{\E}{x}{\zeta} \bullet \func{\phi}{x}}{x}} \qquad \Br{\text{As $ T $ is continuous.}} \\
& = \func{T}{\KetArg{\zeta}{\phi}} \\
& = \func{T}{\KKetArg{\zeta}{\func{q}{\phi}}}, \quad \text{so} \\
    \Ket{\func{T}{\zeta}}
& = T \circ \KKet{\zeta} \circ q.
\end{align*}

For every $ a \in A $, $ \zeta \in \Esi $ and $ \phi \in \Cc{G,A} $, we have
\begin{align*}
    \KetArg{\zeta \bullet a}{\phi}
& = \Int{G}{\gammArg{\E}{x}{\zeta \bullet a} \bullet \func{\phi}{x}}{x} \\
& = \Int{G}{\SqBr{\gammArg{\E}{x}{\zeta} \bullet \alphArg{x}{a}} \bullet \func{\phi}{x}}{x} \\
& = \Int{G}{\gammArg{\E}{x}{\zeta} \bullet \alphArg{x}{a} ~ \func{\phi}{x}}{x} \\
& = \KetArg{\zeta}{\Map{G}{A}{x}{\alphArg{x}{a} ~ \func{\phi}{x}}} \\
& = \KKetArg{\zeta}{\func{q}{\Map{G}{A}{x}{\alphArg{x}{a} ~ \func{\phi}{x}}}} \\
& = \KKetArg{\zeta}{\FUNC{\func{\pi}{a}}{\func{q}{\phi}}}, \quad \text{so} \\
    \Ket{\zeta \bullet a}
& = \KKet{\zeta} \circ \func{\pi}{a} \circ q.
\end{align*}

For every $ r \in G $, $ \zeta \in \Esi $ and $ \phi \in \Cc{G,A} $, we have
\begin{align*}
    \KetArg{\gammArg{\E}{r}{\zeta}}{\phi}
& = \Int{G}{\gammArg{\E}{x}{\gammArg{\E}{r}{\zeta}} \bullet \func{\phi}{x}}{x} \\
& = \Int{G}{\SqBr{\gammArg{\E}{x r}{\zeta} \bullet \om{x}{r}^{*}} \bullet \func{\phi}{x}}{x} \\
& = \Int{G}{\gammArg{\E}{x r}{\zeta} \bullet \om{x}{r}^{*} \func{\phi}{x}}{x} \\
& = \Del{r^{-1}}{} \Int{G}{\gammArg{\E}{x}{\zeta} \bullet \om{x r^{-1}}{r}^{*} \func{\phi}{x r^{-1}}}{x} \qquad
    \Br{\text{By the change of variables $ x \mapsto x r^{-1} $.}} \\
& = \Del{r}{-1} \Int{G}{\gammArg{\E}{x}{\zeta} \bullet \om{x r^{-1}}{r}^{*} \func{\phi}{x r^{-1}}}{x} \\
& = \Del{r}{- \frac{1}{2}}
    \Int{G}{\gammArg{\E}{x}{\zeta} \bullet \SqBr{\Del{r}{- \frac{1}{2}} \om{x r^{-1}}{r}^{*} \func{\phi}{x r^{-1}}}}{x} \\
& = \Del{r}{- \frac{1}{2}} \KetArg{\zeta}{\Map{G}{A}{x}{\Del{r}{- \frac{1}{2}} \om{x r^{-1}}{r}^{*} \func{\phi}{x r^{-1}}}} \\
& = \Del{r}{- \frac{1}{2}} \KKetArg{\zeta}{\func{q}{\Map{G}{A}{x}{\Del{r}{- \frac{1}{2}} \om{x r^{-1}}{r}^{*} \func{\phi}{x r^{-1}}}}} \\
& = \Del{r}{- \frac{1}{2}} \KKetArg{\zeta}{\FUNC{\func{\lambda}{r}^{*}}{\func{q}{\phi}}}, \quad \text{so} \\
    \Ket{\gammArg{\E}{r}{\zeta}}
& = \Del{r}{- \frac{1}{2}} \SqBr{\KKet{\zeta} \circ \func{\lambda}{r}^{*} \circ q}.
\end{align*}
This concludes the proof.
\end{proof}


\Identity{Ket Identity 1} implies that $ \Im{T}{\Esi} \subseteq \Fsi $ for any $ T \in \AdjEqPair{\E}{\F} $. Indeed, if $ \zeta \in \Esi $, then
$$
\forall \phi \in \Cc{G,A}: \quad
\KetArg{\func{T}{\zeta}}{\phi} = \Func{T \circ \KKet{\zeta}}{\func{q}{\phi}}.
$$
As $ T \circ \KKet{\zeta} \in \AdjPair{\LL{2}{G,A}}{\E} $, \autoref{The Extendibility of a Ket Operator Implies the Square-Integrability of the Underlying Vector} implies that $ \func{T}{\zeta} \in \Fsi $.

\Identity{Ket Identity 2} implies that $ \Esi \bullet A \subseteq \Esi $. Indeed, if $ a \in A $ and $ \zeta \in \Esi $, then
$$
\forall \phi \in \Cc{G,A}: \quad
\KetArg{\zeta \bullet a}{\phi} = \FUNC{\KKet{\zeta} \circ \func{\pi}{a}}{\func{q}{\phi}}.
$$
As $ \KKet{\zeta} \circ \func{\pi}{a} \in \AdjPair{\LL{2}{G,A}}{\E} $, \autoref{The Extendibility of a Ket Operator Implies the Square-Integrability of the Underlying Vector} implies that $ \zeta \bullet a \in \Esi $.

Via the same logic, \Identity{Ket Identity 3} implies that $ \Esi $ is invariant under the twisted $ G $-action on $ \E $.


\begin{Lem} \label{Some Norm Inequalities}
For every $ T \in \AdjEqPair{\E}{\F} $, $ r \in G $, $ a \in A $ and $ \zeta \in \Esi $, the following inequalities hold:
\begin{align}
\Norm{\func{T}{\zeta}}{\F,\si}        & \leq \Norm{\zeta}{\E,\si} \Norm{T}{\AdjEqPair{\E}{\F}}. \label{Norm Inequality 1} \\
\Norm{\zeta \bullet a}{\E,\si}        & \leq \Norm{\zeta}{\E,\si} \Norm{a}{A}. \label{Norm Inequality 2} \\
\Norm{\gammArg{\E}{r}{\zeta}}{\E,\si} & \leq \Norm{\zeta}{\E,\si} \cdot \func{\max}{1,\Del{r}{- \frac{1}{2}}}.
                                        \label{Norm Inequality 3}
\end{align}
\end{Lem}

\begin{proof}
For every $ T \in \AdjEqPair{\E}{\F} $ and $ \zeta \in \Esi $, we have $ \func{T}{\zeta} \in \Fsi $, so
\begin{align*}
       \Norm{\func{T}{\zeta}}{\F,\si}
& =    \Norm{\func{T}{\zeta}}{\F} + \BigNorm{\KKet{\func{T}{\zeta}}}{\AdjEqPair{\LL{2}{G,A}}{\F}} \\
& =    \Norm{\func{T}{\zeta}}{\F} + \BigNorm{T \circ \KKet{\zeta}}{\AdjEqPair{\LL{2}{G,A}}{\F}} \qquad
       \Br{\text{By \Identity{Ket Identity 1}.}} \\
& \leq \Norm{T}{\AdjEqPair{\E}{\F}} \Norm{\zeta}{\E} + \Norm{T}{\AdjEqPair{\E}{\F}} \BigNorm{\KKet{\zeta}}{\AdjEqPair{\LL{2}{G,A}}{\F}} \\
& =    \Br{\Norm{\zeta}{\E} + \BigNorm{\KKet{\zeta}}{\AdjEqPair{\LL{2}{G,A}}{\E}}} \Norm{T}{\AdjEqPair{\E}{\F}} \\
& =    \Norm{\zeta}{\E,\si} \Norm{T}{\AdjEqPair{\E}{\F}}.
\end{align*}

For every $ a \in A $ and $ \zeta \in \Esi $, we have $ \zeta \bullet a \in \Esi $, so
\begin{align*}
       \Norm{\zeta \bullet a}{\E,\si}
& =    \Norm{\zeta \bullet a}{\E} + \BigNorm{\KKet{\zeta \bullet a}}{\AdjEqPair{\LL{2}{G,A}}{\E}} \\
& =    \Norm{\zeta \bullet a}{\E} + \BigNorm{\func{\pi}{a} \circ \KKet{\zeta}}{\AdjEqPair{\LL{2}{G,A}}{\E}} \qquad
       \Br{\text{By \Identity{Ket Identity 2}.}} \\
& \leq \Norm{\zeta}{\E} \Norm{a}{A} + \Norm{\func{\pi}{a}}{\AdjEq{\LL{2}{G,A}}} \BigNorm{\KKet{\zeta}}{\AdjEqPair{\LL{2}{G,A}}{\E}} \\
& \leq \Norm{\zeta}{\E} \Norm{a}{A} + \Norm{a}{A} \BigNorm{\KKet{\zeta}}{\AdjEqPair{\LL{2}{G,A}}{\E}} \\
& =    \Br{\Norm{\zeta}{\E} + \BigNorm{\KKet{\zeta}}{\AdjEqPair{\LL{2}{G,A}}{\E}}} \Norm{a}{A} \\
& =    \Norm{\zeta}{\E,\si} \Norm{a}{A}.
\end{align*}

For every $ r \in G $ and $ \zeta \in \Esi $, we have $ \gammArg{\E}{r}{\zeta} \in \Esi $, so
\begin{align*}
       \Norm{\gammArg{\E}{r}{\zeta}}{\E,\si}
& =    \Norm{\gammArg{\E}{r}{\zeta}}{\E} + \BigNorm{\KKet{\gammArg{\E}{r}{\zeta}}}{\AdjEqPair{\LL{2}{G,A}}{\E}} \\
& =    \Norm{\gammArg{\E}{r}{\zeta}}{\E} +
       \Norm{\Del{r}{- \frac{1}{2}} \SqBr{\KKet{\zeta} \circ \func{\lambda}{r}^{*}}}{\AdjEqPair{\LL{2}{G,A}}{\E}} \qquad
       \Br{\text{By \Identity{Ket Identity 3}.}} \\
& =    \Norm{\gammArg{\E}{r}{\zeta}}{\E} +
       \Del{r}{- \frac{1}{2}} \BigNorm{\KKet{\zeta} \circ \func{\lambda}{r}^{*}}{\AdjEqPair{\LL{2}{G,A}}{\E}} \\
& \leq \Norm{\gammArg{\E}{r}{\zeta}}{\E} +
       \Del{r}{- \frac{1}{2}} \BigNorm{\KKet{\zeta}}{\AdjEqPair{\LL{2}{G,A}}{\E}} \Norm{\func{\lambda}{r}^{*}}{\AdjEq{\LL{2}{G,A}}} \\
& =    \Norm{\zeta}{\E} + \Del{r}{- \frac{1}{2}} \BigNorm{\KKet{\zeta}}{\AdjEq{\LL{2}{G,A}}} \qquad
       \Br{\text{As $ \func{\lambda}{r} $ is unitary.}} \\
& \leq \func{\max}{1,\Del{r}{- \frac{1}{2}}} \Norm{\zeta}{\E} +
       \func{\max}{1,\Del{r}{- \frac{1}{2}}} \BigNorm{\KKet{\zeta}}{\AdjEq{\LL{2}{G,A}}} \\
& =    \Br{\Norm{\zeta}{\E} + \BigNorm{\KKet{\zeta}}{\AdjEq{\LL{2}{G,A}}}} \cdot \func{\max}{1,\Del{r}{- \frac{1}{2}}} \\
& =    \Norm{\zeta}{\E,\si} \cdot \func{\max}{1,\Del{r}{- \frac{1}{2}}}.
\end{align*}
This concludes the proof.
\end{proof}


$ \Esi $ can be given the structure of a right $ \Pair{\Cc{G,A}}{\star} $-module. In order to accomplish this, we enlist the aid of two special operators on $ \Cc{G,A} $.


\begin{Def} \label{The Sharp and Flat Operators}
Define operators $ \sharp,\flat: \Cc{G,A} \to \Cc{G,A} $ by
\begin{align*}
\forall f \in \Cc{G,A}: \quad
f^{\sharp} & \df \Map{G}{A}{x}{\Del{x}{- \frac{1}{2}} \om{x}{x^{-1}}^{*} \alphArg{x}{\func{f}{x^{-1}}}}, \\
f^{\flat}  & \df \Map{G}{A}{x}{\Del{x}{- \frac{1}{2}} \alphArg{x}{\func{f}{x^{-1}}} ~ \om{x}{x^{-1}}}.
\end{align*}
\end{Def}



\begin{Lem} \label{The Sharp and Flat Operators Are Inverses of Each Other}
The operators $ \sharp $ and $ \flat $ are inverses of each other.
\end{Lem}

\begin{proof}
For every $ f \in \Cc{G,A} $ and $ x \in G $, we have
\begin{align*}
    \func{f^{\sharp \flat}}{x}
& = \Del{x}{- \frac{1}{2}} \alphArg{x}{\Sharp{f}{x^{-1}}} ~ \om{x}{x^{-1}} \\
& = \Del{x}{- \frac{1}{2}} \alphArg{x}{\Del{x}{\frac{1}{2}} \om{x^{-1}}{x}^{*} \alphArg{x^{-1}}{\func{f}{x}}} ~ \om{x}{x^{-1}} \\
& = \alphArg{x}{\om{x^{-1}}{x}^{*} \alphArg{x^{-1}}{\func{f}{x}}} ~ \om{x}{x^{-1}} \\
& = \alphArgExt{x}{\om{x^{-1}}{x}^{*}} ~ \alphArg{x}{\alphArg{x^{-1}}{\func{f}{x}}} ~ \om{x}{x^{-1}} \\
& = \alphArgExt{x}{\om{x^{-1}}{x}}^{*} \alphArg{x}{\alphArg{x^{-1}}{\func{f}{x}}} ~ \om{x}{x^{-1}} \\
& = \alphArgExt{x}{\om{x^{-1}}{x}}^{*} \om{x}{x^{-1}} ~ \func{f}{x} ~ \om{x}{x^{-1}}^{*} \om{x}{x^{-1}} \\
& = \alphArgExt{x}{\om{x^{-1}}{x}}^{*} \om{x}{x^{-1}} ~ \func{f}{x} \\
& = \om{x}{e} ~ \om{e}{x}^{*} \func{f}{x} \\
& = \func{f}{x} \quad \text{and} \\
    \func{f^{\flat \sharp}}{x}
& = \Del{x}{- \frac{1}{2}} \om{x}{x^{-1}}^{*} \alphArg{x}{\Flat{f}{x^{-1}}} \\
& = \Del{x}{- \frac{1}{2}} \om{x}{x^{-1}}^{*} \alphArg{x}{\Del{x}{\frac{1}{2}} \alphArg{x^{-1}}{\func{f}{x}} ~ \om{x^{-1}}{x}} \\
& = \om{x}{x^{-1}}^{*} \alphArg{x}{\alphArg{x^{-1}}{\func{f}{x}} ~ \om{x^{-1}}{x}} \\
& = \om{x}{x^{-1}}^{*} \alphArg{x}{\alphArg{x^{-1}}{\func{f}{x}}} ~ \alphArgExt{x}{\om{x^{-1}}{x}} \\
& = \om{x}{x^{-1}}^{*} \om{x}{x^{-1}} ~ \func{f}{x} ~ \om{x}{x^{-1}}^{*} \alphArgExt{x}{\om{x^{-1}}{x}} \\
& = \func{f}{x} ~ \om{x}{x^{-1}}^{*} \alphArgExt{x}{\om{x^{-1}}{x}} \\
& = \func{f}{x} ~ \om{e}{x} ~ \om{x}{e}^{*} \\
& = \func{f}{x}.
\end{align*}
The proof is now complete.
\end{proof}



\begin{Rmk} \label{The Sharp and Flat Operators Coincide for a C*-Dynamical System}
When $ \omega $ is trivial, $ \sharp $ and $ \flat $ are the same operator, which Meyer denotes by $ ~ \check{} ~ $ in \cite{Meyer2}.
\end{Rmk}



\begin{Lem} \label{The Sharp and Flat Operators Preserve L2-Norms}
We have $ \Norm{f^{\sharp}}{2} = \Norm{f^{\flat}}{2} = \Norm{f}{2} $ for every $ f \in \Cc{G,A} $.
\end{Lem}

\begin{proof}
Let $ f \in \Cc{G,A} $. Then
\begin{align*}
\forall x \in G: \quad
    \Norm{\Sharp{f}{x}}{A}
& = \Norm{\Del{x}{- \frac{1}{2}} \om{x}{x^{-1}}^{*} \alphArg{x}{\func{f}{x^{-1}}}}{A} \\
& = \Del{x}{- \frac{1}{2}} \Norm{\om{x}{x^{-1}}^{*} \alphArg{x}{\func{f}{x^{-1}}}}{A} \\
& = \Del{x}{- \frac{1}{2}} \Norm{\alphArg{x}{\func{f}{x^{-1}}}}{A} \qquad \Br{\text{As $ \om{x}{x^{-1}} \in \UMult{A} $.}} \\
& = \Del{x}{- \frac{1}{2}} \Norm{\func{f}{x^{-1}}}{A}.
\end{align*}
Similarly, $ \Norm{\Flat{f}{x}}{A} = \Del{x}{- \frac{1}{2}} \Norm{\func{f}{x^{-1}}}{A} $ for every $ x \in G $. Hence,
$$
  \Norm{f^{\sharp}}{2}
= \Norm{f^{\flat}}{2}
= \SqBr{\Int{G}{\Del{x^{-1}}{} \Norm{\func{f}{x^{-1}}}{A}^{2}}{x}}^{\frac{1}{2}}
= \SqBr{\Int{G}{\Norm{\func{f}{x}}{A}^{2}}{x}}^{\frac{1}{2}}
= \Norm{f}{2},
$$
and as $ f $ is arbitrary, we are done.
\end{proof}



\begin{Lem} \label{The Workhorse Lemma}
Let $ f,\phi \in \Cc{G,A} $. Then the integral $ \D \Int{G}{\GammArg{x}{\func{q}{f^{\flat}}} \centerdot \func{\phi}{x}}{x} $ converges in $ \LL{2}{G,A} $ and is equal to $ \FUNC{\func{\rho}{f}}{\func{q}{\phi}} $.
\end{Lem}

\begin{proof}
The integral converges because $ \Map{G}{\LL{2}{G,A}}{x}{\GammArg{x}{\func{q}{f^{\flat}}} \centerdot \func{\phi}{x}} $ is continuous and compactly supported. Knowing that it is well-defined, we have for every $ \psi \in \Cc{G,A} $ that
\begin{align*}
  & ~ \Inner{\func{q}{\psi}}{\Int{G}{\GammArg{x}{\func{q}{f^{\flat}}} \centerdot \func{\phi}{x}}{x}}{\LL{2}{G,A}} \\
= & ~ \Int{G}{\Inner{\func{q}{\psi}}{\GammArg{x}{\func{q}{f^{\flat}}} \centerdot \func{\phi}{x}}{\LL{2}{G,A}}}{x} \qquad
      \Br{\text{As $ \Inner{\cdot}{\cdot}{\LL{2}{G,A}} $ is continuous.}} \\
= & ~ \Int{G}{\SqBr{\Int{G}{\func{\psi}{y}^{*} \om{x}{x^{-1} y}^{*} \alphArg{x}{\func{f^{\flat}}{x^{-1} y}} ~ \func{\phi}{x}}{y}}}{x} \\
= & ~ \Int{G}{
             \SqBr{
                  \Int{G}{
                         \func{\psi}{y}^{*}
                         \om{x}{x^{-1} y}^{*}
                         \alphArg{x}{
                                    \Del{x^{-1} y}{- \frac{1}{2}}
                                    \alphArg{x^{-1} y}{\func{f}{\Br{x^{-1} y}^{-1}}} ~
                                    \om{x^{-1} y}{\Br{x^{-1} y}^{-1}}
                                    } ~
                         \func{\phi}{x}
                         }{y}
                  }
             }{x} \\
= & ~ \Int{G}{
             \SqBr{
                  \Int{G}{
                         \func{\psi}{y}^{*}
                         \om{x}{x^{-1} y}^{*}
                         \alphArg{x}{\Del{y^{-1} x}{\frac{1}{2}} \alphArg{x^{-1} y}{\func{f}{y^{-1} x}} ~ \om{x^{-1} y}{y^{-1} x}} ~
                         \func{\phi}{x}
                         }{y}
                  }
             }{x} \\
= & ~ \Int{G}{
             \SqBr{
                  \Int{G}{
                         \Del{y^{-1} x}{\frac{1}{2}}
                         \func{\psi}{y}^{*}
                         \om{x}{x^{-1} y}^{*}
                         \alphArg{x}{\alphArg{x^{-1} y}{\func{f}{y^{-1} x}} ~ \om{x^{-1} y}{y^{-1} x}} ~
                         \func{\phi}{x}
                         }{y}
                  }
             }{x} \\
= & ~ \Int{G}{
             \SqBr{
                  \Int{G}{
                         \Del{y^{-1} x}{\frac{1}{2}}
                         \func{\psi}{y}^{*}
                         \om{x}{x^{-1} y}^{*}
                         \alphArg{x}{\alphArg{x^{-1} y}{\func{f}{y^{-1} x}}} ~
                         \alphArgExt{x}{\om{x^{-1} y}{y^{-1} x}} ~
                         \func{\phi}{x}
                         }{y}
                  }
             }{x} \\
= & ~ \Int{G}{
             \SqBr{
                  \Int{G}{
                         \Del{y^{-1} x}{\frac{1}{2}}
                         \func{\psi}{y}^{*}
                         \om{x}{x^{-1} y}^{*}
                         \alphArg{x}{\alphArg{x^{-1} y}{\func{f}{y^{-1} x}}} ~
                         \alphArgExt{x}{\om{x^{-1} y}{y^{-1} x}} ~
                         \func{\phi}{x}
                         }{x}
                  }
             }{y} \\
  & ~ \Br{\text{By Fubini's Theorem.}} \\
= & ~ \Int{G}{
             \SqBr{
                  \Int{G}{
                         \Del{x}{\frac{1}{2}}
                         \func{\psi}{y}^{*}
                         \om{y x}{x^{-1}}^{*}
                         \alphArg{y x}{\alphArg{x^{-1}}{\func{f}{x}}} ~
                         \alphArgExt{y x}{\om{x^{-1}}{x}} ~
                         \func{\phi}{y x}
                         }{x}
                  }
             }{y} \\
  & ~ \Br{\text{By the change of variables $ x \mapsto y x $.}} \\
= & ~ \Int{G}{
             \SqBr{
                  \Int{G}{
                         \Del{x}{\frac{1}{2}}
                         \func{\psi}{y}^{*}
                         \om{y x}{x^{-1}}^{*}
                         \om{y x}{x^{-1}} ~
                         \alphArg{y}{\func{f}{x}} ~
                         \om{y x}{x^{-1}}^{*}
                         \alphArgExt{y x}{\om{x^{-1}}{x}} ~
                         \func{\phi}{y x}
                         }{x}
                  }
             }{y} \\
= & ~ \Int{G}{
             \SqBr{
                  \Int{G}{
                         \Del{x}{\frac{1}{2}}
                         \func{\psi}{y}^{*}
                         \alphArg{y}{\func{f}{x}} ~
                         \om{y x}{x^{-1}}^{*}
                         \alphArgExt{y x}{\om{x^{-1}}{x}} ~
                         \func{\phi}{y x}
                         }{x}
                  }
             }{y} \\
= & ~ \Int{G}{
             \SqBr{\Int{G}{\Del{x}{\frac{1}{2}} \func{\psi}{y}^{*} \alphArg{y}{\func{f}{x}} ~ \om{y}{x} ~ \om{y x}{e}^{*} \func{\phi}{y x}}{x}}
             }{y} \\
= & ~ \Int{G}{\SqBr{\Int{G}{\Del{x}{\frac{1}{2}} \func{\psi}{y}^{*} \alphArg{y}{\func{f}{x}} ~ \om{y}{x} ~ \func{\phi}{y x}}{x}}}{y} \\
= & ~ \Int{G}{\SqBr{\func{\psi}{y}^{*} \Int{G}{\Del{x}{\frac{1}{2}} \alphArg{y}{\func{f}{x}} ~ \om{y}{x} ~ \func{\phi}{y x}}{x}}}{y} \\
= & ~ \Inner{\func{q}{\psi}}{\FUNC{\func{\rho}{f}}{\func{q}{\phi}}}{\LL{2}{G,A}},
\end{align*}
where the last line follows from \Identity{An Explicit Formula for the Integrated Form of the Covariant Representation of G and A}. As $ \Im{q}{\Cc{G,A}} $ is dense in $ \LL{2}{G,A} $, we conclude that
$$
\Int{G}{\GammArg{x}{\func{q}{f^{\flat}}} \centerdot \func{\phi}{x}}{x} = \FUNC{\func{\rho}{f}}{\func{q}{\phi}}.
$$
This completes the proof.
\end{proof}



\begin{Cor} \label{A Very Important Corollary}
The following statements are true:
\begin{enumerate}
\item[(i)]
$ \Ket{\func{q}{f}} = \func{\rho}{f^{\sharp}} \circ q $ and $ \KKet{\func{q}{f^{\flat}}} = \func{\rho}{f} \circ q $ for every $ f \in \Cc{G,A} $.

\item[(ii)]
$ \Im{q}{\Cc{G,A}} \subseteq \LL{2}{G,A}_{\si} $, and $ \Im{\Ket{\zeta}}{\Cc{G,A}} \subseteq \Esi $ for every $ \zeta \in \Esi $.

\item[(iii)]
$ \LL{2}{G,A} $ is a square-integrable representation of $ \Quad{G}{A}{\alpha}{\omega} $.
\end{enumerate}
\end{Cor}

\begin{proof}
By \autoref{The Workhorse Lemma}, we have for every $ f \in \Cc{G,A} $ that
$$
\forall \phi \in \Cc{G,A}: \quad
  \KetArg{\func{q}{f}}{\phi}
= \Int{G}{\GammArg{x}{\func{q}{f}} \centerdot \func{\phi}{x}}{x}
= \Int{G}{\GammArg{x}{\func{q}{f^{\sharp \flat}}} \centerdot \func{\phi}{x}}{x}
= \FUNC{\func{\rho}{f^{\sharp}}}{\func{q}{\phi}},
$$
so $ \Ket{\func{q}{f}} = \func{\rho}{f^{\sharp}} \circ q $ and consequently $ \Ket{\func{q}{f^{\flat}}} = \func{\rho}{f^{\flat \sharp}} \circ q = \func{\rho}{f} \circ q $.

As $ \func{\rho}{f^{\sharp}} \in \AdjEq{\LL{2}{G,A}} $ for every $ f \in \Cc{G,A} $, we have $ \Im{q}{\Cc{G,A}} \subseteq \LL{2}{G,A}_{\si} $ by (i) and \autoref{The Extendibility of a Ket Operator Implies the Square-Integrability of the Underlying Vector}. Then by \Identity{Ket Identity 1}, $ \Im{\Ket{\zeta}}{\Cc{G,A}} = \Im{\KKet{\zeta}}{\Im{q}{\Cc{G,A}}} \subseteq \Esi $ for every $ \zeta \in \Esi $.

Finally, as $ \Im{q}{\Cc{G,A}} $ is dense in $ \LL{2}{G,A} $, we find that $ \LL{2}{G,A}_{\si} $ is dense in $ \LL{2}{G,A} $, which makes $ \LL{2}{G,A} $ a square-integrable representation of $ \Quad{G}{A}{\alpha}{\omega} $.
\end{proof}


\subsection{A $ \Pair{\Cc{G,A}}{\star} $-Module Structure for $ \Esi $} 


\begin{ThmDef} \label{A Right-Sided Action of the Twisted Convolution Algebra on the Space of Square-Integrable Elements}
Define a bilinear map $ *_{\E}: \Esi \times \Cc{G,A} \to \Esi $ by
$$
\forall \zeta \in \Esi, ~ \forall f \in \Cc{G,A}: \quad
\zeta *_{\E} f \df \KetArg{\zeta}{f^{\flat}}.
$$
Then $ *_{\E} $ is a right $ \Pair{\Cc{G,A}}{\star} $-action on $ \Esi $.
\end{ThmDef}

\begin{proof}
For every $ \zeta \in \Esi $ and $ f,g \in \Cc{G,A} $, we have
\begin{align*}
    \KKet{\zeta *_{\E} f}
& = \KKet{\KetArg{\zeta}{f^{\flat}}} \\
& = \KKet{\KKetArg{\zeta}{\func{q}{f^{\flat}}}} \\
& = \KKet{\zeta} \circ \KKet{\func{q}{f^{\flat}}} \qquad \Br{\text{By \Identity{Ket Identity 1}.}} \\
& = \KKet{\zeta} \circ \func{\rho}{f}, \quad \text{so} \\
    \KKet{\Br{\zeta *_{\E} f} *_{\E} g}
& = \KKet{\zeta *_{\E} f} \circ \func{\rho}{g} \\
& = \KKet{\zeta} \circ \func{\rho}{f} \circ \func{\rho}{g} \\
& = \KKet{\zeta} \circ \func{\rho}{f \star g} \\
& = \KKet{\zeta *_{\E} \Br{f \star g}}, \quad
    \text{which by \autoref{Every Vector of a Twisted Hilbert C*-Module Gives a Unique Ket Operator} yields} \\
    \Br{\zeta *_{\E} f} *_{\E} g
& = \zeta *_{\E} \Br{f \star g}.
\end{align*}
Therefore, $ *_{\E} $ is indeed a right $ \Pair{\Cc{G,A}}{\star} $-action on $ \Esi $.
\end{proof}


The remaining results in this section are mostly concerned with properties of $ *_{\E} $.


\begin{Lem} \label{A Density Result for the Space of Square-Integrable Elements}
$ \Esi *_{\E} \Cc{G,A} $ is $ \Norm{\cdot}{\E} $-dense in $ \Esi $.
\end{Lem}

\begin{proof}
Recall the net $ \Seq{f_{N}}{N \in \mathcal{N}} $ in the proof of \autoref{Every Vector of a Twisted Hilbert C*-Module Gives a Unique Ket Operator}. Let $ \Seq{e_{i}}{i \in I} $ be an approximate identity for $ A $. Picking $ \zeta \in \Esi $ and $ \epsilon > 0 $, we claim that $ \BigNorm{\zeta - \KetArg{\zeta}{f_{N} e_{i}}}{\E} < \epsilon $ for some pair $ \Pair{N}{i} \in \mathcal{N} \times I $. Before proving this, first note the following assertions:
\begin{itemize}
\item
By the strong continuity of $ \gam{\E} $, there is an $ N \in \mathcal{N} $ such that $ \Norm{\zeta - \gammArg{\E}{x}{\zeta}}{\E} < \dfrac{\epsilon}{2} $ for every $ x \in N $.

\item
There is an $ i \in I $ such that $ \Norm{\zeta - \zeta \bullet e_{i}}{\E} < \dfrac{\epsilon}{2} $.
\end{itemize}
It follows from these that
\begin{align*}
       \BigNorm{\zeta - \KetArg{\zeta}{f_{N} e_{i}}}{\E}
& =    \Norm{\zeta - \Int{G}{\gammArg{\E}{x}{\zeta} \bullet \SqBr{\func{f_{N}}{x} ~ e_{i}}}{x}}{\E} \\
& \leq \Norm{\zeta - \zeta \bullet e_{i}}{\E} +
       \Norm{\zeta \bullet e_{i} - \Int{G}{\gammArg{\E}{x}{\zeta} \bullet \SqBr{\func{f_{N}}{x} ~ e_{i}}}{x}}{\E} \\
& =    \Norm{\zeta - \zeta \bullet e_{i}}{\E} +
       \Norm{
            \underbrace{\Int{G}{\zeta \bullet \SqBr{\func{f_{N}}{x} ~ e_{i}}}{x}}_{\zeta \bullet e_{i}} -
            \Int{G}{\gammArg{\E}{x}{\zeta} \bullet \SqBr{\func{f_{N}}{x} ~ e_{i}}}{x}
            }{\E} \\
& =    \Norm{\zeta - \zeta \bullet e_{i}}{\E} +
       \Norm{\Int{G}{\SqBr{\zeta - \gammArg{\E}{x}{\zeta}} \bullet \SqBr{\func{f_{N}}{x} ~ e_{i}}}{x}}{\E} \\
& \leq \Norm{\zeta - \zeta \bullet e_{i}}{\E} +
       \Int{G}{\Norm{\SqBr{\zeta - \gammArg{\E}{x}{\zeta}} \bullet \SqBr{\func{f_{N}}{x} ~ e_{i}}}{\E}}{x} \\
& \leq \Norm{\zeta - \zeta \bullet e_{i}}{\E} + \Int{G}{\Norm{\zeta - \gammArg{\E}{x}{\zeta}}{\E} \Norm{\func{f_{N}}{x} ~ e_{i}}{A}}{x} \\
& \leq \Norm{\zeta - \zeta \bullet e_{i}}{\E} + \Int{G}{\Norm{\zeta - \gammArg{\E}{x}{\zeta}}{\E} \func{f_{N}}{x}}{x} \\
& =    \Norm{\zeta - \zeta \bullet e_{i}}{\E} + \Int{N}{\Norm{\zeta - \gammArg{\E}{x}{\zeta}}{\E} \func{f_{N}}{x}}{x} \qquad
       \Br{\text{As $ \Supp{f_{N}} \subseteq N $.}} \\
& <    \frac{\epsilon}{2} + \Int{N}{\Br{\frac{\epsilon}{2}} \func{f_{N}}{x}}{x} \\
& =    \frac{\epsilon}{2} + \frac{\epsilon}{2} \Int{N}{\func{f_{N}}{x}}{x} \\
& =    \epsilon, \qquad \Br{\text{As $ \Int{N}{\func{f_{N}}{x}}{x} = 1 $.}}
\end{align*}
and the claim follows.

As $ \epsilon $ is arbitrary, we get $ \zeta \in \Cl{\Im{\Ket{\zeta}}{\Cc{G,A}}}{\E} $, and as $ \zeta $ is arbitrary, $ \Esi \subseteq \Cl{\Im{\Ket{\Esi}}{\Cc{G,A}}}{\E} $. Finally,
$$
           \Esi
\supseteq  \Esi *_{\E} \Cc{G,A}
=          \Im{\Ket{\Esi}}{\Cc{G,A}^{\flat}}
=          \Im{\Ket{\Esi}}{\Cc{G,A}},
$$
thereby concluding the proof.
\end{proof}



\begin{Lem} \label{The Main Norm Inequalities}
For every $ \zeta \in \Esi $ and $ f \in \Cc{G,A} $, the following norm inequalities hold:
\begin{align}
       \Norm{\zeta *_{\E} f}{\E}
& \leq \Norm{\zeta}{\E} \Norm{f^{\flat}}{1}, \label{Main Norm Inequality 1} \\
       \Norm{\zeta *_{\E} f}{\E,\si}
& \leq \BigNorm{\KKet{\zeta}}{\AdjEqPair{\LL{2}{G,A}}{\E}} \cdot 2 \func{\max}{\Norm{f}{1},\Norm{f}{2}}. \label{Main Norm Inequality 2}
\end{align}
\end{Lem}

\begin{proof}
For every $ \zeta \in \Esi $ and $ f \in \Cc{G,A} $, we have
\begin{align*}
       \Norm{\zeta *_{\E} f}{\E}
& =    \Norm{\KetArg{\zeta}{f^{\flat}}}{\E} \\
& \leq \Norm{\zeta}{\E} \Norm{f^{\flat}}{1}, \qquad \Br{\text{By \Inequality{Basic Norm Inequality 1}.}} \\
       \Norm{\zeta *_{\E} f}{\E,\si}
& =    \Norm{\zeta *_{\E} f}{\E} + \BigNorm{\KKet{\zeta *_{\E} f}}{\AdjEqPair{\LL{2}{G,A}}{\E}} \\
& =    \Norm{\KKetArg{\zeta}{\func{q}{f^{\flat}}}}{\E} + \BigNorm{\KKet{\zeta} \circ \func{\rho}{f}}{\AdjEqPair{\LL{2}{G,A}}{\E}} \\
& \leq \BigNorm{\KKet{\zeta}}{\AdjEqPair{\LL{2}{G,A}}{\E}} \Norm{\func{q}{f^{\flat}}}{\LL{2}{G,A}} +
       \BigNorm{\KKet{\zeta}}{\AdjEqPair{\LL{2}{G,A}}{\E}} \Norm{\func{\rho}{f}}{\AdjEq{\LL{2}{G,A}}} \\
& =    \BigNorm{\KKet{\zeta}}{\AdjEqPair{\LL{2}{G,A}}{\E}}
       \Br{\Norm{\func{q}{f^{\flat}}}{\LL{2}{G,A}} + \Norm{\func{\rho}{f}}{\AdjEq{\LL{2}{G,A}}}} \\
& =    \BigNorm{\KKet{\zeta}}{\AdjEqPair{\LL{2}{G,A}}{\E}} \Br{\NNorm{f^{\flat}} + \Norm{\func{\rho}{f}}{\AdjEq{\LL{2}{G,A}}}} \\
& \leq \BigNorm{\KKet{\zeta}}{\AdjEqPair{\LL{2}{G,A}}{\E}} \Br{\Norm{f^{\flat}}{2} + \Norm{f}{1}} \\
& =    \BigNorm{\KKet{\zeta}}{\AdjEqPair{\LL{2}{G,A}}{\E}} \Br{\Norm{f}{2} + \Norm{f}{1}} \qquad
       \Br{\text{By \autoref{The Sharp and Flat Operators Preserve L2-Norms}.}} \\
& \leq \BigNorm{\KKet{\zeta}}{\AdjEqPair{\LL{2}{G,A}}{\E}} \cdot 2 \func{\max}{\Norm{f}{1},\Norm{f}{2}}.
\end{align*}
This finishes the proof.
\end{proof}


\subsection{A Twisted-Equivariant Version of Kasparov's Stabilization Theorem} 

The tools developed earlier in this section now allow us to state and prove a twisted-equivariant version of Kasparov's Stabilization Theorem.


\begin{Prop} \label{A Twisted-Equivariant Version of Kasparov's Stabilization Theorem}
Let $ \E $ be a countably generated Hilbert $ \Quad{G}{A}{\alpha}{\omega} $-module. Then the following statements are equivalent:
\begin{enumerate}
\item[(a)]
$ \E $ is a square-integrable representation of $ \Quad{G}{A}{\alpha}{\omega} $.

\item[(b)]
There is a $ \Hilb{G}{A}{\alpha}{\omega} $-isomorphism $ \E \oplus \LL{2}{G,A}^{\infty} \cong \LL{2}{G,A}^{\infty} $.

\item[(c)]
There is a $ \Hilb{G}{A}{\alpha}{\omega} $-isomorphism from $ \E $ to a $ \Gamma^{\infty} $-invariant orthogonal summand of $ \LL{2}{G,A}^{\infty} $.
\end{enumerate}
\end{Prop}

\begin{proof}
We will follow the structure of the argument in \cite{Meyer1}, which is a variant of that given in \cite{Mingo|Phillips}. \\

\noindent \ul{(a) implies (b)} \\

Suppose that $ \E $ is a square-integrable representation of $ \Quad{G}{A}{\alpha}{\omega} $. There is a sequence $ \Seq{\zeta_{n}}{n \in \N{}{}} $ in $ \Esi $ such that $ \Cl{\Span{\SSet{\zeta_{n}}_{n \in \N{}{}}}}{\E} = \E $. Re-scaling, we may assume that $ \Norm{\BBra{\zeta_{n}}}{\AdjEqPair{\E}{\LL{2}{G,A}}} \leq 1 $ for every $ n \in \N{}{} $, and we may arrange for each member of the sequence to be repeated infinitely often.

Define an operator $ T: \LL{2}{G,A}^{\infty} \to \E \oplus \LL{2}{G,A}^{\infty} $ by
$$
\forall \bm{\Phi} \in \LL{2}{G,A}^{\infty}: \quad
    \func{T}{\sum_{n = 1}^{\infty} \Phi_{n} \cdot \mathbf{e}_{n}}
\df \SqBr{\sum_{n = 1}^{\infty} \frac{1}{2^{n}} \KKetArg{\zeta_{n}}{\Phi_{n}}} \oplus
    \SqBr{\sum_{n = 1}^{\infty} \frac{1}{4^{n}} \Phi_{n} \cdot \mathbf{e}_{n}}.
$$
This is an adjointable operator, whose adjoint $ T^{*}: \E \oplus \LL{2}{G,A}^{\infty} \to \LL{2}{G,A}^{\infty} $ is given by
$$
\forall \eta \in \E, ~ \forall \bm{\Phi} \in \LL{2}{G,A}^{\infty}: \quad
    \func{T^{*}}{\eta \oplus \sum_{n = 1}^{\infty} \Phi_{n} \cdot \mathbf{e}_{n}}
\df \sum_{n = 1}^{\infty} \SqBr{\frac{1}{2^{n}} \BBraArg{\zeta_{n}}{\eta} + \frac{1}{4^{n}} \Phi_{n}} \cdot \mathbf{e}_{n}.
$$
By \autoref{Meyer's BBra-KKet Operators Are Morphisms of Twisted Hilbert C*-Modules}, $ T $ and $ T^{*} $ are twisted-equivariant. Also, $ T^{*} $ has dense range as the set of elements of $ \LL{2}{G,A}^{\infty} $ with finitely many non-zero components is dense in $ \LL{2}{G,A}^{\infty} $ and any such element is equal to $ \func{T^{*}}{0_{\E} \oplus \bm{\Phi}} $ for some $ \bm{\Phi} \in \LL{2}{G,A}^{\infty} $ with finitely many non-zero components too.

We claim that $ T $ has dense range as well. Let $ R \df \Cl{\Range{T}}{\E \oplus \LL{2}{G,A}^{\infty}} $. Pick any $ \zeta \in \SSet{\zeta_{n}}_{n \in \N{}{}} $, and let $ N \df \Set{n \in \N{}{}}{\zeta_{n} = \zeta} $, which is an infinite set. Then
$$
\forall \phi \in \Cc{G,A}, ~ \forall n \in N: \quad
\func{T}{2^{n} \phi \cdot \mathbf{e}_{n}}    =   \KKetArg{\zeta}{\phi} \oplus \frac{1}{2^{n}} \phi \cdot \mathbf{e}_{n},
\quad \text{whence} \quad
\KKetArg{\zeta}{\phi} \oplus \bm{0}_{\infty} \in R.
$$
The proof of \autoref{A Density Result for the Space of Square-Integrable Elements} says that $ \zeta \in \Cl{\Im{\KKet{\zeta}}{\Cc{G,A}}}{\E} $, so $ \zeta \oplus \bm{0}_{\infty} \in R $. Then as $ \zeta $ is arbitrary, $ \zeta_{n} \oplus \bm{0}_{\infty} \in R $ for every $ n \in \N{}{} $, giving $ \E \oplus \bm{0}_{\infty} \subseteq R $. Hence, $ 0_{\E} \oplus \Phi \cdot \mathbf{e}_{n} \in R $ for every $ \Phi \in \LL{2}{G,A} $ and $ n \in \N{}{} $ because
$$
    \func{T}{4^{n} \Phi \cdot \mathbf{e}_{n}}
=   2^{n} \KKetArg{\zeta_{n}}{\Phi} \oplus \Phi \cdot \mathbf{e}_{n}
\in R \qquad \text{and} \qquad
    2^{n} \KKetArg{\zeta_{n}}{\Phi} \oplus \bm{0}_{\infty}
\in R.
$$
Therefore, $ 0_{\E} \oplus \LL{2}{G,A}^{\infty} \subseteq R $, which leads to $ R = \E \oplus \LL{2}{G,A}^{\infty} $.

As $ T $ and $ T^{*} $ have dense range, so does $ T^{*} \circ T $. The same then goes for $ \Abs{T} \df \Br{T^{*} T}^{\frac{1}{2}} $ as
$$
          \Im{\Abs{T}}{\LL{2}{G,A}^{\infty}}
\supseteq \Im{\Abs{T}}{\Im{\Abs{T}}{\LL{2}{G,A}^{\infty}}}
=         \IM{T^{*} \circ T}{\LL{2}{G,A}^{\infty}}.
$$
Observe that
\begin{align*}
\forall \bm{\Phi} \in \LL{2}{G,A}^{\infty}: \quad
    \Inner{\func{\Abs{T}}{\bm{\Phi}}}{\func{\Abs{T}}{\bm{\Phi}}}{\LL{2}{G,A}^{\infty}}
& = \Inner{\Func{T^{*} \circ T}{\bm{\Phi}}}{\bm{\Phi}}{\LL{2}{G,A}^{\infty}} \\
& = \Inner{\func{T}{\bm{\Phi}}}{\func{T}{\bm{\Phi}}}{\E \oplus \LL{2}{G,A}^{\infty}},
\end{align*}
so there is an $ A $-linear isometry $ U: \Range{\Abs{T}} \to \Range{T} $ defined by
$$
\forall \bm{\Phi} \in \LL{2}{G,A}^{\infty}: \quad
\func{U}{\func{\Abs{T}}{\bm{\Phi}}} \df \func{T}{\bm{\Phi}}.
$$
As polynomials in $ T^{*} \circ T $ belong to $ \AdjEq{\LL{2}{G,A}^{\infty}} $, we obtain $ \Abs{T} \in \AdjEq{\LL{2}{G,A}^{\infty}} $. Hence, $ \Range{\Abs{T}} $ is $ \Gamma^{\infty} $-invariant and $ U \circ \Br{\Gamma^{\infty}}_{r} = \Br{\gam{\E} \oplus \Gamma^{\infty}}_{r} \circ U $ for every $ r \in G $. Now, extend $ U $ to a surjective $ A $-linear isometry $ V: \LL{2}{G,A}^{\infty} \to \E \oplus \LL{2}{G,A}^{\infty} $. Then $ V $ is unitary, and as it is twisted-equivariant, we are done. \\

\noindent \ul{(b) implies (c)} \\

This is tautological. \\

\noindent \ul{(c) implies (a)} \\

By \autoref{A Very Important Corollary}, $ \LL{2}{G,A} $ is a square-integrable representation of $ \Quad{G}{A}{\alpha}{\omega} $. By \autoref{A Direct Sum of Square-Integrable Representations Is Square-Integrable}, $ \LL{2}{G,A}^{\infty} $ is then one as well. By abuse of notation, suppose that $ \E $ itself is an orthogonal summand of $ \LL{2}{G,A}^{\infty} $, and let $ P \in \AdjEqPair{\LL{2}{G,A}^{\infty}}{\E} $ denote the associated projection map. By \Identity{Ket Identity 1}, $ \Im{P}{\SqBr{\LL{2}{G,A}^{\infty}}_{\si}} \subseteq \Esi $. As $ P $ is surjective and $ \SqBr{\LL{2}{G,A}^{\infty}}_{\si} $ is dense in $ \LL{2}{G,A}^{\infty} $, we see that $ \Esi $ is dense in $ \E $. Therefore, $ \E $ is a square-integrable representation of $ \Quad{G}{A}{\alpha}{\omega} $.
\end{proof}



\begin{Rmk} \label{A Missing Connection with Integrable Group Actions}
We have not mentioned integrable group actions above, which appear in the statement of the equivariant version of Kasparov's Stabilization Theorem in \cite{Meyer1}. It is not known to us how the content of Rieffel's paper \cite{Rieffel2} may be adapted to the twisted case, so we will not pursue this here.
\end{Rmk}



\section{Approximate Identities} 

This section contains technical results about approximate identities.

Recall the net $ \Seq{f_{N}}{N \in \mathcal{N}} $ in the proof of \autoref{Every Vector of a Twisted Hilbert C*-Module Gives a Unique Ket Operator}, and let $ \Seq{e_{i}}{i \in I} $ be an approximate identity for $ A $. Then $ \Seq{f_{N} e_{i}}{\Pair{N}{i} \in \mathcal{N} \times I} $ is a left approximate identity for $ \Quad{\Cc{G,A}}{\star}{^{*}}{\Norm{\cdot}{1}} $ norm-bounded by $ 1 $. As $ ^{*} $ is isometric for $ \Norm{\cdot}{1} $, we can use a well-known algebraic trick (see \autoref{Constructing a Two-Sided Approximate Identity from Left and Right Ones} below) for converting this into a bounded self-adjoint (with respect to $ ^{*} $) two-sided approximate identity.

However, this is not sufficient, for reasons to be explained in the proof of \autoref{A Pseudo-Approximate Identity for Relatively Continuous Subspaces}. If
$$
\Norm{\cdot}{\flat} \df \Map{\Cc{G,A}}{\RR{\geq 0}{}}{f}{\Norm{f^{\flat}}{1}},
$$
then, as will be shown, $ \Norm{\cdot}{\flat} $ is an algebra norm (i.e., sub-multiplicativity holds) on $ \Trip{\Cc{G,A}}{\star}{^{*}} $ for which $ ^{*} $ is not necessarily isometric. Our goal is to have a \emph{single} bounded self-adjoint two-sided approximate identity for both $ \Quad{\Cc{G,A}}{\star}{^{*}}{\Norm{\cdot}{1}} $ and $ \Quad{\Cc{G,A}}{\star}{^{*}}{\Norm{\cdot}{\flat}} $.

We will construct the desired approximate identity from scratch, but first, let us prove that $ \Norm{\cdot}{\flat} $ is an algebra norm.


\begin{Lem} \label{Some Results on Convolution and Involution}
For every $ f,g \in \Cc{G,A} $, the following hold:
\begin{align}
\Norm{f}{\flat}                       & = \Norm{\Delta^{- \frac{1}{2}} f}{1}, \label{Flat-Norm Identity} \\
\Br{\Delta^{- \frac{1}{2}} f}^{*}     & = \Delta^{\frac{1}{2}} f^{*}, \label{Involution Identity 1} \\
\Br{\Delta^{\frac{1}{2}} f}^{*}       & = \Delta^{- \frac{1}{2}} f^{*}, \label{Involution Identity 2} \\
\Delta^{- \frac{1}{2}} \Br{f \star g} & = \Br{\Delta^{- \frac{1}{2}} f} \star \Br{\Delta^{- \frac{1}{2}} g}. \label{Convolution Identity}
\end{align}
Therefore, $ \Norm{\cdot}{\flat} $ is an algebra norm on $ \Trip{\Cc{G,A}}{\star}{^{*}} $.
\end{Lem}

\begin{proof}
Let $ f \in \Cc{G,A} $. We have seen in the proof of \autoref{The Sharp and Flat Operators Preserve L2-Norms} that $ \Norm{\Flat{f}{x}}{A} = \Del{x}{- \frac{1}{2}} \Norm{\func{f}{x^{-1}}}{A} $ for every $ x \in G $, so
\begin{align*}
    \Norm{f}{\flat}
& = \Norm{f^{\flat}}{1} \\
& = \Int{G}{\Del{x}{- \frac{1}{2}} \Norm{\func{f}{x^{-1}}}{A}}{x} \\
& = \Int{G}{\Del{x^{-1}}{} \Del{x}{\frac{1}{2}} \Norm{\func{f}{x}}{A}}{x} \\
& = \Int{G}{\Del{x}{-1} \Del{x}{\frac{1}{2}} \Norm{\func{f}{x}}{A}}{x} \\
& = \Int{G}{\Del{x}{- \frac{1}{2}} \Norm{\func{f}{x}}{A}}{x} \\
& = \Int{G}{\Norm{\Del{x}{- \frac{1}{2}} \func{f}{x}}{A}}{x} \\
& = \Norm{\Delta^{- \frac{1}{2}} f}{1}.
\end{align*}
This establishes \Identity{Flat-Norm Identity}.

Next, for every $ x \in G $, we have
\begin{align*}
    \Br{\Delta^{- \frac{1}{2}} f}^{*}
& = \Map{G}{A}{x}{\Del{x}{-1} \om{x}{x^{-1}}^{*} \alphArg{x}{\Del{x}{\frac{1}{2}} \func{f}{x^{-1}}}^{*}} \\
& = \Map{G}{A}{x}{\Del{x}{\frac{1}{2}} \SqBr{\Del{x}{-1} \om{x}{x^{-1}}^{*} \alphArg{x}{\func{f}{x^{-1}}}^{*}}} \\
& = \Delta^{\frac{1}{2}} f^{*}, \\
    \Br{\Delta^{\frac{1}{2}} f}^{*}
& = \Map{G}{A}{x}{\Del{x}{-1} \om{x}{x^{-1}}^{*} \alphArg{x}{\Del{x}{- \frac{1}{2}} \func{f}{x^{-1}}}^{*}} \\
& = \Map{G}{A}{x}{\Del{x}{- \frac{1}{2}} \SqBr{\Del{x}{-1} \om{x}{x^{-1}}^{*} \alphArg{x}{\func{f}{x^{-1}}}^{*}}} \\
& = \Delta^{- \frac{1}{2}} f^{*}.
\end{align*}
This yields \Identity{Involution Identity 1} and \Identity{Involution Identity 2}.

Let $ g \in \Cc{G,A} $ also. Then for every $ x \in G $, we have
\begin{align*}
    \Delta^{- \frac{1}{2}} \Br{f \star g}
& = \Map{G}{A}{x}{\Del{x}{- \frac{1}{2}} \Int{G}{\func{f}{y} ~ \alphArg{y}{\func{g}{y^{-1} x}} ~ \om{y}{y^{-1} x}}{y}} \\
& = \Map{G}{A}{x}{
                 \Int{G}{
                        \SqBr{\Del{y}{- \frac{1}{2}} \func{f}{y}}
                        \alphArg{y}{\Del{y^{-1} x}{- \frac{1}{2}} \func{g}{y^{-1} x}} ~
                        \om{y}{y^{-1} x}
                        }{y}
                 } \\
& = \Br{\Delta^{- \frac{1}{2}} f} \star \Br{\Delta^{- \frac{1}{2}} g}.
\end{align*}
This proves \Identity{Convolution Identity}.

Finally, observe that
\begin{align*}
       \Norm{f \star g}{\flat}
& =    \Norm{\Delta^{- \frac{1}{2}} \Br{f \star g}}{1} \qquad \Br{\text{By \Identity{Flat-Norm Identity}.}} \\
& =    \Norm{\Br{\Delta^{- \frac{1}{2}} f} \star \Br{\Delta^{- \frac{1}{2}} g}}{1} \qquad \Br{\text{By \Identity{Convolution Identity}.}} \\
& \leq \Norm{\Delta^{- \frac{1}{2}} f}{1} \Norm{\Delta^{- \frac{1}{2}} g}{1} \qquad
       \Br{\text{As $ \Norm{\cdot}{1} $ is sub-multiplicative.}} \\
& =    \Norm{f}{\flat} \Norm{g}{\flat}, \qquad \Br{\text{By \Identity{Flat-Norm Identity} again.}}
\end{align*}
so $ \Norm{\cdot}{\flat} $ is sub-multiplicative and thus an algebra norm on $ \Trip{\Cc{G,A}}{\star}{^{*}} $.
\end{proof}


Let $ \Seq{g_{N}}{N \in \mathcal{N}} \df \Seq{f_{N} \cdot \func{\min}{\Delta^{\frac{1}{2}},\Delta^{- \frac{1}{2}}}}{N \in \mathcal{N}} $. Then for every $ N \in \mathcal{N} $, the following hold:
$$
\bm{0} \leq g_{N}                        \leq f_{N} \cdot \bm{1}                                     = f_{N}, \quad
\bm{0} \leq g_{N} \Delta^{\frac{1}{2}}   \leq \Br{f_{N} \Delta^{- \frac{1}{2}}} \Delta^{\frac{1}{2}} = f_{N}, \quad
\bm{0} \leq g_{N} \Delta^{- \frac{1}{2}} \leq \Br{f_{N} \Delta^{\frac{1}{2}}} \Delta^{- \frac{1}{2}} = f_{N}.
$$
Furthermore, for every $ f \in \Cont{}{G,A} $,
\begin{alignat*}{3}
     \lim_{N \in \mathcal{N}} \Int{G}{\func{g_{N}}{x} ~ \func{f}{x}}{x}
&  = \lim_{N \in \mathcal{N}} \Int{G}{\func{f_{N}}{x} \cdot \func{\min}{\Del{x}{\frac{1}{2}},\Del{x}{- \frac{1}{2}}} ~ \func{f}{x}}{x}
&& = \func{f}{e}, \\
     \lim_{N \in \mathcal{N}} \Int{G}{\func{g_{N}}{x} ~ \Del{x}{\frac{1}{2}} \func{f}{x}}{x}
&  = \lim_{N \in \mathcal{N}}
     \Int{G}{\func{f_{N}}{x} \cdot \func{\min}{\Del{x}{\frac{1}{2}},\Del{x}{- \frac{1}{2}}} ~ \Del{x}{\frac{1}{2}} \func{f}{x}}{x}
&& = \func{f}{e}, \\
     \lim_{N \in \mathcal{N}} \Int{G}{\func{g_{N}}{x} ~ \Del{x}{- \frac{1}{2}} \func{f}{x}}{x}
&  = \lim_{N \in \mathcal{N}}
     \Int{G}{\func{f_{N}}{x} \cdot \func{\min}{\Del{x}{\frac{1}{2}},\Del{x}{- \frac{1}{2}}} ~ \Del{x}{- \frac{1}{2}} \func{f}{x}}{x}
&& = \func{f}{e}.
\end{alignat*}
We have implicitly used the facts that $ \Seq{f_{N}}{N \in \mathcal{N}} $ is an approximating delta at $ e $ and that $ \func{\Delta}{e} = 1 $.
Hence, the nets $ \Seq{g_{N}}{N \in \mathcal{N}} $, $ \Seq{g_{N} \Delta^{\frac{1}{2}}}{N \in \mathcal{N}} $ and $ \Seq{g_{N} \Delta^{- \frac{1}{2}}}{N \in \mathcal{N}} $ are approximating deltas at $ e $ that are $ \Norm{\cdot}{1} $-bounded by $ 1 $.


\begin{Thm} \label{Three Bounded Left Approximate Identities for the Twisted Convolution Algebra}
Let $ \Seq{e_{i}}{i \in I} $ be an approximate identity in $ A $. Then the nets
$$
\Seq{g_{N} e_{i}}{\Pair{N}{i} \in \mathcal{N} \times I}, \qquad
\Seq{g_{N} \Delta^{\frac{1}{2}} e_{i}}{\Pair{N}{i} \in \mathcal{N} \times I}, \qquad
\Seq{g_{N} \Delta^{- \frac{1}{2}} e_{i}}{\Pair{N}{i} \in \mathcal{N} \times I}
$$
in $ \Cc{G,A} $ are left approximate identities for $ \Quad{\Cc{G,A}}{\star}{^{*}}{\Norm{\cdot}{1}} $ that are norm-bounded by $ 1 $.
\end{Thm}

\begin{proof}
Let $ \Seq{h_{N,i}}{\Pair{N}{i} \in \mathcal{N} \times I} $ denote any of these nets. Then $ \Norm{h_{N,i}}{1} \leq 1 $ for every $ \Pair{N}{i} \in \mathcal{N} \times I $ and
$$
\forall a \in A: \quad
\lim_{\Pair{N}{i} \in \mathcal{N} \times I} \Int{G}{\func{h_{N,i}}{x} ~ a}{x} = a.
$$
Letting $ f \in \Cc{G} \setminus \SSet{\bm{0}} $ and $ a \in A $, we first prove that
$$
\lim_{\Pair{N}{i} \in \mathcal{N} \times I} \Norm{h_{N,i} \star f a - f a}{1} = 0.
$$

To begin, observe for every $ \Pair{N}{i} \in \mathcal{N} \times I $ that
\begin{align*}
     & ~ \Norm{h_{N,i} \star f a - f a}{1} \\
=    & ~ \Int{G}{\Norm{\Func{h_{N,i} \star f a}{x} - \Func{f a}{x}}{A}}{x} \\
=    & ~ \Int{G}{\Norm{\Int{G}{\func{h_{N,i}}{y} ~ \alphArg{y}{\func{f}{y^{-1} x} ~ a} ~ \om{y}{y^{-1} x}}{y} - \func{f}{x} ~ a}{A}}{x} \\
\leq & ~ \Int{G}{
                \Norm{
                     \Int{G}{\func{h_{N,i}}{y} ~ \alphArg{y}{\func{f}{y^{-1} x} ~ a} ~ \om{y}{y^{-1} x}}{y} -
                     \Int{G}{\func{h_{N,i}}{y} ~ \func{f}{x} ~ a}{y}
                     }{A}
                }{x} ~ + \\
     & ~ \Int{G}{\Norm{\Int{G}{\func{h_{N,i}}{y} ~ \func{f}{x} ~ a}{y} - \func{f}{x} ~ a}{A}}{x} \\
=    & ~ \Int{G}{
                \Norm{\Int{G}{\func{h_{N,i}}{y} \SqBr{\alphArg{y}{\func{f}{y^{-1} x} ~ a} ~ \om{y}{y^{-1} x} - \func{f}{x} ~ a}}{y}}{A}
                }{x} ~ + \\
     & ~ \Int{G}{\Norm{\func{f}{x} \SqBr{\Int{G}{\func{h_{N,i}}{y} ~ a}{y} - a}}{A}}{x} \\
\leq & ~ \Int{G}{
                \SqBr{
                     \Int{G}{\Norm{\func{h_{N,i}}{y}}{A} \Norm{\alphArg{y}{\func{f}{y^{-1} x} ~ a} ~ \om{y}{y^{-1} x} - \func{f}{x} ~ a}{A}}{y}
                     }
                }{x} ~ + \\
     & ~ \Int{G}{\Abs{\func{f}{x}} \Norm{\Int{G}{\func{h_{N,i}}{y} ~ a}{y} - a}{A}}{x} \\
=    & ~ \Int{G}{
                \SqBr{
                     \Int{G}{\Norm{\func{h_{N,i}}{y}}{A} \Norm{\alphArg{y}{\func{f}{y^{-1} x} ~ a} ~ \om{y}{y^{-1} x} - \func{f}{x} ~ a}{A}}{x}
                     }
                }{y} ~ + \\
     & ~ \Int{G}{\Abs{\func{f}{x}} \Norm{\Int{G}{\func{h_{N,i}}{y} ~ a}{y} - a}{A}}{x} \\
     & ~ \Br{\text{By Fubini's Theorem.}} \\
=    & ~ \Int{G}{
                \SqBr{
                     \Int{G}{\Norm{\func{h_{N,i}}{y}}{A} \Norm{\alphArg{y}{\func{f}{y^{-1} x} ~ a} ~ \om{y}{y^{-1} x} - \func{f}{x} ~ a}{A}}{x}
                     }
                }{y} ~ + \\
     & ~ \Norm{f}{1} \Norm{\Int{G}{\func{h_{N,i}}{y} ~ a}{y} - a}{A}.
\end{align*}
Let $ \epsilon > 0 $ and $ S \df \Supp{f} $. Fix a compact subset $ K $ of $ G $ with $ e $ in its interior. By continuity, find $ K S $-indexed sequences $ \Seq{V_{x}}{x \in K S} $ and $ \Seq{W_{x}}{x \in K S} $ of subsets of $ G $ so that for every $ x \in K S $:
\begin{itemize}
\item
$ V_{x} $ is the intersection of $ K S $ with an open neighborhood of $ x $.

\item
$ W_{x} $ is the intersection of $ K^{\circ} $ with an open neighborhood of $ e $.

\item
$ \Norm{\alphArg{y}{\func{f}{y^{-1} z} ~ a} ~ \om{y}{y^{-1} z} - \func{f}{x} ~ a}{A} < \dfrac{\epsilon}{4 \func{\mu}{K S}} $ for every $ \Pair{z}{y} \in V_{x} \times W_{x} $, whence
\begin{equation}
\forall \Pair{z}{y} \in V_{x} \times W_{x}: \quad
  \Norm{\alphArg{y}{\func{f}{y^{-1} z} ~ a} ~ \om{y}{y^{-1} z} - \func{f}{z} ~ a}{A}
< \frac{\epsilon}{2 \func{\mu}{K S}}. \label{A Most Crucial Approximation}
\end{equation}
\end{itemize}
By the compactness of $ K S $, there exist $ x_{1},\ldots,x_{n} \in K S $ such that $ \D K S = \bigcup_{k = 1}^{n} V_{x_{k}} $. Pick any $ N \in \mathcal{N} $ contained in $ \D \bigcap_{k = 1}^{n} W_{x_{k}} $, and let $ \Pair{x}{y} \in N S \times N $. As $ \D N S \subseteq K S $, there is a $ k \in \SSet{1,\ldots,n} $ such that $ x \in V_{x_{k}} $, and as $ e,y \in W_{x_{k}} $, \Inequality{A Most Crucial Approximation} gives us
$$
\Norm{\alphArg{y}{\func{f}{y^{-1} x} ~ a} ~ \om{y}{y^{-1} x} - \func{f}{x} ~ a}{A} < \frac{\epsilon}{2 \func{\mu}{K S}}.
$$
We chose $ \Pair{x}{y} $ arbitrarily, so
\begin{align*}
     & ~ \Int{G}{
                \Int{G}{\Norm{\func{h_{N,i}}{y}}{A} \Norm{\alphArg{y}{\func{f}{y^{-1} x} ~ a} ~ \om{y}{y^{-1} x} - \func{f}{x} ~ a}{A}}{x}
                }{y} \\
=    & ~ \Int{N}{
                \Int{N S}{\Norm{\func{h_{N,i}}{y}}{A} \Norm{\alphArg{y}{\func{f}{y^{-1} x} ~ a} ~ \om{y}{y^{-1} x} - \func{f}{x} ~ a}{A}}{x}
                }{y} \\
     & ~ \Br{\text{As the integrand vanishes outside of $ N S \times N $.}} \\
\leq & ~ \Int{N}{\Int{N S}{\Norm{\func{h_{N,i}}{y}}{A} \SqBr{\frac{\epsilon}{2 \func{\mu}{K S}}}}{x}}{y} \\
=    & ~ \Int{N}{\Norm{\func{h_{N,i}}{y}}{A} \SqBr{\frac{\epsilon}{2 \func{\mu}{K S}}} \func{\mu}{N S}}{y} \\
=    & ~ \SqBr{\frac{\epsilon}{2 \func{\mu}{K S}}} \func{\mu}{N S} \Int{N}{\Norm{\func{h_{N,i}}{y}}{A}}{y} \\
\leq & ~ \SqBr{\frac{\epsilon}{2 \func{\mu}{K S}}} \func{\mu}{N S} \\
\leq & ~ \SqBr{\frac{\epsilon}{2 \func{\mu}{K S}}} \func{\mu}{K S} \\
=    & ~ \frac{\epsilon}{2}.
\end{align*}
Next, pick $ U \in \mathcal{N} $ and $ i_{0} \in I $ so that for every $ N \in \mathcal{N} $ contained in $ U $ and every $ i \in I_{\geq i_{0}} $,
$$
\Norm{\Int{G}{\func{h_{N,i}}{y} ~ a}{y} - a}{A} < \dfrac{\epsilon}{2 \Br{\Norm{f}{1} + 1}}.
$$
Hence, for every $ N \in \mathcal{N} $ contained in $ \D U \cap \bigcap_{k = 1}^{n} W_{x_{k}} $ and every $ i \in I_{\geq i_{0}} $, we have
$$
  \Norm{h_{N,i} \star f a - f a}{1}
< \frac{\epsilon}{2} + \frac{\epsilon}{2}
= \epsilon.
$$
As $ \epsilon $ is arbitrary, we obtain $ \D \lim_{\Pair{N}{i} \in \mathcal{N} \times I} \Norm{h_{N,i} \star f a - f a}{1} = 0 $.

By Urysohn's Lemma, $ \Cc{G} \odot A $ is $ \Norm{\cdot}{1} $-dense in $ \Cc{G,A} $. Let $ f \in \Cc{G,A} $ and $ \epsilon > 0 $. Find $ f_{1},\ldots,f_{n} \in \Cc{G} $ and $ a_{1},\ldots,a_{n} \in A $ so that $ \D \Norm{f - \sum_{k = 1}^{n} f_{k} a_{k}}{1} < \frac{\epsilon}{3} $. Then for every $ \Pair{N}{i} \in \mathcal{N} \times I $,
\begin{align*}
     & ~ \Norm{h_{N,i} \star f - f}{1} \\
\leq & ~ \Norm{h_{N,i} \star f - \sum_{k = 1}^{n} h_{N,i} \star f_{k}  a_{k}}{1} +
         \Norm{\sum_{k = 1}^{n} h_{N,i} \star f_{k} a_{k} - \sum_{k = 1}^{n} f_{k} a_{k}}{1} +
         \Norm{\sum_{k = 1}^{n} f_{k} a_{k} - f}{1} \\
=    & ~ \Norm{h_{N,i} \star \Br{f - \sum_{k = 1}^{n} f_{k} a_{k}}}{1} +
         \Norm{\sum_{k = 1}^{n} \Br{h_{N,i} \star f_{k} a_{k} - f_{k} a_{k}}}{1} +
         \Norm{\sum_{k = 1}^{n} f_{k} a_{k} - f}{1} \\
\leq & ~ \Norm{h_{N,i}}{1} \Norm{f - \sum_{k = 1}^{n} f_{k} a_{k}}{1} +
         \sum_{k = 1}^{n} \Norm{h_{N,i} \star f_{k} a_{k} - f_{k} a_{k}}{1} +
         \Norm{\sum_{k = 1}^{n} f_{k} a_{k} - f}{1} \\
<    & ~ 1 \cdot \frac{\epsilon}{3} +
         \sum_{k = 1}^{n} \Norm{h_{N,i} \star f_{k} a_{k} - f_{k} a_{k}}{1} +
         \frac{\epsilon}{3} \\
=    & ~ \frac{2 \epsilon}{3} + \sum_{k = 1}^{n} \Norm{h_{N,i} \star f_{k} a_{k} - f_{k} a_{k}}{1}.
\end{align*}
By the foregoing argument, we can pick $ \Pair{N_{0}}{i_{0}} \in \mathcal{N} \times I $ so that the second term in the last line is $ < \dfrac{\epsilon}{3} $ for every $ \Pair{N}{i} \in \Br{\mathcal{N} \times I}_{\geq \Pair{N_{0}}{i_{0}}} $. As $ \epsilon $ is arbitrary, we obtain $ \D \lim_{\Pair{N}{i} \in \mathcal{N} \times I} \Norm{h_{N,i} \star f - f}{1} = 0 $ for every $ f \in \Cc{G,A} $.
\end{proof}


The net $ \Seq{g_{N} e_{i}}{\Pair{N}{i} \in \mathcal{N} \times I} $ plays an important role, so we will for brevity denote it by $ \Seq{u_{j}}{j \in \mathcal{N} \times I} $.

The algebraic trick we mentioned at the start of this section appears as Proposition 2.6 of \cite{Doran|Wichmann}, which we shall state and prove next.


\begin{Thm} \label{Constructing a Two-Sided Approximate Identity from Left and Right Ones}
Let $ X $ be a normed algebra. If $ \Seq{e_{i}}{i \in I} $ and $ \Seq{f_{j}}{j \in J} $ are, respectively, left and right approximate identities for $ X $, and both are norm-bounded by $ M > 0 $, then $ \Seq{e_{i} + f_{j} - f_{j} e_{i}}{\Pair{i}{j} \in I \times J} $ is a two-sided approximate identity for $ X $ that is norm-bounded by $ M \Br{M + 2} $.
\end{Thm}

\begin{proof}
As $ \D \sup_{i \in I} ~ \Norm{e_{i}}{X} \leq M $ and $ \D \sup_{j \in J} ~ \Norm{f_{j}}{X} \leq M $, the Triangle Inequality yields
$$
     \sup_{\Pair{i}{j} \in I \times J} \Norm{e_{i} + f_{j} - f_{j} e_{i}}{X}
\leq \sup_{\Pair{i}{j} \in I \times J} \Br{\Norm{e_{i}}{X} + \Norm{f_{j}}{X} + \Norm{f_{j}}{X} \Norm{e_{i}}{X}}
\leq M + M + M^{2}
=    M \Br{M + 2}.
$$
This proves the boundedness assertion.

Let $ x \in X $ and $ \epsilon > 0 $. Observe that
\begin{align*}
\forall \Pair{i}{j} \in I \times J: \quad
       \Norm{x \Br{e_{i} + f_{j} - f_{j} e_{i}} - x}{X}
& =    \Norm{x e_{i} + x f_{j} - x f_{j} e_{i} - x}{X} \\
& =    \Norm{x e_{i} - x f_{j} e_{i} + x f_{j} - x}{X} \\
& =    \Norm{\Br{x - x f_{j}} e_{i} + x f_{j} - x}{X} \\
& \leq \Norm{\Br{x - x f_{j}} e_{i}}{X} + \Norm{x f_{j} - x}{X} \\
& \leq \Norm{x - x f_{j}}{X} \Norm{e_{i}}{X} + \Norm{x f_{j} - x}{X} \\
& \leq M \Norm{x - x f_{j}}{X} + \Norm{x f_{j} - x}{X} \\
& =    \Br{M + 1} \Norm{x f_{j} - x}{X}, \\
       \Norm{\Br{e_{i} + f_{j} - f_{j} e_{i}} x - x}{X}
& =    \Norm{e_{i} x + f_{j} x - f_{j} e_{i} x - x}{X} \\
& =    \Norm{e_{i} x - x + f_{j} x - f_{j} e_{i} x}{X} \\
& =    \Norm{e_{i} x - x + f_{j} \Br{x - e_{i} x}}{X} \\
& \leq \Norm{e_{i} x - x}{X} + \Norm{f_{j} \Br{x - e_{i} x}}{X} \\
& \leq \Norm{e_{i} x - x}{X} + \Norm{f_{j}}{X} \Norm{x - e_{i}}{X} \\
& \leq \Norm{e_{i} x - x}{X} + M \Norm{x - e_{i}}{X} \\
& =    \Br{M + 1} \Norm{e_{i} x - x}{X}.
\end{align*}
Pick $ \Pair{i_{0}}{j_{0}} \in I \times J $ so that for every $ \Pair{i}{j} \in \Br{I \times J}_{\geq \Pair{i_{0}}{j_{0}}} $,
\begin{align*}
    \Norm{e_{i} x - x}{X},\Norm{x f_{j} - x}{X}
& < \dfrac{\epsilon}{M + 1}, \quad \text{in which case} \quad \\
    \Norm{x \Br{e_{i} + f_{j} - f_{j} e_{i}} - x}{X},\Norm{\Br{e_{i} + f_{j} - f_{j} e_{i}} x - x}{X}
& < \epsilon.
\end{align*}
As $ \epsilon $ is arbitrary, we therefore obtain
$$
  \lim_{\Pair{i}{j} \in I \times J} \Norm{x \Br{e_{i} + f_{j} - f_{j} e_{i}} - x}{X}
= \lim_{\Pair{i}{j} \in I \times J} \Norm{\Br{e_{i} + f_{j} - f_{j} e_{i}} x - x}{X}
= 0,
$$
which concludes the proof.
\end{proof}



\begin{Prop} \label{A Self-Adjoint Two-Sided Approximate Identity for the Twisted Convolution Algebra for Two Different Norms}
The net $ \Seq{v_{j}}{j \in \mathcal{N} \times I} \df \Seq{u_{j} + u_{j}^{*} - u_{j}^{*} \star u_{j}}{j \in \mathcal{N} \times I} $ is a single self-adjoint two-sided approximate identity for both $ \Quad{\Cc{G,A}}{\star}{^{*}}{\Norm{\cdot}{1}} $ and $ \Quad{\Cc{G,A}}{\star}{^{*}}{\Norm{\cdot}{\flat}} $ norm-bounded by $ 3 $.
\end{Prop}

\begin{proof}
Self-adjointness is evident for both normed $ * $-algebras, being a purely algebraic condition.

By \autoref{Three Bounded Left Approximate Identities for the Twisted Convolution Algebra}, $ \Seq{u_{j}}{j \in \mathcal{N} \times I} $ is a left approximate identity for $ \Quad{\Cc{G,A}}{\star}{^{*}}{\Norm{\cdot}{1}} $. Consequently,
$$
\forall f \in \Cc{G,A}: \quad
  \lim_{j \in \mathcal{N} \times I} \Norm{f \star u_{j}^{*} - f}{1}
= \lim_{j \in \mathcal{N} \times I} \Norm{\Br{f \star u_{j}^{*} - f}^{*}}{1}
= \lim_{j \in \mathcal{N} \times I} \Norm{u_{j} \star f^{*} - f^{*}}{1}
= 0.
$$
As $ \Norm{u_{j}^{*}}{1} = \Norm{u_{j}}{1} \leq 1 $ for every $ j \in \mathcal{N} \times I $, we see that $ \Seq{u_{j}^{*}}{j \in \mathcal{N} \times I} $ is a right approximate identity for $ \Quad{\Cc{G,A}}{\star}{^{*}}{\Norm{\cdot}{1}} $ that is norm-bounded by $ 1 $. By \autoref{Constructing a Two-Sided Approximate Identity from Left and Right Ones}, $ \Seq{v_{j}}{j \in \mathcal{N} \times I} $ is therefore a two-sided approximate identity for $ \Quad{\Cc{G,A}}{\star}{^{*}}{\Norm{\cdot}{1}} $ that is norm-bounded by $ 3 $.

According to \autoref{Three Bounded Left Approximate Identities for the Twisted Convolution Algebra}, $ \Seq{\Delta^{- \frac{1}{2}} u_{j}}{j \in \mathcal{N} \times I} $ is a left approximate identity for $ \Quad{\Cc{G,A}}{\star}{^{*}}{\Norm{\cdot}{1}} $, so by \Identity{Convolution Identity},
$$
  \lim_{j \in \mathcal{N} \times I} \Norm{u_{j} \star f - f}{\flat}
= \lim_{j \in \mathcal{N} \times I} \Norm{\Delta^{- \frac{1}{2}} \Br{u_{j} \star f - f}}{1}
= \lim_{j \in \mathcal{N} \times I} \Norm{\Br{\Delta^{- \frac{1}{2}} u_{j}} \star \Br{\Delta^{- \frac{1}{2}} f} - \Delta^{- \frac{1}{2}} f}{1}
= 0
$$
for every $ f \in \Cc{G,A} $. Then as $ \Norm{u_{j}}{\flat} = \Norm{\Delta^{- \frac{1}{2}} u_{j}}{1} \leq 1 $ for every $ j \in \mathcal{N} \times I $, we find that $ \Seq{u_{j}}{j \in \mathcal{N} \times I} $ is a left approximate identity for $ \Quad{\Cc{G,A}}{\star}{^{*}}{\Norm{\cdot}{\flat}} $ that is norm-bounded by $ 1 $.

According to \autoref{Three Bounded Left Approximate Identities for the Twisted Convolution Algebra}, $ \Seq{\Delta^{\frac{1}{2}} u_{j}}{j \in \mathcal{N} \times I} $ is a left approximate identity for $ \Quad{\Cc{G,A}}{\star}{^{*}}{\Norm{\cdot}{1}} $, so by \Identity{Convolution Identity},
\begin{align*}
    \lim_{j \in \mathcal{N} \times I} \Norm{f \star u_{j}^{*} - f}{\flat}
& = \lim_{j \in \mathcal{N} \times I} \Norm{\Delta^{- \frac{1}{2}} \Br{f \star u_{j}^{*} - f}}{1} \\
& = \lim_{j \in \mathcal{N} \times I}
    \Norm{\Br{\Delta^{- \frac{1}{2}} f} \star \Br{\Delta^{- \frac{1}{2}} u_{j}^{*}} - \Delta^{- \frac{1}{2}} f}{1} \\
& = \lim_{j \in \mathcal{N} \times I}
    \Norm{\SqBr{\Br{\Delta^{- \frac{1}{2}} f} \star \Br{\Delta^{- \frac{1}{2}} u_{j}^{*}} - \Delta^{- \frac{1}{2}} f}^{*}}{1} \\
& = \lim_{j \in \mathcal{N} \times I}
    \Norm{\Br{\Delta^{- \frac{1}{2}} u_{j}^{*}}^{*} \star \Br{\Delta^{- \frac{1}{2}} f}^{*} - \Br{\Delta^{- \frac{1}{2}} f}^{*}}{1} \\
& = \lim_{j \in \mathcal{N} \times I}
    \Norm{\Br{\Delta^{\frac{1}{2}} u_{j}} \star \Br{\Delta^{\frac{1}{2}} f^{*}} - \Delta^{\frac{1}{2}} f^{*}}{1} \qquad
    \Br{\text{By \Identity{Involution Identity 1}.}} \\
& = 0
\end{align*}
for every $ f \in \Cc{G,A} $. Another application of \Identity{Involution Identity 1} gives
$$
\forall j \in \mathcal{N} \times I: \quad
     \Norm{u_{j}^{*}}{\flat}
=    \Norm{\Delta^{- \frac{1}{2}} u_{j}^{*}}{1}
=    \Norm{\Br{\Delta^{- \frac{1}{2}} u_{j}^{*}}^{*}}{1}
=    \Norm{\Delta^{\frac{1}{2}} u_{j}}{1}
\leq 1,
$$
so $ \Seq{u_{j}^{*}}{j \in \mathcal{N} \times I} $ is a right approximate identity for $ \Quad{\Cc{G,A}}{\star}{^{*}}{\Norm{\cdot}{\flat}} $ that is norm-bounded by $ 1 $.

By \autoref{Constructing a Two-Sided Approximate Identity from Left and Right Ones}, $ \Seq{v_{j}}{j \in \mathcal{N} \times I} $ is therefore a two-sided approximate identity for $ \Quad{\Cc{G,A}}{\star}{^{*}}{\Norm{\cdot}{\flat}} $ that is norm-bounded by $ 3 $.
\end{proof}



\begin{Cor} \label{A Special Approximate Identity for the Reduced Twisted Crossed Product}
$ \Seq{\func{\rho}{v_{j}}}{j \in \mathcal{N} \times I} $ is a self-adjoint two-sided approximate identity for $ \RTCP{G}{A}{\alpha}{\omega} $ that is norm-bounded by $ 3 $.
\end{Cor}

\begin{proof}
This follows from \autoref{A Self-Adjoint Two-Sided Approximate Identity for the Twisted Convolution Algebra for Two Different Norms} and the fact that $ \rho $ is a $ * $-homomorphism with $ \Norm{\rho}{} \leq 1 $.
\end{proof}



\section{Relative Continuity} 

\textbf{In this section, $ \E $ is a Hilbert $ \Quad{G}{A}{\alpha}{\omega} $-module.}

\subsection{Definition} 

As mentioned earlier on, relative continuity was first defined by Ruy Exel in \cite{Exel}, in the context of a $ C^{*} $-dynamical system $ \Trip{G}{A}{\alpha} $ with abelian $ G $, as a relation $ R $ on $ A_{\si} $, and he proved that
$$
\forall a,b \in A_{\si}: \quad
\Pair{a}{b} \in R \iff \BBraKKet{a}{b} \in \RCP{G}{A}{\alpha}.
$$

In \cite{Meyer2}, Meyer defined relative continuity for the case of non-abelian $ G $ by adopting the relation $ \BBraKKet{\cdot}{\cdot} \subseteq \RCP{G}{A}{\alpha} $ as the \emph{defining condition}.


\begin{Def} \label{Relative Continuity}
A linear subspace $ \R $ of $ \E $ is called \emph{relatively continuous} if and only if
$$
\R \subseteq \Esi \qquad \text{and} \qquad
          \BBraKKet{\R}{\R}
\df       \Set{\BBra{\zeta} \circ \KKet{\eta}}{\zeta,\eta \in \R}
\subseteq \RTCP{G}{A}{\alpha}{\omega}.
$$
\end{Def}



\begin{Eg} \label{The Set of Compactly Supported Continuous A-Valued Functions on G Is a Relatively Continuous Subspace}
We have already shown in \autoref{A Very Important Corollary} that $ \KKet{\func{q}{f}} = \func{\rho}{f^{\sharp}} $ for every $ f \in \Cc{G,A} $, so we have
$$
\forall f \in \Cc{G,A}: \quad
  \BBra{\func{q}{f}}
= \KKet{\func{q}{f}}^{*}
= \func{\rho}{f^{\sharp}}^{*}
= \func{\rho}{\Br{f^{\sharp}}^{*}}.
$$
Hence,
$$
\forall f,g \in \Cc{G,A}: \quad
    \BBraKKet{\func{q}{f}}{\func{q}{g}}
=   \func{\rho}{\Br{f^{\sharp}}^{*}} \circ \func{\rho}{g^{\sharp}}
=   \func{\rho}{\Br{f^{\sharp}}^{*} \star g^{\sharp}}
\in \RTCP{G}{A}{\alpha}{\omega},
$$
which proves that $ \Im{q}{\Cc{G,A}} $ is a dense relatively continuous subspace of $ \LL{2}{G,A} $.
\end{Eg}


Having relatively continuous subspaces allows us to construct generalized fixed-point algebras. We will describe the construction later.

\subsection{Square-Integrable Completeness} 

$ \E $ could have no dense relatively continuous subspaces or it could have many of them. In order to gain finer control in the latter case, Meyer introduced a structural condition on relatively continuous subspaces, which we call \emph{square-integrable completeness}.


\begin{Def} \label{S.i.-Completeness}
We say that a relatively continuous subspace $ \R $ of $ \E $ is \emph{square-integrably complete}, or simply \emph{s.i.-complete}, if and only if the following conditions hold:
\begin{itemize}
\item
$ \R *_{\E} \Cc{G,A} \subseteq \R $.

\item
$ \R $ is complete with respect to $ \Norm{\cdot}{\E,\si} $.
\end{itemize}
\end{Def}



\begin{Prop} \label{The Compactness of G Implies the Uniqueness of Dense S.i.-Complete Relatively Continuous Subspace}
Suppose that $ G $ in $ \Quad{G}{A}{\alpha}{\omega} $ is compact. Then the only dense s.i.-complete relatively continuous subspace of $ \E $ is itself.
\end{Prop}

\begin{proof}
By \autoref{The Compactness of G Implies the Equivalence of the Hilbert-Module Norm and the S.i.-Norm}, $ \Esi = \E $. Let $ \zeta,\eta \in \E $, and define a function
$$
f_{\zeta,\eta} \df \Map{G}{A}{x}{\Del{x}{- \frac{1}{2}} \Inner{\zeta}{\gammArg{\E}{x}{\eta}}{\E}},
$$
which belongs to $ \Cc{G,A} $ as $ G $ is compact. Then for every $ \phi \in \Cc{G,A} $, \Identity{An Explicit Formula for the Integrated Form of the Covariant Representation of G and A} yields
\begin{align*}
  & ~ \FUNC{\func{\rho}{f_{\zeta,\eta}}}{\func{q}{\phi}} \\
= & ~ \func{q}{\Map{G}{A}{x}{\Int{G}{\Del{y}{\frac{1}{2}} \alphArg{x}{\func{f_{\zeta,\eta}}{y}} ~ \om{x}{y} ~ \func{\phi}{x y}}{y}}} \\
= & ~ \func{q}{
              \Map{G}{A}{x}{
                           \Int{G}{
                                  \Del{y}{\frac{1}{2}}
                                  \alphArg{x}{\Del{y}{- \frac{1}{2}} \Inner{\zeta}{\gammArg{\E}{y}{\eta}}{\E}} ~
                                  \om{x}{y} ~
                                  \func{\phi}{x y}
                                  }{y}
                           }
              } \\
= & ~ \func{q}{\Map{G}{A}{x}{\Int{G}{\alphArg{x}{\Inner{\zeta}{\gammArg{\E}{y}{\eta}}{\E}} ~ \om{x}{y} ~ \func{\phi}{x y}}{y}}} \\
= & ~ \func{q}{
              \Map{G}{A}{x}{
                           \Int{G}{
                                  \Inner{\gammArg{\E}{x}{\zeta}}{\gammArg{\E}{x}{\gammArg{\E}{y}{\eta}}}{\E} ~ \om{x}{y} ~ \func{\phi}{x y}
                                  }{y}
                           }
              } \\
= & ~ \func{q}{
              \Map{G}{A}{x}{
                           \Int{G}{
                                  \Inner{\gammArg{\E}{x}{\zeta}}{\gammArg{\E}{x y}{\eta} \bullet \om{x}{y}^{*}}{\E} ~
                                  \om{x}{y} ~
                                  \func{\phi}{x y}
                                  }{y}
                           }
              } \\
= & ~ \func{q}{
              \Map{G}{A}{x}{
                           \Int{G}{
                                  \Inner{\gammArg{\E}{x}{\zeta}}{\gammArg{\E}{x y}{\eta} \bullet \om{x}{y}^{*} \om{x}{y} ~ \func{\phi}{x y}}
                                        {\E}
                                  }{y}
                           }
              } \\
= & ~ \func{q}{\Map{G}{A}{x}{\Int{G}{\Inner{\gammArg{\E}{x}{\zeta}}{\gammArg{\E}{x y}{\eta} \bullet \func{\phi}{x y}}{\E}}{y}}} \\
= & ~ \func{q}{\Map{G}{A}{x}{\Int{G}{\Inner{\gammArg{\E}{x}{\zeta}}{\gammArg{\E}{y}{\eta} \bullet \func{\phi}{y}}{\E}}{y}}} \qquad
      \Br{\text{By the change of variables $ y \mapsto x^{-1} y $.}} \\
= & ~ \func{q}{\Map{G}{A}{x}{\Inner{\gammArg{\E}{x}{\zeta}}{\Int{G}{\gammArg{\E}{y}{\eta} \bullet \func{\phi}{y}}{y}}{\E}}} \qquad
      \Br{\text{By continuity.}} \\
= & ~ \func{q}{\Map{G}{A}{x}{\Inner{\gammArg{\E}{x}{\zeta}}{\KetArg{\eta}{\phi}}{\E}}} \\
= & ~ \func{q}{\BraArg{\zeta}{\KetArg{\eta}{\phi}}} \\
= & ~ \func{q}{\BraArg{\zeta}{\KKetArg{\eta}{\func{q}{\phi}}}} \\
= & ~ \func{\BBraKKet{\zeta}{\eta}}{\func{q}{\phi}}.
\end{align*}
As $ \Im{q}{\Cc{G,A}} $ is dense in $ \LL{2}{G,A} $, we get $ \BBraKKet{\zeta}{\eta} = \func{\rho}{f_{\zeta,\eta}} \in \RTCP{G}{A}{\alpha}{\omega} $, and as $ \zeta $ and $ \eta $ are arbitrary, we find that $ \E $ is a relatively continuous subspace of itself.

Now, if $ \S $ is any dense s.i.-complete relatively continuous subspace of $ \E $, then
\begin{align*}
    \S
& = \Cl{\S}{\E,\si} \qquad \Br{\text{As $ \S $ is s.i.-complete.}} \\
& = \Cl{\S}{\E} \qquad \Br{\text{By \autoref{The Compactness of G Implies the Equivalence of the Hilbert-Module Norm and the S.i.-Norm}.}} \\
& = \E. \qquad \Br{\text{As $ \S $ is dense in $ \E $.}}
\end{align*}
Therefore, the only dense s.i.-complete relatively continuous subspace of $ \E $ is itself.
\end{proof}



\begin{Prop} \label{A Pseudo-Approximate Identity for Relatively Continuous Subspaces}
Let $ \R $ be a relatively continuous subspace of $ \E $. If $ \Seq{v_{j}}{j \in \mathcal{N} \times I} $ is the net defined in \autoref{A Self-Adjoint Two-Sided Approximate Identity for the Twisted Convolution Algebra for Two Different Norms}, then for every $ \zeta \in \R $,
$$
  \lim_{j \in \mathcal{N} \times I} \Norm{\zeta *_{\E} v_{j} - \zeta}{\E}
= \lim_{j \in \mathcal{N} \times I} \BigNorm{\KKet{\zeta *_{\E} v_{j}} - \KKet{\zeta}}{\AdjEqPair{\LL{2}{G,A}}{\E}}
= \lim_{j \in \mathcal{N} \times I} \Norm{\zeta *_{\E} v_{j} - \zeta}{\E,\si}
= 0.
$$
\end{Prop}

\begin{proof}
Let $ \zeta \in \R $ and $ \epsilon > 0 $. Using \autoref{A Density Result for the Space of Square-Integrable Elements}, we can find $ \eta \in \Esi $ and $ f \in \Cc{G,A} $ that satisfy $ \Norm{\zeta - \eta *_{\E} f}{\E} < \dfrac{\epsilon}{9} $. Pick $ j_{0} \in \mathcal{N} \times I $ so that $ \Norm{f \star v_{j} - f}{\flat} \leq \dfrac{\epsilon}{3 \Br{\Norm{\eta}{\E} + 1}} $ for every $ j \in \Br{\mathcal{N} \times I}_{\geq j_{0}} $. Then for every such $ j $,
\begin{align*}
       \Norm{\zeta *_{\E} v_{j} - \zeta}{\E}
& \leq \Norm{\zeta *_{\E} v_{j} - \eta *_{\E} \Br{f \star v_{j}}}{\E} +
       \Norm{\eta *_{\E} \Br{f \star v_{j}} - \eta *_{\E} f}{\E} +
       \Norm{\eta *_{\E} f - \zeta}{\E} \\
& =    \Norm{\zeta *_{\E} v_{j} - \Br{\eta *_{\E} f} *_{\E} v_{j}}{\E} +
       \Norm{\eta *_{\E} \Br{f \star v_{j}} - \eta *_{\E} f}{\E} +
       \Norm{\eta *_{\E} f - \zeta}{\E} \\
& =    \Norm{\Br{\zeta - \eta *_{\E} f} *_{\E} v_{j}}{\E} +
       \Norm{\eta *_{\E} \Br{f \star v_{j} - f}}{\E} +
       \Norm{\eta *_{\E} f - \zeta}{\E} \\
& \leq \Norm{\zeta - \eta *_{\E} f}{\E} \Norm{v_{j}}{\flat} +
       \Norm{\eta}{\E} \Norm{f \star v_{j} - f}{\flat} +
       \Norm{\eta *_{\E} f - \zeta}{\E} \qquad \Br{\text{By \Inequality{Main Norm Inequality 1}.}} \\
& <    \frac{\epsilon}{9} \cdot 3 + \Norm{\eta}{\E} \cdot \frac{\epsilon}{3 \Br{\Norm{\eta}{\E} + 1}} + \frac{\epsilon}{9} \\
& <    \frac{\epsilon}{3} + \frac{\epsilon}{3} + \frac{\epsilon}{3} \\
& =    \epsilon.
\end{align*}
As $ \epsilon $ is arbitrary, we obtain $ \D \lim_{j \in \mathcal{N} \times I} ~ \Norm{\zeta *_{\E} v_{j} - \zeta}{\E} = 0 $.

Next, \autoref{A Special Approximate Identity for the Reduced Twisted Crossed Product} and the relative continuity of $ \R $ both yield
\begin{align*}
  & ~ \lim_{j \in \mathcal{N} \times I} \BigNorm{\KKet{\zeta *_{\E} v_{j}} - \KKet{\zeta}}{\AdjEqPair{\LL{2}{G,A}}{\E}} \\
= & ~ \lim_{j \in \mathcal{N} \times I} \BigNorm{\KKet{\zeta} \circ \func{\rho}{v_{j}} - \KKet{\zeta}}{\AdjEqPair{\LL{2}{G,A}}{\E}} \\
= & ~ \lim_{j \in \mathcal{N} \times I}
      \BigNorm{
              \SqBr{\KKet{\zeta} \circ \func{\rho}{v_{j}} - \KKet{\zeta}}^{*} \circ \SqBr{\KKet{\zeta} \circ \func{\rho}{v_{j}} - \KKet{\zeta}}
              }{\AdjEq{\LL{2}{G,A}}}^{\frac{1}{2}} \\
= & ~ \lim_{j \in \mathcal{N} \times I}
      \BigNorm{
              \SqBr{\func{\rho}{v_{j}}^{*} \circ \BBra{\zeta} - \BBra{\zeta}} \circ \SqBr{\KKet{\zeta} \circ \func{\rho}{v_{j}} - \KKet{\zeta}}
              }{\AdjEq{\LL{2}{G,A}}}^{\frac{1}{2}} \\
= & ~ \lim_{j \in \mathcal{N} \times I}
      \BigNorm{
              \SqBr{\func{\rho}{v_{j}^{*}} \circ \BBra{\zeta} - \BBra{\zeta}} \circ \SqBr{\KKet{\zeta} \circ \func{\rho}{v_{j}} - \KKet{\zeta}}
              }{\AdjEq{\LL{2}{G,A}}}^{\frac{1}{2}} \\
= & ~ \lim_{j \in \mathcal{N} \times I}
      \BigNorm{
              \SqBr{\func{\rho}{v_{j}} \circ \BBra{\zeta} - \BBra{\zeta}} \circ \SqBr{\KKet{\zeta} \circ \func{\rho}{v_{j}} - \KKet{\zeta}}
              }{\AdjEq{\LL{2}{G,A}}}^{\frac{1}{2}} \\
= & ~ \lim_{j \in \mathcal{N} \times I}
      \BigNorm{
              \func{\rho}{v_{j}} \circ \BBraKKet{\zeta}{\zeta} \circ \func{\rho}{v_{j}} -
              \BBraKKet{\zeta}{\zeta} \circ \func{\rho}{v_{j}} -
              \func{\rho}{v_{j}} \circ \BBraKKet{\zeta}{\zeta} +
              \BBraKKet{\zeta}{\zeta}
              }{\AdjEq{\LL{2}{G,A}}}^{\frac{1}{2}} \\
= & ~ 0.
\end{align*}

If $ \Seq{v_{j}}{j \in \mathcal{N} \times I} $ had only been a bounded self-adjoint approximate identity for $ \Quad{\Cc{G,A}}{\star}{^{*}}{\Norm{\cdot}{1}} $, then although $ \Seq{\func{\rho}{v_{j}}}{j \in \mathcal{N} \times I} $ would be a bounded two-sided approximate identity for $ \RTCP{G}{A}{\alpha}{\omega} $, the $ \Norm{\cdot}{\flat} $-boundedness of $ \Seq{v_{j}}{j \in \mathcal{N} \times I} $ would not be guaranteed to ensure $ \D \lim_{j \in \mathcal{N} \times I} ~ \Norm{\zeta *_{\E} v_{j} - \zeta}{\E} = 0 $.

Combining the foregoing arguments, we get
$$
  \lim_{j \in \mathcal{N} \times I} \Norm{\zeta *_{\E} v_{j} - \zeta}{\E,\si}
= \lim_{j \in \mathcal{N} \times I} \Norm{\zeta *_{\E} v_{j} - \zeta}{\E} +
  \lim_{j \in \mathcal{N} \times I} \BigNorm{\KKet{\zeta *_{\E} v_{j}} - \KKet{\zeta}}{\AdjEqPair{\LL{2}{G,A}}{\E}}
= 0.
$$
As $ \zeta $ is arbitrary, this concludes the proof.
\end{proof}



\begin{Prop} \label{A Density Result for S.i.-Complete Relatively Continuous Subspaces}
Let $ \R $ be a relatively continuous subspace of $ \E $. Then
\begin{alignat*}{3}
             \Cl{\R}{\E}
&  \subseteq \Cl{\R *_{\E} \Cc{G,A}}{\E}
&& \subseteq \Cl{\Span{\R *_{\E} \Cc{G,A}}}{\E}, \\
             \Cl{\KKet{\R}}{\AdjEqPair{\LL{2}{G,A}}{\E}}
&  \subseteq \Cl{\KKet{\R *_{\E} \Cc{G,A}}}{\AdjEqPair{\LL{2}{G,A}}{\E}}
&& \subseteq \Cl{\KKet{\Span{\R *_{\E} \Cc{G,A}}}}{\AdjEqPair{\LL{2}{G,A}}{\E}}, \\
             \Cl{\R}{\E,\si}
&  \subseteq \Cl{\R *_{\E} \Cc{G,A}}{\E,\si}
&& \subseteq \Cl{\Span{\R *_{\E} \Cc{G,A}}}{\E,\si},
\end{alignat*}
with equalities occurring if $ \R $ is s.i.-complete.
\end{Prop}

\begin{proof}
The second inclusions are obvious, and if $ \R $ is s.i.-complete, then $ \Span{\R *_{\E} \Cc{G,A}} \subseteq \R $. Hence, it suffices to prove the first inclusions to obtain the full result. However, we have from \autoref{A Pseudo-Approximate Identity for Relatively Continuous Subspaces} that
\begin{align*}
\forall \zeta \in \R: \quad
    \lim_{j \in \mathcal{N} \times I} \Norm{\zeta *_{\E} v_{j} - \zeta}{\E}
& = \lim_{j \in \mathcal{N} \times I} \BigNorm{\KKet{\zeta *_{\E} v_{j}} - \KKet{\zeta}}{\AdjEqPair{\LL{2}{G,A}}{\E}} \\
& = \lim_{j \in \mathcal{N} \times I} \Norm{\zeta *_{\E} v_{j} - \zeta}{\E,\si} \\
& = 0.
\end{align*}
Therefore,
$$
\R        \subseteq \Cl{\R *_{\E} \Cc{G,A}}{\E}, \qquad
\KKet{\R} \subseteq \Cl{\KKet{\R *_{\E} \Cc{G,A}}}{\AdjEqPair{\LL{2}{G,A}}{\E}}, \qquad
\R        \subseteq \Cl{\R *_{\E} \Cc{G,A}}{\E,\si},
$$
or equivalently,
\begin{align*}
\Cl{\R}{\E}                                 & \subseteq \Cl{\R *_{\E} \Cc{G,A}}{\E}, \\
\Cl{\KKet{\R}}{\AdjEqPair{\LL{2}{G,A}}{\E}} & \subseteq \Cl{\KKet{\R *_{\E} \Cc{G,A}}}{\AdjEqPair{\LL{2}{G,A}}{\E}}, \\
\Cl{\R}{\E,\si}                             & \subseteq \Cl{\R *_{\E} \Cc{G,A}}{\E,\si}.
\end{align*}
The proof is now complete.
\end{proof}



\section{Concrete Representations of Hilbert $ C^{*} $-Modules} 

\textbf{In this section, $ \E,\L $ are Hilbert $ \Quad{G}{A}{\alpha}{\omega} $-modules and $ \A $ an $ \L $-essential $ C^{*} $-subalgebra}, i.e., a $ C^{*} $-subalgebra $ \A $ of $ \AdjEq{\L} $ such that $ \Span{\Im{\A}{\L}} $ is dense in $ \L $. By the Cohen Factorization Theorem, we have, in fact, $ \Im{\A}{\L} = \L $.

Observe that $ \RTCP{G}{A}{\alpha}{\omega} $ is an $ \LL{2}{G,A} $-essential $ C^{*} $-algebra.


\begin{Def} \label{Concrete Hilbert C*-Modules}
A \emph{concrete Hilbert $ \Trip{\E}{\L}{\A} $-module} is a closed linear subspace $ \M $ of $ \AdjEqPair{\L}{\E} $, where $ \M \circ \A \subseteq \M $ and $ \M^{*} \circ \M \subseteq \A $, and we say that $ \M $ is \emph{essential} if and only if $ \Span{\Im{\M}{\L}} $ is dense in $ \E $.
\end{Def}


Any concrete Hilbert $ \Trip{\E}{\L}{\A} $-module $ \M $ can be `essentialized' by shrinking $ \E $ appropriately. Indeed, if $ \E' = \Cl{\Span{\Im{\M}{\L}}}{\E} $, then $ \E' $ is a Hilbert $ \Quad{G}{A}{\alpha}{\omega} $-submodule of $ \E $, and $ \M $ is an essential concrete Hilbert $ \Trip{\E'}{\L}{\A} $-module.

Concrete Hilbert $ C^{*} $-modules provide us with a means of concretely realizing a Hilbert module over an $ \L $-essential $ C^{*} $-subalgebra as a module of twisted-equivariant adjointable operators between Hilbert $ \Quad{G}{A}{\alpha}{\omega} $-modules.


\begin{Prop} \label{Concrete Hilbert C*-Modules Are Hilbert C*-Modules}
Let $ \M $ be a concrete Hilbert $ \Trip{\E}{\L}{\A} $-module. Then $ \M $ is a Hilbert $ \A $-module with the right $ \A $-action
$$
\forall P \in \M, ~ \forall L \in \A: \quad
P \bullet L \df P \circ L
$$
and the $ \A $-inner product
$$
\forall P,Q \in \M: \quad
\Inner{P}{Q}{\M} \df P^{*} \circ Q.
$$
The Hilbert $ \A $-module norm on $ \M $ is the restriction of $ \Norm{\cdot}{\AdjEqPair{\L}{\E}} $ to $ \M $. Furthermore,
\begin{equation} \label{Identity}
  \M
= \M \circ \A
= \M \circ \M^{*} \circ \M
\end{equation}
and
$$
  \Im{\M}{\L}
= \IM{\M \circ \M^{*}}{\E}
= \IM{\M \circ \M^{*} \circ \M}{\L}.
$$
Consequently, $ \M $ is essential if and only if $ \Span{\IM{\M \circ \M^{*}}{\E}} $ is dense in $ \E $.
\end{Prop}

\begin{proof}
We omit the easy proof that $ \bullet $ and $ \Inner{\cdot}{\cdot}{\M} $ obey the axioms of a Hilbert $ C^{*} $-module.

Observe that
$$
\forall P \in \M: \quad
    \Norm{P}{\M}
\df \sqrt{\Norm{\Inner{P}{P}{\M}}{\A}}
=   \sqrt{\Norm{P^{*} \circ P}{\A}}
=   \sqrt{\Norm{P^{*} \circ P}{\AdjEq{\L}}}
=   \Norm{P}{\AdjEqPair{\L}{\E}}.
$$
Hence, the Hilbert $ \A $-module norm on $ \M $ is the restriction of $ \Norm{\cdot}{\AdjEqPair{\L}{\E}} $ to $ \M $.

As $ \M $ is a Hilbert $ \A $-module, we have
$$
\M \circ \M^{*} \circ \M \subseteq \M \circ \A \subseteq \M.
$$
As every Hilbert $ C^{*} $-module $ \X $ has the property that each element equals $ \xi \bullet \Inner{\xi}{\xi}{\X} $ for some $ \xi \in \X $ (cf. Proposition 2.31 of \cite{Raeburn|Williams}), we also have $ \M \subseteq \M \circ \M^{*} \circ \M $. Hence,
$$
          \M
=         \M \circ \A
=         \M \circ \M^{*} \circ \M \qquad \text{and} \qquad
          \Im{\M}{\L}
=         \IM{\M \circ \M^{*} \circ \M}{\L}
\subseteq \IM{\M \circ \M^{*}}{\E}
\subseteq \Im{\M}{\L}.
$$
Therefore,
$$
  \Im{\M}{\L}
= \IM{\M \circ \M^{*}}{\E}
= \IM{\M \circ \M^{*} \circ \M}{\L},
$$
so $ \M $ is essential if and only if $ \Span{\IM{\M \circ \M^{*}}{\E}} $ is dense in $ \E $.
\end{proof}


Before proceeding further, let us state a useful result by E. Lance about unitary operators on Hilbert $ C^{*} $-modules.


\begin{Thm}[\cite{Lance}] \label{Lance's Theorem}
Let $ B $ be a $ C^{*} $-algebra and $ T: \X \to \Y $ an operator between Hilbert $ B $-modules. Then the following are equivalent:
\begin{enumerate}
\item[(i)]
$ T $ is unitary, i.e., $ T $ is adjointable, $ T^{*} \circ T = \Id_{\X} $ and $ T \circ T^{*} = \Id_{\Y} $.

\item[(ii)]
$ T $ is a $ B $-linear surjective isometry.
\end{enumerate}
\end{Thm}


Every Hilbert $ \A $-module $ \X $ can be represented as a concrete Hilbert $ \Trip{\X \Otimes{\A} \L}{\L}{\A} $-module, where $ \X \Otimes{\A} \L $ denotes the completed $ \A $-balanced tensor product of $ \X $ and $ \L $. In order to show this, we must lay some groundwork first.

Let $ \X $ be a Hilbert $ \A $-module. We can form the $ \A $-balanced algebraic tensor product $ \X \Odot{\A} \L $ because $ \L $ is a left $ \A $-module. Next, define an $ A $-valued sesquilinear form $ \Inner{\cdot}{\cdot}{\X \Odot{\A} \L} $ on $ \X \Odot{\A} \L $ by
$$
\forall \xi_{1},\xi_{2} \in \X, ~ \forall \Phi_{1},\Phi_{2} \in \L: \quad
  \Inner{\xi_{1} \Odot{\A} \Phi_{1}}{\xi_{2} \Odot{\A} \Phi_{2}}{\X \Odot{\A} \L}
= \Inner{\Phi_{1}}{\func{\Inner{\xi_{1}}{\xi_{2}}{\M}}{\Phi_{2}}}{\L}.
$$
It is a non-trivial fact that $ \X \Odot{\A} \L $ is a pre-Hilbert $ A $-module for $ \Inner{\cdot}{\cdot}{\X \Odot{\A} \L} $. If we complete $ \X \Odot{\A} \L $ with respect to the norm induced by $ \Inner{\cdot}{\cdot}{\X \Odot{\A} \L} $, we immediately get the Hilbert $ A $-module $ \X \Otimes{\A} \L $. Equipping $ \X $ with the trivial $ G $-action, $ \X \Otimes{\A} \L $ becomes a Hilbert $ \Quad{G}{A}{\alpha}{\omega} $-module.

Now, the class map
$$
\Map{\HilbMod{\A}}{\Hilb{G}{A}{\alpha}{\omega}}{\X}{\X \Otimes{\A} \L}
$$
is functorial because any $ \HilbMod{\A} $-morphism $ T: \X \to \Y $ induces a $ \Hilb{G}{A}{\alpha}{\omega} $-morphism $ T \Otimes{\A} \Id_{\L}: \X \Otimes{\A} \L \to \Y \Otimes{\A} \L $. Upon exploiting the isomorphism $ \A \Otimes{\A} \L \cong \A \cdot \L = \L $, we acquire for every Hilbert $ \A $-module $ \X $ an operator
$$
\Lambda_{\X}: \X \to \AdjEqPair{\A \Otimes{\A} \L}{\X \Otimes{\A} \L} \stackrel{\cong}{\longrightarrow} \AdjEqPair{\L}{\X \Otimes{\A} \L},
$$
which is explicitly given by
$$
\forall \xi,\xi_{1},\xi_{2} \in \X, ~ \forall \Phi \in \L: \quad
\FUNC{\func{\Lambda_{\X}}{\xi}}{\Phi} = \xi \Odot{\A} \Phi \qquad \text{and} \qquad
\FUNC{\func{\Lambda_{\X}}{\xi_{1}}^{*}}{\xi_{2} \Odot{\A} \Phi} = \func{\Inner{\xi_{1}}{\xi_{2}}{\X}}{\Phi}.
$$


\begin{Prop} \label{A Hilbert Module over an Essential C*-Subalgebra Is a Concrete Hilbert C*-Module}
Let $ \X $ be any Hilbert $ \A $-module. Then $ \Range{\Lambda_{\X}} $ is an essential concrete Hilbert $ \Trip{\X \Otimes{\A} \L}{\L}{\A} $-module. Viewing $ \Range{\Lambda_{\X}} $ as a Hilbert $ \A $-module (as per \autoref{Concrete Hilbert C*-Modules Are Hilbert C*-Modules}), we get a $ \HilbMod{\A} $-isomorphism $ \Lambda_{\X}: \X \to \Range{\Lambda_{\X}} $.

Let $ \M $ be any essential concrete Hilbert $ \Trip{\E}{\L}{\A} $-module. Viewing $ \M $ as a Hilbert $ \A $-module, there is a unitary operator $ U \in \AdjEqPair{\M \Otimes{\A} \L}{\E} $, defined on elementary tensors by $ P \Odot{\A} \Phi \mapsto \func{P}{\Phi} $ for every $ P \in \M $ and $ \Phi \in \L $, such that $ U \circ \func{\Lambda_{\M}}{P} = P $ for every $ P \in \M $.
\end{Prop}

\begin{proof}
To establish that $ \Im{\Lambda_{\X}}{\X} $ is a concrete Hilbert $ \Trip{\X \Otimes{\A} \L}{\L}{\A} $-module, we need to prove that it is a closed linear subspace of $ \AdjEqPair{\L}{\X \Otimes{\A} \L} $ and that $ \Im{\Lambda_{\X}}{\X} \circ \A \subseteq \Im{\Lambda_{\X}}{\X} $ and $ \Im{\Lambda_{\X}}{\X}^{*} \circ \Im{\Lambda_{\X}}{\X} \subseteq \A $.

Firstly, we have $ \Im{\Lambda_{\X}}{\X} \circ \A \subseteq \Im{\Lambda_{\X}}{\X} $: For every $ \xi \in \X $, $ L \in \A $ and $ \Phi \in \L $,
\begin{align*}
    \FUNC{\func{\Lambda_{\X}}{\xi} \circ L}{\Phi}
& = \FUNC{\func{\Lambda_{\X}}{\xi}}{\func{L}{\Phi}} \\
& = \xi \Odot{\A} \func{L}{\Phi} \\
& = \xi \Odot{\A} \Br{L \cdot \Phi} \\
& = \Br{\xi \bullet L} \Odot{\A} \Phi \qquad \Br{\text{As the tensor product is $ \A $-balanced.}} \\
& = \FUNC{\func{\Lambda_{\X}}{\xi \bullet L}}{\Phi}.
\end{align*}

Secondly, we have $ \Im{\Lambda_{\X}}{\X}^{*} \circ \Im{\Lambda_{\X}}{\X} \subseteq \A $: For every $ \xi_{1},\xi_{2} \in \X $ and $ \Phi \in \L $,
$$
  \FUNC{\func{\Lambda_{\X}}{\xi_{1}}^{*} \circ \func{\Lambda_{\X}}{\xi_{2}}}{\Phi}
= \FUNC{\func{\Lambda_{\X}}{\xi_{1}}^{*}}{\xi_{2} \Odot{\A} \Phi}
= \func{\Inner{\xi_{1}}{\xi_{2}}{\X}}{\Phi},
$$
so $ \func{\Lambda_{\X}}{\xi_{1}}^{*} \circ \func{\Lambda_{\X}}{\xi_{2}} = \Inner{\xi_{1}}{\xi_{2}}{\X} \in \A $.

Thirdly, for every $ \xi \in \X $, we have
$$
  \Norm{\func{\Lambda_{\X}}{\xi}}{\AdjEqPair{\L}{\X \Otimes{\A} \L}}
= \sqrt{\Norm{\func{\Lambda_{\X}}{\xi}^{*} \circ \func{\Lambda_{\X}}{\xi}}{\AdjEq{\L}}}
= \sqrt{\Norm{\Inner{\xi}{\xi}{\X}}{\AdjEq{\L}}}
= \sqrt{\Norm{\Inner{\xi}{\xi}{\X}}{\A}}
= \Norm{\xi}{\X}.
$$
Hence, $ \Lambda_{\X} $ is isometric, and the completeness of $ \X $ results in $ \Im{\Lambda_{\X}}{\X} $ being a closed linear subspace of $ \AdjEqPair{\L}{\X \Otimes{\A} \L} $. Consequently, $ \Im{\Lambda_{\X}}{\X} $ is a concrete Hilbert $ \Trip{\X \Otimes{\A} \L}{\L}{\A} $-module.

Fourthly, $ \xi \Odot{\A} \Phi = \FUNC{\func{\Lambda_{\X}}{\xi}}{\Phi} $ for every $ \xi \in \X $ and $ \Phi \in \L $, which implies that $ \Span{\IM{\Im{\Lambda_{\X}}{\X}}{\L}} $ is dense in $ \X \Otimes{\A} \L $. Therefore, $ \Im{\Lambda_{\X}}{\X} $ is essential.

Viewing $ \Im{\Lambda_{\X}}{\X} $ as a Hilbert $ \A $-module, \autoref{Lance's Theorem} now says that $ \Lambda_{\X}: \X \to \Im{\Lambda_{\X}}{\X} $ --- being an $ \A $-linear surjective isometry --- is unitary. It is thus an isomorphism of Hilbert $ \A $-modules.

Now, consider a concrete Hilbert $ \Trip{\E}{\L}{\A} $-module $ \M $. For any $ n $ elements $ P_{1},\ldots,P_{n} $ of $ \M $ and any $ n $ elements $ \Phi_{1},\ldots,\Phi_{n} $ of $ \L $, observe that
\begin{align*}
    \Norm{\sum_{k = 1}^{n} P_{k} \Odot{\A} \Phi_{k}}{\M \Otimes{\A} \L}
& = \Norm{\Inner{\sum_{k = 1}^{n} P_{k} \Odot{\A} \Phi_{k}}{\sum_{l = 1}^{n} P_{l} \Odot{\A} \Phi_{l}}{\M \Otimes{\A} \L}}{A}^{\frac{1}{2}} \\
& = \Norm{\sum_{k,l = 1}^{n} \Inner{P_{k} \Odot{\A} \Phi_{k}}{P_{l} \Odot{\A} \Phi_{l}}{\M \Otimes{\A} \L}}{A}^{\frac{1}{2}} \\
& = \Norm{\sum_{k,l = 1}^{n} \Inner{\Phi_{k}}{\func{\Inner{P_{k}}{P_{l}}{\M}}{\Phi_{l}}}{\L}}{A}^{\frac{1}{2}} \qquad
    \Br{\text{By the definition of $ \Inner{\cdot}{\cdot}{\M \Otimes{\A} \L} $.}} \\
& = \Norm{\sum_{k,l = 1}^{n} \Inner{\Phi_{k}}{\Func{P_{k}^{*} \circ P_{l}}{\Phi_{l}}}{\L}}{A}^{\frac{1}{2}} \\
& = \Norm{\sum_{k,l = 1}^{n} \Inner{\func{P_{k}}{\Phi_{k}}}{\func{P_{l}}{\Phi_{l}}}{\E}}{A}^{\frac{1}{2}} \\
& = \Norm{\Inner{\sum_{k = 1}^{n} \func{P_{k}}{\Phi_{k}}}{\sum_{l = 1}^{n} \func{P_{l}}{\Phi_{l}}}{\E}}{A}^{\frac{1}{2}} \\
& = \Norm{\sum_{k = 1}^{n} \func{P_{k}}{\Phi_{k}}}{\E}.
\end{align*}
Hence, $ U $ is well-defined and isometric. The twisted-equivariance of any $ P \in \M $ then implies that
$$
  \func{U}{\func{\Br{\tr \Otimes{\A} \gam{\L}}_{r}}{P \Odot{\A} \Phi}}
= \func{U}{P \Odot{\A} \gammArg{\L}{r}{\Phi}}
= \func{P}{\gammArg{\L}{r}{\Phi}}
= \gammArg{\E}{r}{\func{P}{\Phi}}
= \gammArg{\E}{r}{\func{U}{P \Odot{\A} \Phi}}
$$
for every $ r \in G $ and $ \Phi \in \L $, so $ U $ is twisted-equivariant also. Assuming $ \M $ to be essential, we have
$$
\Range{U} = \Cl{\Span{\Im{\M}{\L}}}{\E} = \E,
$$
which yields the surjectivity of $ U $. As $ U $ is $ A $-linear, it follows from \autoref{Lance's Theorem} that $ U $ is unitary. Finally,
$$
\forall P \in \M, ~ \forall \Phi \in \L: \quad
  \func{U}{\FUNC{\func{\Lambda_{\M}}{P}}{\Phi}}
= \func{U}{P \Odot{\A} \Phi}
= \func{P}{\Phi},
$$
whence we conclude that $ U \circ \func{\Lambda_{\M}}{P} = P $.
\end{proof}



\begin{Prop} \label{A *-Representation of Compact Operators on a Concrete Hilbert C*-Module}
Let $ \M $ be any concrete Hilbert $ \Trip{\E}{\L}{\A} $-module. The closed linear extension of
$$
\Map{\ket{\M} \bra{\M}}{\M \circ \M^{*}}{\ket{P} \bra{Q}}{P \circ Q^{*}}
$$
is then a faithful $ * $-representation of $ \Comp{\M} $ on $ \E $ (with range $ \Cl{\Span{\M \circ \M^{*}}}{\AdjEq{\E}} $) that is essential (i.e., the image of $ \Comp{\M} $ in $ \AdjEq{\E} $ is $ \E $-essential) if and only if $ \M $ is essential.

If $ \M $ is essential, then we may extend this $ * $-representation to a strictly continuous and injective unital $ * $-homomorphism $ \Theta: \Adj{\M} \to \AdjEq{\E} $ whose range is
$$
M \df \Set{S \in \AdjEq{\E}}{S \circ \M \subseteq \M ~ \textnormal{and} ~ S^{*} \circ \M \subseteq \M}.
$$
\end{Prop}

\begin{proof}
Note that $ M $ is a $ C^{*} $-subalgebra of $ \AdjEq{\E} $. Define a $ * $-homomorphism $ \Psi: M \to \Adj{\M} $ by
$$
\forall S \in M: \quad
\func{\Psi}{S} \df \Map{\M}{\M}{P}{S \circ P}.
$$
Letting $ D \df \Cl{\Span{\M \circ \M^{*}}}{\AdjEq{\E}} \subseteq \AdjEq{\E} $, we intend to prove the following three assertions:
\begin{enumerate}
\item[(a)]
$ D $ is an ideal of $ M $.

\item[(b)]
$ \Psi|_{D} $ is injective.

\item[(c)]
The range of $ \Psi|_{D} $ is $ \Comp{\M} $, and $ \Br{\Psi|_{D}}^{-1}: \Comp{\M} \to D $ is the closed linear extension of the map
$$
\Map{\ket{\M} \bra{\M}}{\M \circ \M^{*}}{\ket{P} \bra{Q}}{P \circ Q^{*}}.
$$
\end{enumerate}

To prove (a), note that $ D \circ \M \subseteq \M $ by \Identity{Identity}, and as $ D^{*} = D $, we get $ D^{*} \circ \M \subseteq \M $ too. Hence, $ D \subseteq M $. Furthermore,
\begin{gather*}
\forall S \in M: \quad
          \Br{\M \circ \M^{*}} \circ S
=         \M \circ \Br{\M^{*} \circ S}
=         \M \circ \Br{S^{*} \circ \M}^{*}
\subseteq \M \circ \M^{*} \quad \text{and} \\
          S \circ \Br{\M \circ \M^{*}}
=         \Br{S \circ \M} \circ \M^{*}
\subseteq \M \circ \M^{*}.
\end{gather*}
Therefore, $ D \circ S \subseteq D $ and $ S \circ D \subseteq D $ for every $ S \in M $, which implies that $ D $ is an ideal of $ M $.

To prove (b), suppose that $ S \in D $ satisfies $ \func{\Psi}{S} = 0_{\Adj{\M}} $. Then
$$
  S \circ \M
= \IM{\func{\Psi}{S}}{\M}
= \Im{0_{\Adj{\M}}}{\M}
= \SSet{0_{\M}}.
$$
It follows that
$$
  S \circ \Br{\M \circ \M^{*}}
= \Br{S \circ \M} \circ \M^{*}
= \SSet{0_{\AdjEq{\E}}},
$$
so $ S \circ D = \SSet{0_{\AdjEq{\E}}} $. Therefore,
$$
  S \circ S^{*} \in S \circ D^{*}
= S \circ D
= \SSet{0_{\AdjEq{\E}}},
$$
from which we obtain $ S = 0_{\AdjEq{\E}} $. This establishes the injectivity of $ \Psi|_{D} $.

To prove (c), observe for every $ P,Q,R \in \M $ that
$$
  \FUNC{\func{\Psi}{P \circ Q^{*}}}{R}
= P \circ Q^{*} \circ R
= P \circ \Inner{Q}{R}{\M}
= \Func{\ket{P} \bra{Q}}{R},
$$
which gives us $ \func{\Psi}{P \circ Q^{*}} = \ket{P} \bra{Q} $. Hence, $ \Im{\Psi}{\Span{\M \circ \M^{*}}} $ is a dense $ * $-subalgebra of $ \Comp{\M} $, and using the continuity of $ \Psi $, we get
$$
  \Cl{\Im{\Psi}{D}}{\Adj{\M}}
= \Cl{\Im{\Psi}{\Span{\M \circ \M^{*}}}}{\Adj{\M}}
= \Comp{\M}.
$$
By (b), $ \Psi|_{D} $ is an injective, thus isometric, $ * $-homomorphism from $ D $ to $ \Adj{\M} $, so its range is closed. Therefore, $ \Im{\Psi}{D} = \Cl{\Im{\Psi}{D}}{\Adj{\M}} = \Comp{\M} $, and $ \Br{\Psi|_{D}}^{-1}: \Comp{\M} \to D $ is the closed linear extension of
$$
\Map{\ket{\M} \bra{\M}}{\M \circ \M^{*}}{\ket{P} \bra{Q}}{P \circ Q^{*}}.
$$

As $ \Im{\INV{\Psi|_{D}}}{\Comp{\M}} = D \df \Cl{\Span{\M \circ \M^{*}}}{\AdjEq{\E}} $, simple closure arguments yield
$$
          \Span{\Im{\M}{\L}}
=         \Span{\IM{\M \circ \M^{*}}{\E}}
\subseteq \Span{\Im{D}{\E}}
\subseteq \Cl{\Span{\IM{\M \circ \M^{*}}{\E}}}{\E}
=         \Cl{\Span{\Im{\M}{\L}}}{\E}.
$$
Hence, $ \Cl{\Span{\IM{\M \circ \M^{*}}{\E}}}{\E} = \Cl{\Span{\Im{\M}{\L}}}{\E} $, and so by \autoref{Concrete Hilbert C*-Modules Are Hilbert C*-Modules}, $ \INV{\Psi|_{D}}: \Comp{\M} \to \AdjEq{\E} $ is an essential $ * $-representation of $ \Comp{\M} $ on $ \E $ if and only if $ \M $ is essential.

Suppose that $ \M $ is essential; it is practically $ C^{*} $-folklore that $ \INV{\Psi|_{D}}: \Comp{\M} \to \AdjEq{\E} $ can be extended to a unique, strictly continuous and \emph{injective} unital $ * $-homomorphism $ \Theta: \Adj{\M} \to \AdjEq{\E} $. For every $ \Xi \in \Adj{\M} $, we have
\begin{align*}
            \func{\Theta}{\Xi} \circ \M
& =         \func{\Theta}{\Xi} \circ D \circ \M \qquad \Br{\text{As $ D \circ \M = \M $ by \Identity{Identity}.}} \\
& =         \func{\Theta}{\Xi} \circ \Im{\Theta}{\Comp{\M}} \circ \M \qquad
            \Br{\text{As $ \Im{\Theta}{\Comp{\M}} = \Im{\INV{\Psi|_{D}}}{\Comp{\M}} = D $.}} \\
& =         \Im{\Theta}{\Xi \circ \Comp{\M}} \circ \M \\
& \subseteq \Im{\Theta}{\Comp{\M}} \circ \M \qquad \Br{\text{As $ \Comp{\M} $ is an ideal of $ \Adj{\M} $.}} \\
& =         D \circ \M \\
& =         \M,
\end{align*}
so $ \func{\Theta}{\Xi}^{*} \circ \M = \func{\Theta}{\Xi^{*}} \circ \M \subseteq \M $ as well. Therefore, $ \Range{\Theta} \subseteq M $. To show that $ \Range{\Theta} = M $, it suffices to establish that $ \Theta \circ \Psi = \Id_{M} $. Let $ S \in M $. Then every $ K \in \Comp{\M} $ and $ \zeta \in \E $, we have
\begin{align*}
  & ~ \FUNC{\Func{\Theta \circ \Psi}{S}}{\FUNC{\func{\Theta}{K}}{\zeta}} \\
= & ~ \FUNC{\func{\Theta}{\func{\Psi}{S}}}{\FUNC{\func{\Theta}{K}}{\zeta}} \\
= & ~ \FUNC{\func{\Theta}{\func{\Psi}{S}} \circ \func{\Theta}{K}}{\zeta} \\
= & ~ \FUNC{\func{\Theta}{\func{\Psi}{S} \circ K}}{\zeta} \qquad \Br{\text{As $ \Theta $ is a homomorphism.}} \\
= & ~ \FUNC{\func{\Theta}{\func{\Psi}{S} \circ \func{\Psi}{\func{\Theta}{K}}}}{\zeta} \qquad
      \Br{\text{As $ \Psi \circ \Theta|_{\Comp{\M}} = \Id_{\Comp{\M}} $.}} \\
= & ~ \FUNC{\func{\Theta}{\func{\Psi}{S \circ \func{\Theta}{K}}}}{\zeta} \qquad \Br{\text{As $ \Psi $ is a homomorphism.}} \\
= & ~ \FUNC{S \circ \func{\Theta}{K}}{\zeta} \qquad
      \Br{\text{As $ S \circ \func{\Theta}{K} \in M \circ D \subseteq D $ by (a), and $ \Theta \circ \Psi|_{D} = \Id_{D} $.}} \\
= & ~ \func{S}{\FUNC{\func{\Theta}{K}}{\zeta}}.
\end{align*}
Hence, $ \Func{\Theta \circ \Psi}{S} $ and $ S $ coincide on the dense linear subspace $ \Span{\IM{\Im{\Theta}{\Comp{\M}}}{\E}} = \Span{\Im{D}{\E}} $ of $ \E $, so $ \Func{\Theta \circ \Psi}{S} = S $ by continuity. As $ S $ is arbitrary, we conclude that $ \Theta \circ \Psi = \Id_{M} $.
\end{proof}



\section{Constructing Generalized Fixed-Point Algebras} 

\textbf{In this section, $ \E $ is a Hilbert $ \Quad{G}{A}{\alpha}{\omega} $-module and $ \R $ a relatively continuous subspace of $ \E $. Also, fix $ \L \df \LL{2}{G,A} $ and $ \A \df \RTCP{G}{A}{\alpha}{\omega} $}.

Our generalized fixed-point algebras will be constructed from $ \E $ and $ \R $. As a first step, define $ \Imp{\E}{\R} $ as the following subset of $ \AdjEqPair{\L}{\E} $:
$$
\Imp{\E}{\R} \df \Cl{\Span{\KKet{\R} \cup \Br{\KKet{\R} \circ \A}}}{\AdjEqPair{\L}{\E}}.
$$
As $ \Im{\rho}{\Cc{G,A}} $ is by construction dense in $ \A $, it is true that
$$
  \Imp{\E}{\R}
= \Cl{\Span{\KKet{\R} \cup \Br{\KKet{\R} \circ \Im{\rho}{\Cc{G,A}}}}}{\AdjEqPair{\L}{\E}}
= \Cl{\Span{\KKet{\R} \cup \Br{\KKet{\R *_{\E} \Cc{G,A}}}}}{\AdjEqPair{\L}{\E}}.
$$
Furthermore, if $ \R $ is s.i.-complete, so that $ \R *_{\E} \Cc{G,A} \subseteq \R $, then $ \Imp{\E}{\R} = \Cl{\KKet{\R}}{\AdjEqPair{\L}{\E}} $.


\begin{Prop} \label{Constructing a Concrete Hilbert C*-Module from a Relatively Continuous Subspace}
$ \Imp{\E}{\R} $ is a concrete Hilbert $ \Trip{\E}{\L}{\A} $-module. If $ \R $ is dense in $ \E $, then $ \Imp{\E}{\R} $ is essential.
\end{Prop}

\begin{proof}
By construction, $ \Imp{\E}{\R} $ is a closed subspace of $ \AdjEqPair{\L}{\E} $. Furthermore,
$$
\Imp{\E}{\R} \circ \A \subseteq \Imp{\E}{\R} \qquad \text{and} \qquad
\Imp{\E}{\R}^{*} \circ \Imp{\E}{\R} \subseteq \A.
$$
Therefore, $ \Imp{\E}{\R} $ is a concrete Hilbert $ \Trip{\E}{\L}{\A} $-module.

Now, suppose that $ \mathcal{R} $ is dense in $ \E $. \autoref{A Density Result for S.i.-Complete Relatively Continuous Subspaces} then implies that $ \mathcal{R} *_{\E} \Cc{G,A} $ is also dense in $ \E $. From the definition of $ \Imp{\E}{\R} $, we have
\begin{align*}
            \R *_{\E} \Cc{G,A}
& =         \Im{\KKet{R}}{\Im{q}{\Cc{G,A}^{\flat}}} \\
& =         \Im{\KKet{R}}{\Im{q}{\Cc{G,A}}} \qquad \Br{\text{As $ \Cc{G,A}^{\flat} = \Cc{G,A} $.}} \\
& \subseteq \IM{\Imp{\E}{\R}}{\L}. \qquad \Br{\text{As $ \KKet{\mathcal{R}} \subseteq \Imp{\E}{\R} $ and $ \Im{q}{\Cc{G,A}} \subseteq \L $.}}
\end{align*}
Therefore, $ \Span{\IM{\Imp{\E}{\R}}{\L}} $ is dense in $ \E $, so $ \Imp{\E}{\R} $ is essential.
\end{proof}


We can finally construct the generalized fixed-point algebra in our twisted setting:
\begin{itemize}
\item
By \autoref{Constructing a Concrete Hilbert C*-Module from a Relatively Continuous Subspace}, $ \Imp{\E}{\R} $ is a concrete Hilbert $ \Trip{\E}{\L}{\A} $-module, so by \autoref{Concrete Hilbert C*-Modules Are Hilbert C*-Modules}, $ \Imp{\E}{\R} $ is a Hilbert $ \A $-module, with the right $ \A $-action defined by right-composition by elements of $ \A $, and the $ \A $-inner product $ \Inner{\cdot}{\cdot}{\Imp{\E}{\R}} $ by
$$
\forall P,Q \in \Imp{\E}{\R}: \quad
\Inner{P}{Q}{\Imp{\E}{\R}} \df P^{*} \circ Q.
$$

\item
Hence, $ \Imp{\E}{\R} $ is a full Hilbert $ J $-module, with $ J \df \Cl{\Span{\Imp{\E}{\R}^{*} \circ \Imp{\E}{\R}}}{\A} $ an ideal of $ \A $.

\item
The generalized fixed-point algebra is defined as $ \FixT{\E}{\R} \df \Cl{\Span{\Imp{\E}{\R} \circ \Imp{\E}{\R}^{*}}}{\AdjEq{\E}} $. By \autoref{A *-Representation of Compact Operators on a Concrete Hilbert C*-Module}, $ \FixT{\E}{\R} $ and $ \Comp{\Imp{\E}{\R}} $ are $ * $-isomorphic.

\item
As $ \Imp{\E}{\R} $ is a $ \Pair{\Comp{\Imp{\E}{\R}}}{J} $-imprimitivity bimodule, $ \FixT{\E}{\R} $ is Morita-Rieffel equivalent to $ J $.
\end{itemize}
In the absence of twisting (i.e., $ \omega $ is trivial), our construction becomes identical to that of Meyer.


\begin{Prop} \label{The R-Map}
Let $ \M $ be any concrete Hilbert $ \Trip{\E}{\L}{\A} $-module. Define
\begin{align*}
\RM{\E}{\M}{}  & \df \Set{\zeta \in \Esi}{\KKet{\zeta} \in \M}, \\
\RM{\E}{\M}{0} & \df \Span{\Set{\func{P}{\func{q}{f}}}{P \in \M ~ \text{and} ~ f \in \Cc{G,A}}}.
\end{align*}
Then the following statements hold:
\begin{itemize}
\item
$ \RM{\E}{\M}{0} \subseteq \RM{\E}{\M}{} $.

\item
Both $ \RM{\E}{\M}{} $ and $ \RM{\E}{\M}{0} $ are relatively continuous subspaces of $ \E $, the former being s.i.-complete.

\item
Both $ \KKet{\RM{\E}{\M}{}} $ and $ \KKet{\RM{\E}{\M}{0}} $ are dense in $ \M $.

\item
$ \Imp{\E}{\RM{\E}{\M}{0}} = \Imp{\E}{\RM{\E}{\M}{}} = \M $.
\end{itemize}
\end{Prop}

\begin{proof}
Note that
$$
          \BBraKKet{\RM{\E}{\M}{}}{\RM{\E}{\M}{}}
=         \KKet{\RM{\E}{\M}{}}^{*} \circ \KKet{\RM{\E}{\M}{}}
\subseteq \M^{*} \circ \M
\subseteq \A.
$$
This implies that $ \RM{\E}{\M}{} $ is a linear subspace of $ \Esi $, so it is a relatively continuous subspace of $ \E $.

To prove that $ \RM{\E}{\M}{} $ is s.i.-complete, we must first show that it is closed under the right action $ *_{\E} $ of $ \Cc{G,A} $. Indeed,
$$
\forall \zeta \in \RM{\E}{\M}{}, ~ \forall f \in \Cc{G,A}: \quad
          \KKet{\zeta *_{\E} f}
=         \KKet{\zeta} \circ \func{\rho}{f}
\in       \M \circ \A
\subseteq \M,
$$
so $ \RM{\E}{\M}{} *_{\E} \Cc{G,A} \subseteq \RM{\E}{\M}{} $.

Next, we show that $ \RM{\E}{\M}{} $ is $ \Norm{\cdot}{\E,\si} $-complete. If $ \Seq{\zeta_{n}}{n \in \N{}} $ is a $ \Norm{\cdot}{\E,\si} $-Cauchy sequence in $ \RM{\E}{\M}{} $, then the s.i.-completeness of $ \Esi $ furnishes a $ \zeta \in \Esi $ such that $ \D \lim_{n \to \infty} \Norm{\zeta_{n} - \zeta}{\E,\si} = 0 $. In particular, $ \D \lim_{n \to \infty} \BigNorm{\KKet{\zeta_{n}} - \KKet{\zeta}}{\AdjEqPair{\L}{\E}} = 0 $. However, $ \M $ is a closed subspace of $ \AdjEqPair{\L}{\E} $, so $ \KKet{\zeta} \in \M $, which gives $ \zeta \in \RM{\E}{\M}{} $. The s.i.-completeness of $ \RM{\E}{\M}{} $ is therefore established.

As $ \Im{q}{\Cc{G,A}} \subseteq \L_{\si} $, we have $ \func{P}{\func{q}{f}} \in \Esi $ for every $ P \in \M $ and $ f \in \Cc{G,A} $, so by \Identity{Ket Identity 1},
$$
          \KKet{\func{P}{\func{q}{f}}}
=         P \circ \KKet{\func{q}{f}}
=         P \circ \func{\rho}{f^{\sharp}}
\in       \M \circ \A
\subseteq \M.
$$
Hence, $ \RM{\E}{\M}{0} \subseteq \RM{\E}{\M}{} $, making $ \RM{\E}{\M}{0} $ a relatively continuous subspace of $ \E $.

The computation in the previous paragraph also shows that
$$
\KKet{\RM{\E}{\M}{0}} = \Span{\M \circ \Im{\rho}{\Cc{G,A}}}.
$$
As $ \Im{\rho}{\Cc{G,A}} $ is dense in $ \A $, and as the right $ \A $-action on $ \M $ is non-degenerate, it follows that $ \KKet{\RM{\E}{\M}{0}} $ is dense in $ \M $. The same can then be said of $ \KKet{\RM{\E}{\M}{}} $ as $ \RM{\E}{\M}{0} \subseteq \RM{\E}{\M}{} $. Now,
\begin{align*}
            \KKet{\RM{\E}{\M}{0}}
& \subseteq \Imp{\E}{\RM{\E}{\M}{0}}
  =         \Cl{\Span{\KKet{\RM{\E}{\M}{0}} \cup \Br{\KKet{\RM{\E}{\M}{0}} \circ \A}}}{\AdjEqPair{\L}{\E}}
  \subseteq \Cl{\M}{\AdjEqPair{\L}{\E}}
  =         \M, \\
            \KKet{\RM{\E}{\M}{}}
& \subseteq \Imp{\E}{\RM{\E}{\M}{}}
  =         \Cl{\Span{\KKet{\RM{\E}{\M}{}} \cup \Br{\KKet{\RM{\E}{\M}{}} \circ \A}}}{\AdjEqPair{\L}{\E}}
  \subseteq \Cl{\M}{\AdjEqPair{\L}{\E}}
  =         \M.
\end{align*}
Taking closures therefore yields $ \Imp{\E}{\RM{\E}{\M}{0}} = \Imp{\E}{\RM{\E}{\M}{}} = \M $.
\end{proof}



\section{Categorical Results} 

\textbf{In this section, we continue to fix $ \L = \LL{2}{G,A} $ and $ \A = \RTCP{G}{A}{\alpha}{\omega} $.}

We will construct a category naturally equivalent to the category of all Hilbert $ \A $-modules, where morphisms are adjointable operators.

In the non-twisted case, this natural equivalence is already implicit in Meyer's paper \cite{Meyer2}, though it should be noted that no results on functoriality or naturality appear there. As has been our style, we will be rather pedantic about these matters and pay close attention to them.

\subsection{Continuously Square-Integrable Twisted Hilbert $ C^{*} $-Modules} 


\begin{Def} \label{Continuously Square-Integrable Twisted Hilbert C*-Modules}
A \emph{continuously square-integrable (c.s.i.)} Hilbert $ \Quad{G}{A}{\alpha}{\omega} $-module is a pair $ \Pair{\E}{\R} $, where $ \E $ is a Hilbert $ \Quad{G}{A}{\alpha}{\omega} $-module and $ \R $ a dense s.i.-complete relatively continuous subspace. Write $ \csiHilb{G}{A}{\alpha}{\omega} $ for the category of c.s.i. Hilbert $ \Quad{G}{A}{\alpha}{\omega} $-modules. If $ \Pair{\E}{\R} $ and $ \Pair{\F}{\S} $ are c.s.i. Hilbert $ \Quad{G}{A}{\alpha}{\omega} $-modules, then a morphism from $ \Pair{\E}{\R} $ to $ \Pair{\F}{\S} $ is a Hilbert $ \Quad{G}{A}{\alpha}{\omega} $-module morphism $ T: \E \to \F $ such that $ \Im{T}{\R} \subseteq \S $ and $ \Im{T^{*}}{\S} \subseteq \R $.
\end{Def}


\begin{Prop} \label{Preparation for the Classification Theorem}
Let $ \E $ be a Hilbert $ \Quad{G}{A}{\alpha}{\omega} $-module. Then the map $ \M \mapsto \RM{\E}{\M}{} $ is a bijection from the set of concrete Hilbert $ \Trip{\E}{\L}{\A} $-modules to the set of s.i.-complete relatively continuous subspaces of $ \E $. Its inverse is given by $ \R \mapsto \Imp{\E}{\R} $.

A concrete Hilbert $ \Trip{\E}{\L}{\A} $-module $ \M $ is essential if and only if $ \RM{\E}{\M}{} $ is dense in $ \E $.
\end{Prop}

\begin{proof}
\autoref{The R-Map} asserts that $ \Imp{\E}{\RM{\E}{\M}{}} = \M $, so the map $ \M \mapsto \RM{\E}{\M}{} $ is injective.

As for surjectivity, let $ \R $ be an s.i.-complete relatively continuous subspace. By \autoref{Constructing a Concrete Hilbert C*-Module from a Relatively Continuous Subspace}, $ \Imp{\E}{\R} $ is a concrete Hilbert $ \Trip{\E}{\L}{\A} $-module. Our claim is that $ \R = \RM{\E}{\Imp{\E}{\R}}{} $. Observe that $ \R \subseteq \RM{\E}{\Imp{\E}{\R}}{} $ because $ \KKet{\R} \subseteq \Imp{\E}{\R} $, so it remains to prove the reverse inclusion.

Let $ \zeta \in \RM{\E}{\Imp{\E}{\R}}{} $. As $ \Imp{\E}{\R} = \Cl{\KKet{\R}}{\AdjEqPair{\L}{\E}} $, we can find a sequence $ \Seq{\zeta_{n}}{n \in \N{}} $ in $ \R $ satisfying
$$
\lim_{n \to \infty} \BigNorm{\KKet{\zeta_{n}} - \KKet{\zeta}}{\AdjEqPair{\L}{\E}} = 0.
$$
Recalling the net $ \Seq{v_{j}}{j \in \mathcal{N} \times I} $ in \autoref{A Self-Adjoint Two-Sided Approximate Identity for the Twisted Convolution Algebra for Two Different Norms}, we have by \Inequality{Main Norm Inequality 2} that
$$
\forall j \in \mathcal{N} \times I: \quad
\lim_{n \to \infty} \Norm{\zeta_{n} *_{\E} v_{j} - \zeta *_{\E} v_{j}}{\E,\si} = 0, \quad \text{so} \quad
\zeta *_{\E} v_{j} \in \R,
$$
given that $ \R $ is s.i.-complete. Then from \autoref{A Pseudo-Approximate Identity for Relatively Continuous Subspaces}, we obtain
$$
\lim_{j \in \mathcal{N} \times I} \Norm{\zeta *_{\E} v_{j} - \zeta}{\E,\si} = 0,
$$
so $ \zeta \in \R $ by s.i.-completeness again. As $ \zeta $ is arbitrary, $ \R \subseteq \RM{\E}{\Imp{\E}{\R}}{} $.

If a concrete Hilbert $ \Trip{\E}{\L}{\A} $-module $ \M $ is essential (i.e., $ \Span{\Im{\M}{\L}} $ is dense in $ \E $), then $ \RM{\E}{\M}{0} $ is dense in $ \E $ as $ \Span{\Im{\M}{\Im{q}{\Cc{G,A}}}} $ is dense in $ \Span{\Im{\M}{\L}} $. Hence, $ \RM{\E}{\M}{} $ is dense in $ \E $. Conversely, if $ \RM{\E}{\M}{} $ is dense in $ \E $, then $ \M = \Imp{\E}{\RM{\E}{\M}{}} $ is essential by \autoref{Constructing a Concrete Hilbert C*-Module from a Relatively Continuous Subspace}.
\end{proof}


\subsection{Functoriality} 


\begin{Prop} \label{The F-Functor}
There is a functor $ \FF $ from $ \csiHilb{G}{A}{\alpha}{\omega} $ to $ \HilbMod{\A} $ defined by
$$
\func{\FF}{\E,\R} \df \Imp{\E}{\R}
$$
for every c.s.i. Hilbert $ \Quad{G}{A}{\alpha}{\omega} $-module $ \Pair{\E}{\R} $ and
$$
\func{\FF}{T} \df \Map{\Imp{\E}{\R}}{\Imp{\F}{\S}}{P}{T \circ P}
$$
for every $ \csiHilb{G}{A}{\alpha}{\omega} $-morphism $ T: \Pair{\E}{\R} \to \Pair{\F}{\S} $.
\end{Prop}

\begin{proof}
Let $ \Pair{\E}{\R} $ be a c.s.i. Hilbert $ \Quad{G}{A}{\alpha}{\omega} $-module. Then by \autoref{Concrete Hilbert C*-Modules Are Hilbert C*-Modules} and \autoref{Constructing a Concrete Hilbert C*-Module from a Relatively Continuous Subspace}, $ \Imp{\E}{\R} $ can be viewed as a Hilbert $ \A $-module.

Next, let $ T: \Pair{\E}{\R} \to \Pair{\F}{\S} $ be a $ \csiHilb{G}{A}{\alpha}{\omega} $-morphism. Observe that
$$
          T \circ \KKet{\R}
=         \KKet{\Im{T}{\R}}
\subseteq \KKet{\S}.
$$
Then as $ \Imp{\E}{\R} = \Cl{\KKet{\R}}{\AdjEqPair{\L}{\E}} $ and $ \Imp{\F}{\S} = \Cl{\KKet{S}}{\AdjEqPair{\L}{\F}} $, we obtain $ T \circ \Imp{\E}{\R} \subseteq \Imp{\F}{\S} $. Similarly,
$$
          T^{*} \circ \KKet{\S}
=         \KKet{\Im{T^{*}}{\S}}
\subseteq \KKet{\R},
$$
so $ T^{*} \circ \Imp{\F}{\S} \subseteq \Imp{\E}{\R} $. It is easily seen that $ \Map{\Imp{\F}{\S}}{\Imp{\E}{\R}}{Q}{T^{*} \circ Q} $ is the adjoint of $ \func{\FF}{T} $. Therefore, $ \func{\FF}{T} $ is a $ \HilbMod{\A} $-morphism.

Finally, as $ \FF $ obeys the Law of Composition for Functors, it is a functor.
\end{proof}



\begin{Prop} \label{The G-Functor}
There is a functor $ \GG $ from $ \HilbMod{\A} $ to $ \csiHilb{G}{A}{\alpha}{\omega} $ defined by
$$
\func{\GG}{\X} \df \Pair{\X \Otimes{\A} \L}{\RM{\X \Otimes{\A} \L}{\Range{\Lambda_{\X}}}{}}
$$
for every Hilbert $ \A $-module $ \X $ and
$$
\func{\GG}{T} \df \SSet{T \Otimes{\A} \Id_{\L}: \X \Otimes{\A} \L \to \Y \Otimes{\A} \L}
$$
for every $ \HilbMod{\A} $-morphism $ T: \X \to \Y $.
\end{Prop}

\begin{proof}
Let $ \X $ be a Hilbert $ \A $-module. By \autoref{A Hilbert Module over an Essential C*-Subalgebra Is a Concrete Hilbert C*-Module}, $ \Range{\Lambda_{\X}} $ is an essential concrete Hilbert $ \Trip{\X \Otimes{\A} \L}{\L}{\A} $-module, so \autoref{Preparation for the Classification Theorem} implies that $ \RM{\X \Otimes{\A} \L}{\Range{\Lambda_{\X}}}{} $ is a dense s.i.-complete relatively continuous subspace of $ \X \Otimes{\A} \L $. Therefore, $ \func{\GG}{\X} $ is a c.s.i. Hilbert $ \Quad{G}{A}{\alpha}{\omega} $-module.

Next, let $ T: \X \to \Y $ be a $ \HilbMod{\A} $-morphism. Note that $ T \Otimes{\A} \Id_{\L} \in \AdjEqPair{\X \Otimes{\A} \L}{\Y \Otimes{\A} \L} $, which gives us
$$
\IM{T \Otimes{\A} \Id_{\L}}{\Br{\X \Otimes{\A} \L}_{\si}} \subseteq \Br{\Y \Otimes{\A} \L}_{\si}.
$$
For each $ \zeta \in \RM{\X \Otimes{\A} \L}{\Range{\Lambda_{\X}}}{} \subseteq \Br{\X \Otimes{\A} \L}_{\si} $, we have $ \KKet{\zeta} = \func{\Lambda_{\X}}{\xi} $ for some $ \xi \in \X $, so
$$
    \KKet{\Func{T \Otimes{\A} \Id_{\L}}{\zeta}}
=   \Br{T \Otimes{\A} \Id_{\L}} \circ \KKet{\zeta}
=   \Br{T \Otimes{\A} \Id_{\L}} \circ \func{\Lambda_{\X}}{\xi}
=   \func{\Lambda_{\Y}}{\func{T}{\xi}}
\in \Range{\Lambda_{\Y}}.
$$
Hence,
$$
\IM{T \Otimes{\A} \Id_{\L}}{\RM{\X \Otimes{\A} \L}{\Range{\Lambda_{\X}}}{}} \subseteq \RM{\Y \Otimes{\A} \L}{\Range{\Lambda_{\Y}}}{}.
$$
The adjoint of $ \func{\GG}{T} $ is $ T^{*} \Otimes{\A} \Id_{\L}: \Y \Otimes{\A} \L \to \X \Otimes{\A} \L $, and similarly,
$$
\IM{T^{*} \Otimes{\A} \Id_{\L}}{\RM{\Y \Otimes{\A} \L}{\Range{\Lambda_{\Y}}}{}} \subseteq \RM{\X \Otimes{\A} \L}{\Range{\Lambda_{\X}}}{}.
$$
Therefore, $ \func{\GG}{T} $ is a $ \csiHilb{G}{A}{\alpha}{\omega} $-morphism.

Finally, as $ \GG $ obeys the Law of Composition for Functors, it is a functor.
\end{proof}


\subsection{An Equivalence of Categories} 

We have finally arrived at our main result, which is a consequence of the previous two propositions.


\begin{Cor} \label{The Classification Theorem}
There is an equivalence between $ \HilbMod{\A} $ and $ \csiHilb{G}{A}{\alpha}{\omega} $.
\end{Cor}

\begin{proof}
We must show the following:
\begin{enumerate}
\item[(i)]
There is a natural isomorphism between $ \GG \FF $ and the identity functor on $ \csiHilb{G}{A}{\alpha}{\omega} $.

\item[(ii)]
There is a natural isomorphism between $ \FF \GG $ and the identity functor on $ \HilbMod{\A} $.
\end{enumerate}

\mbox{}

\noindent \ul{Proof of (i)} \\

Let $ \Pair{\E}{\R} $ be a c.s.i. Hilbert $ \Quad{G}{A}{\alpha}{\omega} $-module. Then
$$
  \func{\GG \FF}{\E,\R}
= \func{\GG}{\Imp{\E}{\R}}
= \Pair{\Imp{\E}{\R} \Otimes{\A} \L}{\RM{\Imp{\E}{\R} \Otimes{\A} \L}{\Range{\Lambda_{\Imp{\E}{\R}}}}{}}.
$$
As $ \R $ is dense in $ \E $, \autoref{Constructing a Concrete Hilbert C*-Module from a Relatively Continuous Subspace} says that $ \Imp{\E}{\R} $ is essential. Consequently, according to \autoref{A Hilbert Module over an Essential C*-Subalgebra Is a Concrete Hilbert C*-Module}, there is a unitary operator $ U_{\Pair{\E}{\R}} \in \AdjEqPair{\Imp{\E}{\R} \Otimes{\A} \L}{\E} $ such that
$$
\func{U_{\Pair{\E}{\R}}}{\sum_{k = 1}^{n} P_{k} \otimes \Phi_{k}} = \sum_{k = 1}^{n} \func{P_{k}}{\Phi_{k}}
$$
for any $ n $ elements $ P_{1},\ldots,P_{n} \in \Imp{\E}{\R} $ and any $ n $ elements $ \Phi_{1},\ldots,\Phi_{n} \in \L $. We claim that $ U_{\Pair{\E}{\R}} $ is a $ \csiHilb{G}{A}{\alpha}{\omega} $-isomorphism, for which (because $ U_{\Pair{\E}{\R}}^{*} = \Inv{U_{\Pair{\E}{\R}}} $) it suffices to establish
$$
\Im{U_{\Pair{\E}{\R}}}{\RM{\Imp{\E}{\R} \Otimes{\A} \L}{\Range{\Lambda_{\Imp{\E}{\R}}}}{}} = \R.
$$

Observe for every $ \zeta,\eta \in \RM{\Imp{\E}{\R} \Otimes{\A} \L}{\Range{\Lambda_{\Imp{\E}{\R}}}}{} $ that
\begin{align*}
      \BBraKKet{\func{U_{\Pair{\E}{\R}}}{\zeta}}{\func{U_{\Pair{\E}{\R}}}{\eta}}
& =   \KKet{\func{U_{\Pair{\E}{\R}}}{\zeta}}^{*} \circ \KKet{\func{U_{\Pair{\E}{\R}}}{\eta}} \\
& =   \SqBr{U_{\Pair{\E}{\R}} \circ \KKet{\zeta}}^{*} \circ \SqBr{U_{\Pair{\E}{\R}} \circ \KKet{\eta}} \\
& =   \KKet{\zeta}^{*} \circ U_{\Pair{\E}{\R}}^{*} \circ U_{\Pair{\E}{\R}} \circ \KKet{\eta} \\
& =   \KKet{\zeta}^{*} \circ \KKet{\eta} \\
& =   \BBraKKet{\zeta}{\eta} \\
& \in \A, \quad \text{so}
\end{align*}
$$
\Im{U_{\Pair{\E}{\R}}}{\RM{\Imp{\E}{\R} \Otimes{\A} \L}{\Range{\Lambda_{\Imp{\E}{\R}}}}{}} \qquad \text{and} \qquad
\Im{U_{\Pair{\E}{\R}}}{\RM{\Imp{\E}{\R} \Otimes{\A} \L}{\Range{\Lambda_{\Imp{\E}{\R}}}}{0}}
$$
are relatively continuous subspaces of $ \E $. Furthermore, $ U_{\Pair{\E}{\R}} $ maps $ \Br{\Imp{\E}{\R} \Otimes{\A} \L}_{\si} $ isometrically to $ \Esi $ with respect to the norms $ \Norm{\cdot}{\Imp{\E}{\R} \Otimes{\A} \L,\si} $ and $ \Norm{\cdot}{\E,\si} $. Consequently,
$$
\Im{U_{\Pair{\E}{\R}}}{\RM{\Imp{\E}{\R} \Otimes{\A} \L}{\Range{\Lambda_{\Imp{\E}{\R}}}}{}}
$$
is complete with respect to $ \Norm{\cdot}{\E,\si} $. In addition,
\begin{align*}
          & ~ \KKet{\Im{U_{\Pair{\E}{\R}}}{\RM{\Imp{\E}{\R} \Otimes{\A} \L}{\Range{\Lambda_{\Imp{\E}{\R}}}}{}} *_{\E} \Cc{G,A}} \\
=         & ~ \KKet{\Im{U_{\Pair{\E}{\R}}}{\RM{\Imp{\E}{\R} \Otimes{\A} \L}{\Range{\Lambda_{\Imp{\E}{\R}}}}{}}} \circ \Im{\rho}{\Cc{G,A}} \\
=         & ~ U_{\Pair{\E}{\R}} \circ \KKet{\RM{\Imp{\E}{\R} \Otimes{\A} \L}{\Range{\Lambda_{\Imp{\E}{\R}}}}{}} \circ \Im{\rho}{\Cc{G,A}} \\
=         & ~ U_{\Pair{\E}{\R}} \circ
              \KKet{\RM{\Imp{\E}{\R} \Otimes{\A} \L}{\Range{\Lambda_{\Imp{\E}{\R}}}}{} *_{\Imp{\E}{\R} \Otimes{\A} \L} \Cc{G,A}} \\
\subseteq & ~ U_{\Pair{\E}{\R}} \circ \KKet{\RM{\Imp{\E}{\R} \Otimes{\A} \L}{\Range{\Lambda_{\Imp{\E}{\R}}}}{}} \qquad
              \Br{\text{By s.i.-completeness.}} \\
=         & ~ \KKet{\Im{U_{\Pair{\E}{\R}}}{\RM{\Imp{\E}{\R} \Otimes{\A} \L}{\Range{\Lambda_{\Imp{\E}{\R}}}}{}}}.
\end{align*}
By \autoref{Every Vector of a Twisted Hilbert C*-Module Gives a Unique Ket Operator}, this means that
$$
          \Im{U_{\Pair{\E}{\R}}}{\RM{\Imp{\E}{\R} \Otimes{\A} \L}{\Range{\Lambda_{\Imp{\E}{\R}}}}{}} *_{\E} \Cc{G,A}
\subseteq \Im{U_{\Pair{\E}{\R}}}{\RM{\Imp{\E}{\R} \Otimes{\A} \L}{\Range{\Lambda_{\Imp{\E}{\R}}}}{}}.
$$
Hence, $ \Im{U_{\Pair{\E}{\R}}}{\RM{\Imp{\E}{\R} \Otimes{\A} \L}{\Range{\Lambda_{\Imp{\E}{\R}}}}{}} $ is s.i.-complete, so if we can show that
$$
\Imp{\E}{\Im{U_{\Pair{\E}{\R}}}{\RM{\Imp{\E}{\R} \Otimes{\A} \L}{\Range{\Lambda_{\Imp{\E}{\R}}}}{}}} = \Imp{\E}{\R},
$$
then $ \Im{U_{\Pair{\E}{\R}}}{\RM{\Imp{\E}{\R} \Otimes{\A} \L}{\Range{\Lambda_{\Imp{\E}{\R}}}}{}} = \R $ by \autoref{Preparation for the Classification Theorem}. On one hand,
\begin{align*}
          & ~ \Imp{\E}{\Im{U_{\Pair{\E}{\R}}}{\RM{\Imp{\E}{\R} \Otimes{\A} \L}{\Range{\Lambda_{\Imp{\E}{\R}}}}{}}} \\
=         & ~ \Cl{\KKet{\Im{U_{\Pair{\E}{\R}}}{\RM{\Imp{\E}{\R} \Otimes{\A} \L}{\Range{\Lambda_{\Imp{\E}{\R}}}}{}}}}{\AdjEqPair{\L}{\E}} \\
=         & ~ \Cl{U_{\Pair{\E}{\R}} \circ \KKet{\RM{\Imp{\E}{\R} \Otimes{\A} \L}{\Range{\Lambda_{\Imp{\E}{\R}}}}{}}}{\AdjEqPair{\L}{\E}} \\
\subseteq & ~ \Cl{U_{\Pair{\E}{\R}} \circ \Range{\Lambda_{\Imp{\E}{\R}}}}{\AdjEqPair{\L}{\E}} \\
=         & ~ \Cl{\Imp{\E}{\R}}{\AdjEqPair{\L}{\E}} \qquad
              \Br{\text{By the second part of \autoref{A Hilbert Module over an Essential C*-Subalgebra Is a Concrete Hilbert C*-Module}.}} \\
=         & ~ \Imp{\E}{\R}.
\end{align*}
On the other hand,
\begin{align*}
          & ~ \Imp{\E}{\Im{U_{\Pair{\E}{\R}}}{\RM{\Imp{\E}{\R} \Otimes{\A} \L}{\Range{\Lambda_{\Imp{\E}{\R}}}}{}}} \\
\supseteq & ~ \Imp{\E}{\Im{U_{\Pair{\E}{\R}}}{\RM{\Imp{\E}{\R} \Otimes{\A} \L}{\Range{\Lambda_{\Imp{\E}{\R}}}}{0}}} \\
\supseteq & ~ \KKet{\Im{U_{\Pair{\E}{\R}}}{\RM{\Imp{\E}{\R} \Otimes{\A} \L}{\Range{\Lambda_{\Imp{\E}{\R}}}}{0}}} \\
=         & ~ \KKet{\Im{U_{\Pair{\E}{\R}}}{\Imp{\E}{\R} \Odot{\A} \Im{q}{\Cc{G,A}}}} \qquad
              \Br{\text{By the definition of $ \Lambda_{\Imp{\E}{\R}} $.}} \\
=         & ~ \KKet{\Span{\Im{\Imp{\E}{\R}}{\Im{q}{\Cc{G,A}}}}} \qquad \Br{\text{By the definition of $ U_{\Pair{\E}{\R}} $.}} \\
\supseteq & ~ \KKet{\Span{\Im{\KKet{\R}}{\Im{q}{\Cc{G,A}}}}} \\
=         & ~ \KKet{\Span{\R *_{\E} \Cc{G,A}}}. \qquad \Br{\text{As $ \Cc{G,A}^{\flat} = \Cc{G,A} $.}}
\end{align*}
However, we know from \autoref{A Density Result for S.i.-Complete Relatively Continuous Subspaces} that
$$
  \Cl{\KKet{\Span{\R *_{\E} \Cc{G,A}}}}{\AdjEqPair{\L}{\E}}
= \Cl{\KKet{\R}}{\AdjEqPair{\L}{\E}}
= \Imp{\E}{\R}, \quad \text{so}
$$
$$
\Imp{\E}{\R} \subseteq \Imp{\E}{\Im{U_{\Pair{\E}{\R}}}{\RM{\Imp{\E}{\R} \Otimes{\A} \L}{\Range{\Lambda_{\Imp{\E}{\R}}}}{}}}.
$$
Therefore,
$$
\Imp{\E}{\Im{U_{\Pair{\E}{\R}}}{\RM{\Imp{\E}{\R} \Otimes{\A} \L}{\Range{\Lambda_{\Imp{\E}{\R}}}}{}}} = \Imp{\E}{\R}
$$
as claimed, making $ U_{\Pair{\E}{\R}} $ a $ \csiHilb{G}{A}{\alpha}{\omega} $-isomorphism.

We now show that for any $ \csiHilb{G}{A}{\alpha}{\omega} $-morphism $ T: \Pair{\E}{\R} \to \Pair{\F}{\S} $, the diagram
$$
\begin{tikzcd}
\Pair{\Imp{\E}{\R} \Otimes{\A} \L}{\RM{\Imp{\E}{\R} \Otimes{\A} \L}{\Range{\Lambda_{\Imp{\E}{\R}}}}{}}
\arrow{r}{U_{\Pair{\E}{\R}}} \arrow[swap]{d}{\func{\GG \FF}{T}} &
\Pair{\E}{\R} \arrow{d}{T} \\
\Pair{\Imp{\F}{\S} \Otimes{\A} \L}{\RM{\Imp{\F}{\S} \Otimes{\A} \L}{\Range{\Lambda_{\Imp{\F}{\S}}}}{}}
\arrow{r}[swap]{U_{\Pair{\F}{\S}}} & \Pair{\F}{\S}
\end{tikzcd}
$$
commutes. Indeed,
\begin{align*}
\forall P \in \Imp{\E}{\R}, ~ \forall \Phi \in \L: \quad
    \FUNC{U_{\Pair{\F}{\S}} \circ \func{\GG \FF}{T}}{P \Odot{\A} \Phi}
& = \func{U_{\Pair{\F}{\S}}}{\FUNC{\func{\GG \FF}{T}}{P \Odot{\A} \Phi}} \\
& = \func{U_{\Pair{\F}{\S}}}{\FUNC{\func{\FF}{T} \Otimes{\A} \Id_{\L}}{P \Odot{\A} \Phi}} \\
& = \func{U_{\Pair{\F}{\S}}}{\FUNC{\func{\FF}{T}}{P} \Odot{\A} \Phi} \\
& = \func{U_{\Pair{\F}{\S}}}{\Br{T \circ P} \Odot{\A} \Phi} \\
& = \Func{T \circ P}{\Phi} \\
& = \func{T}{\func{P}{\Phi}} \\
& = \func{T}{\func{U_{\Pair{\E}{\R}}}{P \Odot{\A} \Phi}} \\
& = \Func{T \circ U_{\Pair{\E}{\R}}}{P \Odot{\A} \Phi},
\end{align*}
so by continuity and the denseness of $ \Imp{\E}{\R} \Odot{\A} \L $ in $ \Imp{\E}{\R} \Otimes{\A} \L $, we obtain
$$
U_{\Pair{\F}{\S}} \circ \func{\GG \FF}{T} = T \circ U_{\Pair{\E}{\R}}.
$$

\noindent \ul{Proof of (ii)} \\

For every Hilbert $ \A $-module $ \X $, we have
$$
  \func{\FF \GG}{\X}
= \func{\FF}{\Pair{\X \Otimes{\A}{\L}}{\RM{\X \Otimes{\A} \L}{\Range{\Lambda_{\X}}}{}}}
= \Imp{\X \Otimes{\A}{\L}}{\RM{\X \Otimes{\A} \L}{\Range{\Lambda_{\X}}}{}}
= \Range{\Lambda_{\X}},
$$
and $ \Lambda_{\X}: \X \to \Range{\Lambda_{\X}} $ is, by \autoref{A Hilbert Module over an Essential C*-Subalgebra Is a Concrete Hilbert C*-Module}, a $ \HilbMod{\A} $-isomorphism.

We must show that for any $ \HilbMod{\A} $-morphism $ T: \X \to \Y $, the diagram
$$
\begin{tikzcd}
\X \arrow{r}{T}       \arrow[swap]{d}{\Lambda_{\X}} & \Y \arrow{d}{\Lambda_{\Y}} \\
\Im{\Lambda_{\X}}{\X} \arrow{r}[swap]{\func{\FF \GG}{T}}  & \Im{\Lambda_{\Y}}{\Y}
\end{tikzcd}
$$
commutes. Indeed,
\begin{align*}
\forall \xi \in \X, ~ \forall \Phi \in \L: \quad
    \FUNC{\FUNC{\func{\FF \GG}{T} \circ \Lambda_{\X}}{\xi}}{\Phi}
& = \FUNC{\FUNC{\func{\FF \GG}{T}}{\func{\Lambda_{\X}}{\xi}}}{\Phi} \\
& = \FUNC{\FUNC{\func{\FF}{T \Otimes{} \Id_{\L}}}{\func{\Lambda_{\X}}{\xi}}}{\Phi} \\
& = \FUNC{\Br{T \Otimes{\A} \Id_{\L}} \circ \func{\Lambda_{\X}}{\xi}}{\Phi} \\
& = \Func{T \Otimes{\A} \Id_{\L}}{\FUNC{\func{\Lambda_{\X}}{\xi}}{\Phi}} \\
& = \Func{T \Otimes{\A} \Id_{\L}}{\xi \Odot{\A} \Phi} \\
& = \func{T}{\xi} \Odot{\A} \Phi \\
& = \FUNC{\func{\Lambda_{\Y}}{\func{T}{\xi}}}{\Phi} \\
& = \FUNC{\Func{\Lambda_{\Y} \circ T}{\xi}}{\Phi}, \quad \text{so} \\
    \func{\FF \GG}{T} \circ \Lambda_{\X}
& = \Lambda_{\Y} \circ T.
\end{align*}
The proof is finally complete.
\end{proof}



\begin{Eg} \label{Hilbert Modules over Non-Commutative Tori}
Consider \autoref{Twisted C*-Dynamical Systems Associated with Higher-Dimensional Non-Commutative Tori}. The $ d $-dimensional non-commutative torus $ A_{\Theta} $ is then defined as the full twisted crossed product $ \FTCP{\Z{d}}{\C{}}{\tr}{\omega_{\Theta}} $. As $ \Z{d} $ is an amenable discrete group, we have
$$
\FTCP{\Z{d}}{\C{}}{\tr}{\omega_{\Theta}} \cong \RTCP{\Z{d}}{\C{}}{\tr}{\omega_{\Theta}}
$$
by a 1968 result of Zeller-Meier \cite{Zeller-Meier}, so $ \HilbMod{A_{\Theta}} $ and $ \csiHilb{\Z{d}}{\C{}}{\tr}{\omega_{\Theta}} $ are equivalent. Therefore, every Hilbert $ A_{\Theta} $-module can be fully constructed from a Hilbert space endowed with a twisted $ \Z{d} $-action and a dense s.i.-complete relatively continuous subspace.
\end{Eg}



\section{Further Results on S.i.-Completeness} 

\textbf{In this section, $ \E $ is a Hilbert $ \Quad{G}{A}{\alpha}{\omega} $-module and $ \R $ a relatively continuous subspace of $ \E $. Also, continue to fix $ \L \df \LL{2}{G,A} $ and $ \A \df \RTCP{G}{A}{\alpha}{\omega} $.}


\begin{Def} \label{The S.i.-Completion of a Relatively Continuous Subspace}
Define the \emph{s.i.-completion} of $ \R $, denoted by $ \Cl{\R}{\si} $, as $ \Cl{\Span{\R \cup \Br{\R *_{\E} \Cc{G,A}}}}{\E,\si} $. Equivalently, it is the smallest linear subspace of $ \Esi $ containing $ \R $ that is both $ \Norm{\cdot}{\E,\si} $-closed and invariant under the right action $ *_{\E} $ of $ \Cc{G,A} $.
\end{Def}


It is not immediately clear from \autoref{The S.i.-Completion of a Relatively Continuous Subspace} that $ \Cl{\R}{\si} $ is a relatively continuous subspace of $ \E $. Our next result shows that this is indeed the case and even gives an explicit formula for it.


\begin{Thm} \label{The S.i.-Completion of a Relatively Continuous Subspace Is Still Relatively Continuous}
$ \Cl{\R}{\si} $ is relatively continuous subspace of $ \E $ and equals $ \RM{\E}{\Imp{\E}{\R}}{} $.
\end{Thm}

\begin{proof}
Let $ \zeta,\eta \in \Cl{\R}{\si} $. Then there are sequences $ \Seq{\zeta_{n}}{n \in \N{}} $ and $ \Seq{\eta_{n}}{n \in \N{}} $ in $ \Span{\R \cup \Br{\R *_{\E} \Cc{G,A}}} $ such that
$$
  \lim_{n \to \infty} \Norm{\zeta_{n} - \zeta}{\E,\si}
= \lim_{n \to \infty} \Norm{\eta_{n} - \eta}{\E,\si}
= 0.
$$
In particular,
$$
  \lim_{n \to \infty} \BigNorm{\KKet{\zeta_{n}} - \KKet{\zeta}}{\AdjEqPair{\LL{2}{G,A}}{\E}}
= \lim_{n \to \infty} \BigNorm{\KKet{\eta_{n}} - \KKet{\eta}}{\AdjEqPair{\LL{2}{G,A}}{\E}}
= 0.
$$
Hence, $ \D \BBraKKet{\zeta}{\eta} = \lim_{n \to \infty} \BBraKKet{\zeta_{n}}{\eta_{n}} \in \RTCP{G}{A}{\alpha}{\omega} $, and as $ \zeta $ and $ \eta $ are arbitrary, we conclude that $ \Cl{\R}{\si} $ is a relatively continuous subspace of $ \E $.

By \autoref{Preparation for the Classification Theorem}, $ \RM{\E}{\Imp{\E}{\R}}{} $ is an s.i.-complete relatively continuous subspace containing $ \R $. Let $ \S $ be another subspace of $ \E $ with the same properties. Then $ \S = \RM{\E}{\Imp{\E}{\S}}{} $, and as $ \R \subseteq \S $, we have $ \Imp{\E}{\R} \subseteq \Imp{\E}{\S} $, which yields
$$
          \RM{\E}{\Imp{\E}{\R}}{}
\subseteq \RM{\E}{\Imp{\E}{\S}}{}
=         \S.
$$
In particular, $ \RM{\E}{\Imp{\E}{\R}}{} \subseteq \Cl{\R}{\si} $. By definition, $ \Cl{\R}{\si} \subseteq \RM{\E}{\Imp{\E}{\R}}{} $, so $ \Cl{\R}{\si} = \RM{\E}{\Imp{\E}{\R}}{} $.
\end{proof}



\begin{Thm} \label{An S.i.-Complete Relatively Continuous Subspace Is Invariant Under the Right C*-Algebra Action and the Twisted Group Action}
Suppose that $ \R $ is s.i.-complete. Then $ \R $ is invariant under both the right $ A $-action and the twisted $ G $-action on $ \E $. Furthermore, $ \Pair{\R}{\Norm{\cdot}{\E,\si}} $ is an essential right $ A $-module, i.e., $ \R \bullet A = \R $.
\end{Thm}

\begin{proof}
By \autoref{Preparation for the Classification Theorem}, there is a concrete Hilbert $ \Trip{\E}{\L}{\A} $-module $ \M $ satisfying $ \R = \RM{\E}{\M}{} $. Let $ \zeta \in \R $. As $ \KKet{\zeta} \circ L \in \M \circ \A \subseteq \M $ for every $ L \in \A $, we get the operator from $ \A $ to $ \M $ below:
$$
\Map{\A}{\M}{L}{\KKet{\zeta} \circ L}.
$$
We contend that this extends to an operator from $ \Mult{\A} $ to $ \M $ that is continuous with respect to the strict topology on $ \Mult{\A} $ and the operator-norm topology on $ \M $, where we are viewing $ \Mult{\A} $ as the idealizer of $ \A $ in $ \Adj{\L} $ (this is justified as the inclusion $ \A \hookrightarrow \Adj{\L} $ is non-degenerate):
$$
\Mult{\A} = \Set{T \in \Adj{\L}}{T \circ \A \subseteq \A ~ \text{and} ~ \A \circ T \subseteq \A}.
$$
By Cohen's Factorization Theorem, $ \M \circ \A = \M $, so $ \KKet{\zeta} = P \circ L_{0} $ for some $ P \in \M $ and $ L_{0} \in \A $. Hence, $ \KKet{\zeta} \circ T = P \circ L_{0} \circ T \in \M \circ \A \subseteq \M $ for every $ T \in \Mult{\A} $, which gives us an operator
$$
\Map{\Mult{\A}}{\M}{T}{\KKet{\zeta} \circ T}.
$$
This operator clearly extends the one earlier. To show that it has the stated continuity conditions, observe that if $ \Seq{T_{n}}{n \in \N{}} $ is any sequence in $ \Mult{\A} $ that converges strictly to some $ T \in \Mult{A} $, then $ \D \lim_{n \to \infty} L_{0} \circ T_{n} = L_{0} \circ T $ in $ \A $, so $ \D \lim_{n \to \infty} \KKet{\zeta} \circ T_{n} = \KKet{\zeta} \circ T $ in $ \M $.

Now, let $ a \in A $ and $ r \in G $. From \Identity{Ket Identity 2} and \Identity{Ket Identity 3}, we have
$$
\KKet{\zeta \bullet a}        = \KKet{\zeta} \circ \func{\pi}{a} \qquad \text{and} \qquad
\KKet{\gammArg{\E}{r}{\zeta}} = \Del{r}{- \frac{1}{2}} \SqBr{\KKet{\zeta} \circ \func{\lambda}{r}^{*}}.
$$
Some rather straightforward calculations reveal that $ \func{\pi}{a},\func{\lambda}{r} \in \Mult{\A} $, so $ \KKet{\zeta \bullet a},\KKet{\gammArg{\E}{r}{\zeta}} \in \M $. Therefore, $ \zeta \bullet a,\gammArg{\E}{r}{\zeta} \in \RM{\E}{\M}{} $, and as $ \zeta $, $ a $ and $ r $ are arbitrary, we find that $ \R $ is invariant under both the right $ A $-action and the twisted $ G $-action on $ \E $.

Finally, let $ \Seq{e_{i}}{i \in I} $ be an approximate identity for $ A $. Then
\begin{align*}
\forall \zeta \in \R: \quad
    \lim_{i \in I} \Norm{\zeta \bullet e_{i} - \zeta}{\E,\si}
& = \lim_{i \in I} \Norm{\zeta \bullet e_{i} - \zeta}{\E} +
    \lim_{i \in I} \BigNorm{\KKet{\zeta \bullet e_{i}} - \KKet{\zeta}}{\AdjEqPair{\L}{\E}} \\
& = \lim_{i \in I} \Norm{\zeta \bullet e_{i} - \zeta}{\E} +
    \lim_{i \in I} \BigNorm{\KKet{\zeta} \circ \func{\pi}{e_{i}} - \KKet{\zeta}}{\AdjEqPair{\L}{\E}} \\
& = 0. \qquad \Br{\text{As $ \func{\pi}{e_{i}} $ converges strictly to $ \Id_{\L} $.}}
\end{align*}
By Cohen's Factorization Theorem once more, we obtain $ \R \bullet A = \R $.
\end{proof}



\begin{Prop} \label{The Final Proposition}
Suppose that $ \Pair{\E}{\R} $ is a c.s.i. Hilbert $ \Quad{G}{A}{\alpha}{\omega} $-module, and let $ \M \df \Imp{\E}{\R} $. If $ \Theta: \Adj{\M} \to \AdjEq{\E} $ is as defined in the statement of \autoref{A *-Representation of Compact Operators on a Concrete Hilbert C*-Module}, then
$$
\Range{\Theta} = \text{Set of all $ \csiHilb{G}{A}{\alpha}{\omega} $-endomorphisms on $ \Pair{\E}{\R} $}.
$$
\end{Prop}

\begin{proof}
Let $ M $ be as defined in the statement of \autoref{A *-Representation of Compact Operators on a Concrete Hilbert C*-Module}. \\

\noindent \textbf{Claim 1:} Every $ \csiHilb{G}{A}{\alpha}{\omega} $-endomorphism on $ \Pair{\E}{\R} $ is an element of $ M $.

\begin{proof}[Proof of Claim 1]
Let $ T $ be a $ \csiHilb{G}{A}{\alpha}{\omega} $-endomorphism on $ \Pair{\E}{\R} $. Then
$$
\Im{T}{\R}     \subseteq \R \qquad \text{and} \qquad
\Im{T^{*}}{\R} \subseteq \R.
$$
By \Identity{Ket Identity 1}, $ T \circ \KKet{\R} = \KKet{\Im{T}{\R}} \subseteq \KKet{\R} $, and as $ \M = \Cl{\KKet{\R}}{\AdjEqPair{\L}{\E}} $ by the s.i.-completeness of $ \R $, we have
$$
          T \circ \M
=         T \circ \Cl{\KKet{\R}}{\AdjEqPair{\L}{\E}}
\subseteq \Cl{T \circ \KKet{\R}}{\AdjEqPair{\L}{\E}}
\subseteq \Cl{\KKet{\R}}{\AdjEqPair{\L}{\E}}
=         \M.
$$
Similarly, $ T^{*} \circ \M \subseteq \M $. Therefore, $ T \in M $ by the definition of $ M $.
\end{proof}

\noindent \textbf{Claim 2:} Every element of $ M $ is a $ \csiHilb{G}{A}{\alpha}{\omega} $-endomorphism on $ \Pair{\E}{\R} $.

\begin{proof}[Proof of Claim 2]
The proof of this is more complex because it involves a tight interplay between the norms $ \Norm{\cdot}{\E} $ and $ \Norm{\cdot}{\E,\si} $. We first show that $ \Im{\M}{\Im{q}{\Cc{G,A}}} \subseteq \R $. Let $ T \in \M $ and $ \phi \in \Cc{G,A} $. As $ \KKet{\R} $ is dense in $ \M $, there is a sequence $ \Seq{\zeta_{n}}{n \in \N{}} $ in $ \R $ such that
$$
\lim_{n \to \infty} \BigNorm{\KKet{\zeta_{n}} - T}{\AdjEqPair{\L}{\E}} = 0.
$$
Then
$$
\lim_{n \to \infty} \BigNorm{\KKetArg{\zeta_{n}}{\func{q}{\phi}} - \func{T}{\func{q}{\phi}}}{\E} = 0.
$$
Now,
$$
\forall n \in \N{}: \quad
          \KKetArg{\zeta_{n}}{\func{q}{\phi}}
=         \zeta_{n} *_{\E} \phi^{\sharp}
\in       \R *_{\E} \Cc{G,A}
\subseteq \R
\subseteq \Esi, \quad \text{so}
$$
\begin{align*}
\forall m,n \in \N{}: \qquad
     & ~ \BigNorm{\KKetArg{\zeta_{m}}{\func{q}{\phi}} - \KKetArg{\zeta_{n}}{\func{q}{\phi}}}{\E,\si} \\
=    & ~ \BigNorm{\Func{\KKet{\zeta_{m}} - \KKet{\zeta_{n}}}{\func{q}{\phi}}}{\E,\si} \\
=    & ~ \BigNorm{\Func{\KKet{\zeta_{m}} - \KKet{\zeta_{n}}}{\func{q}{\phi}}}{\E} +
         \BigNorm{\KKet{\Func{\KKet{\zeta_{m}} - \KKet{\zeta_{n}}}{\func{q}{\phi}}}}{\AdjEqPair{\L}{\E}} \\
=    & ~ \BigNorm{\Func{\KKet{\zeta_{m}} - \KKet{\zeta_{n}}}{\func{q}{\phi}}}{\E} +
         \BigNorm{\Br{\KKet{\zeta_{m}} - \KKet{\zeta_{n}}} \circ \KKet{\func{q}{\phi}}}{\AdjEqPair{\L}{\E}} \\
\leq & ~ \BigNorm{\KKet{\zeta_{m}} - \KKet{\zeta_{n}}}{\AdjEqPair{\L}{\E}} \Norm{\func{q}{\phi}}{\L} +
         \BigNorm{\KKet{\zeta_{m}} - \KKet{\zeta_{n}}}{\AdjEqPair{\L}{\E}} \BigNorm{\KKet{\func{q}{\phi}}}{\AdjEq{\L}} \\
=    & ~ \BigNorm{\KKet{\zeta_{m}} - \KKet{\zeta_{n}}}{\AdjEqPair{\L}{\E}}
         \Br{\Norm{\func{q}{\phi}}{\L} + \BigNorm{\KKet{\func{q}{\phi}}}{\AdjEq{\L}}} \\
=    & ~ \BigNorm{\KKet{\zeta_{m}} - \KKet{\zeta_{n}}}{\AdjEqPair{\L}{\E}} \Norm{\func{q}{\phi}}{\L,\si}.
\end{align*}
As $ \Seq{\KKet{\zeta_{n}}}{n \in \N{}} $ is Cauchy in $ \M $, it follows that $ \Seq{\KKetArg{\zeta_{n}}{\func{q}{\phi}}}{n \in \N{}} $ is $ \Norm{\cdot}{\E,\si} $-Cauchy in $ \R $. However, $ \R $ is $ \Norm{\cdot}{\E,\si} $-complete, so there exists an $ \eta \in \R $ such that
$$
\lim_{n \to \infty} \BigNorm{\KKetArg{\zeta_{n}}{\func{q}{\phi}} - \eta}{\E,\si} = 0.
$$
Convergence with respect to $ \Norm{\cdot}{\E,\si} $ implies the same with respect to $ \Norm{\cdot}{\E} $, which means that
$$
\lim_{n \to \infty} \BigNorm{\KKetArg{\zeta_{n}}{\func{q}{\phi}} - \eta}{\E} = 0.
$$
Therefore, $ \func{T}{\func{q}{\phi}} = \eta \in \R $, and consequently, $ \Im{\M}{\Im{q}{\Cc{G,A}}} \subseteq \R $ as $ T $ and $ \phi $ are arbitrary.

By our arguments thus far, we have
$$
\forall S \in M: \quad
          \Im{S}{\Im{\M}{\Im{q}{\Cc{G,A}}}}
=         \Im{\Br{S \circ \M}}{\Im{q}{\Cc{G,A}}}
\subseteq \Im{\M}{\Im{q}{\Cc{G,A}}}
\subseteq \R.
$$
Our next goal is to show that $ \Im{\M}{\Im{q}{\Cc{G,A}}} $ is $ \Norm{\cdot}{\E,\si} $-dense in $ \R $.

Indeed, by \autoref{A Density Result for S.i.-Complete Relatively Continuous Subspaces}, $ \R *_{\E} \Cc{G,A} = \Im{\KKet{\R}}{\Im{q}{\Cc{G,A}}} $ is $ \Norm{\cdot}{\E,\si} $-dense in $ \R $, and as
$$
          \Im{\KKet{\R}}{\Im{q}{\Cc{G,A}}}
\subseteq \Im{\M}{\Im{q}{\Cc{G,A}}}
\subseteq \R,
$$
we find that $ \Im{\M}{\Im{q}{\Cc{G,A}}} $ is $ \Norm{\cdot}{\E,\si} $-dense in $ \R $.

Let $ S \in M $. We wish to prove that $ \Im{S}{\R} \subseteq \R $. Toward this end, let $ \zeta \in \R $, and pick a sequence $ \Seq{\zeta_{n}}{n \in \N{}} $ in $ \Im{\M}{\Im{q}{\Cc{G,a}}} $ where
$$
\lim_{n \to \infty} \Norm{\zeta_{n} - \zeta}{\E,\si} = 0.
$$
As mentioned earlier, $ \Norm{\cdot}{\E,\si} $-convergence implies $ \Norm{\cdot}{\E} $-convergence, so by the continuity of $ S $,
$$
\lim_{n \to \infty} \Norm{\func{S}{\zeta_{n}} - \func{S}{\zeta}}{\E} = 0.
$$
Furthermore,
\begin{align*}
\forall m,n \in \N{}: \quad
       \Norm{\func{S}{\zeta_{m}} - \func{S}{\zeta_{n}}}{\E,\si}
& =    \Norm{\func{S}{\zeta_{m} - \zeta_{n}}}{\E,\si} \\
& =    \Norm{\func{S}{\zeta_{m} - \zeta_{n}}}{\E} + \BigNorm{\KKet{\func{S}{\zeta_{m} - \zeta_{n}}}}{\AdjEqPair{\L}{\E}} \\
& =    \Norm{\func{S}{\zeta_{m} - \zeta_{n}}}{\E} + \BigNorm{S \circ \KKet{\zeta_{m} - \zeta_{n}}}{\AdjEqPair{\L}{\E}} \\
& \leq \Norm{S}{\AdjEq{\E}} \Norm{\zeta_{m} - \zeta_{n}}{\E} +
       \Norm{S}{\AdjEq{\E}} \BigNorm{\KKet{\zeta_{m} - \zeta_{n}}}{\AdjEqPair{\L}{\E}} \\
& =    \Norm{S}{\AdjEq{\E}} \Br{\Norm{\zeta_{m} - \zeta_{n}}{\E} + \BigNorm{\KKet{\zeta_{m} - \zeta_{n}}}{\AdjEqPair{\L}{\E}}} \\
& =    \Norm{S}{\AdjEq{\E}} \Norm{\zeta_{m} - \zeta_{n}}{\E,\si}.
\end{align*}
As $ \Seq{\zeta_{n}}{n \in \N{}} $ is $ \Norm{\cdot}{\E,\si} $-Cauchy in $ \R $, it follows that $ \Seq{\func{S}{\zeta_{n}}}{n \in \N{}} $ is $ \Norm{\cdot}{\E,\si} $-Cauchy in $ \R $ also. Thanks to the $ \Norm{\cdot}{\E,\si} $-completeness of $ \R $, there exists an $ \eta \in \R $ satisfying
$$
\lim_{n \to \infty} \Norm{\func{S}{\zeta_{n}} - \eta}{\E,\si} = 0.
$$
By now, it should be clear that this yields
$$
\lim_{n \to \infty} \Norm{\func{S}{\zeta_{n}} - \eta}{\E} = 0.
$$
Hence, $ \func{S}{\zeta} = \eta \in \R $, which shows that $ \Im{S}{\R} \subseteq \R $. The proof that $ \Im{S^{*}}{\R} \subseteq \R $ is similar. Therefore, $ S $ is a $ \csiHilb{G}{A}{\alpha}{\omega} $-endomorphism on $ \Pair{\E}{\R} $, and as $ S $ is arbitrary, the claim is settled.
\end{proof}

The range of $ \Theta $ is indeed the set of $ \csiHilb{G}{A}{\alpha}{\omega} $-endomorphisms on $ \Pair{\E}{\R} $.
\end{proof}



\section{Limitations and Concluding Remarks}

A major restrictive assumption that we have made in this thesis is the continuity of our maps. The maps $ \alpha $ and $ \omega $ present in $ \Quad{G}{A}{\alpha}{\omega} $ are strongly continuous and strictly continuous respectively, while the twisted action $ \gam{\E} $ on a Hilbert $ \Quad{G}{A}{\alpha}{\omega} $-module $ \E $ is strongly continuous. However, twisted $ C^{*} $-dynamical systems, when studied in full generality, are only assumed to be measurable, so one might ask: Why not work with measurable twisted $ C^{*} $-dynamical systems in the first place? The answer to this question is that in the absence of continuity, difficulties arise in trying to prove results such as \autoref{A Self-Adjoint Two-Sided Approximate Identity for the Twisted Convolution Algebra for Two Different Norms} and \autoref{Preparation for the Classification Theorem}. Left approximate identities for $ \RTCP{G}{A}{\alpha}{\omega} $, when $ \Quad{G}{A}{\alpha}{\omega} $ is measurable, definitely exist (see \cite{Packer|Raeburn}), but in general, they do not assume the nice form that we have used. As mentioned in \cite{Busby|Smith}, the usual tensor product of an approximate delta for $ \LL{1}{G} $ with an approximate identity for $ A $ does not always work. Therefore, the techniques employed here would have to be completely revamped to handle the measurable case, not to mention the special attention that has to be paid to basic measure-theoretical issues.

At the time of writing, it is not known how to make $ A $ a Hilbert $ \Quad{G}{A}{\alpha}{\omega} $-module. I consider this to be the most important problem. If we try to set $ \gam{A} \df \alpha $, just as in a $ C^{*} $-dynamical system, then we obtain an inconsistency because
$$
\forall r,s \in G, ~ \forall a \in A: \quad
\gammArg{A}{r}{\gammArg{A}{s}{a}} = \gammArg{A}{r s}{a} ~ \om{r}{s}^{*} \qquad \text{but} \qquad
\alphArg{r}{\alphArg{s}{a}}       = \om{r}{s} ~ \alphArg{r s}{a} ~ \om{r}{s}^{*}.
$$
The obstruction is caused by an extra $ \om{r}{s} $ (or a lack thereof). If we can overcome this problem, then it is possible to use our results to give a Rieffel-type definition of properness (see \autoref{Rieffel-Properness}) for a twisted $ C^{*} $-dynamical system. One might suggest that (4) of \autoref{Twisted Hilbert C*-Modules} be modified to read
$$
\forall r,s \in G: \quad
\gamm{\E}{r} \circ \gamm{\E}{s} = \Ad{\om{r}{s}} \circ \gamm{\E}{r s},
$$
but this assumes the existence of a left $ A $-action on $ \E $, which we do not have. This has been proposed by E. B\'edos and R. Conti in \cite{Bedos|Conti}, but their definition does not lead to the nice property in \autoref{Morphisms of Twisted Hilbert C*-Modules Are Closed Under the Operator-Adjoint} that morphisms are closed under operator-adjoints. In any case, these authors were not presenting a categorical viewpoint in their work.


\end{document}